\documentclass{article}

\usepackage[francais]{babel}
\usepackage[utf8]{inputenc}
\usepackage{amsmath,amsfonts,bbold,stmaryrd,txfonts,graphicx}
\usepackage{fullpage}

\newtheorem{prout}{Théorème} [subsection]
\newtheorem{prop}[prout]{Proposition}
\newtheorem{lemme}[prout]{Lemme}
\newtheorem{defin}[prout]{Définition}
\newtheorem{defins}[prout]{Définitions}
\newtheorem{cor}[prout]{Corollaire}

\newtheorem{sprout}{Théorème} [section]
\newtheorem{sprop}[sprout]{Proposition}

\newtheorem{scor}[sprout]{Corollaire}

\newcommand{\I}{\mathcal I}
\newcommand{\ib}{\overline{\mathcal I}}

\newcommand{\fix}{\mathrm{Fix}}
\newcommand{\stab}{\mathrm{Stab}}
\newcommand{\aut}{\mathrm{Aut}}
\newcommand{\card}{\mathrm{Card}}
\newcommand{\codim}{\mathrm{Codim}}

\newcommand{\eng}[2]{\langle  #1 \; | \; #2 \rangle } 

\newcommand{\F}{\mathcal F}
\newcommand{\cl}{\text{Cl}}
\newcommand{\V}{\mathcal{V}}
\newcommand{\D}{\mathcal D}
\newcommand{\M}{\mathcal M}
\newcommand{\U}{\mathcal U}
\newcommand{\A}{\mathcal A}
\newcommand{\B}{\mathcal B}
\newcommand{\fl}[1]{\overrightarrow{#1}}
\newcommand{\chap}{\widehat}
\newcommand{\tilda}{\widetilde}
\newcommand{\infini}{\infty}
\newcommand{\vide}{\emptyset}
\newcommand{\fleche}{\rightarrow}
\newcommand{\implique}{\Rightarrow}
\newcommand{\equi}{\:\Leftrightarrow\:}
\newcommand{\parallele}{\varparallel}
\newcommand{\barre}{\overline}
\newcommand{\tq}{\text{ tq }}
\newcommand{\et}{\text{ et }}

\newcommand{\vect}{\mathrm{Vect}}
\newcommand{\aff}{\mathrm{Aff}}
\newcommand{\conv}{\text{Conv}}
\newcommand{\proj}{\text{proj}}

\newcommand{\demo}{\textit{Démonstration: }}
\newcommand{\pv}{\textit{Preuve du lemme: }}
\newcommand{\rema}{\textit{Remarque: }}
\newcommand{\N}{\mathbb N}
\newcommand{\isom}{\xrightarrow{\sim}}
\newcommand{\R}{\mathbb R}
\newcommand{\ens}[2]{ \{ #1 \; | \; #2 \} }

\newcommand{\Z}{\mathbb Z}
\newcommand{\cqfd}{\hfill $\boxempty$}
\newcommand {\inv} {^{-1}}

\newcommand {\fonc}[5] { $#1$:
$\begin{array}{rrcl}
#2 & \rightarrow & #3 \\
#4 & \mapsto     & #5
\end{array}$ }

\renewcommand{\abstract}[1] { \begin{center}\textbf{Abstract} \end{center} \hspace{.5cm} #1} 
\newcommand{\resume}[1]  { \begin{center}\textbf{Résumé} \end{center} \hspace{.5cm} #1}

\begin{document}
\author{Cyril Charignon}
\title{Compactifications polygonales d'un immeuble affine}
\maketitle
\resume{A partir d'une décomposition en cônes de l'espace directeur d'un appartement d'un immeuble affine localement fini, on définit une compactification de l'immeuble semblable à la compactification de Satake d'un espace symétrique. Comme cas particuliers de cette construction, on retrouve la compactification polygonale classique telle que décrite dans \cite{landvogt} ou ses généralisation décrites dans \cite{werner}. Un intérêt de la construction présentée ici est qu'elle est totalement géométrique: elle est indépendante de l'existence d'un groupe agissant sur l'immeuble. On prouve au passage plusieurs résultats permettant d'identifier certaines parties de l'immeuble qui sont incluses dans un appartement, par exemple on prouve que deux facettes de quartier sont, quittes à être réduites, incluses dans un même appartement.}

\abstract{We define a compactification of an affine building $\I$ indexed by a family of partitions of the director space $\vec A$ of one of its appartments $A$. This compactification is similar to Satake's compatification of a symetric space, and it generalizes the quite well known polygonal compactification of an affine building in the sense that it is independant of the action of a group on the building, and that it allows some variations depending on the choice of the partition of $\vec A$. The different choices will mainly lead to different subgroups of the Weyl group acting on the border of $A$. Along the proofs, we get some results to help one find subsets of the building wich are included in an apartment, for exemple we prove that two sector facets can always be reduced so that they fit in one apartment.}

\tableofcontents

\section{Introduction}
 Le but de cet article est de présenter une construction purement géométrique de la compactification polygonale d'un immeuble affine localement fini, comme celle décrite par Landvogt \cite{landvogt}. On s'affranchit totalement de l'usage d'un groupe agissant sur l'immeuble, ce qui permet de définir une compactification pour les quelques immeubles qui pourraient ne pas être associés à un groupe muni d'une BN-paire.\\
 On propose également une généralisation: on définit toute une famille de compactifications, indexée par l'ensemble des partitions de l'espace directeur d'un appartement en cônes assujettis à certaines conditions.
Parmi elles se trouvent la compactification polygonale classique définie dans \cite{landvogt} ou \cite{guivarch-remy} ainsi que celles décrites par Annette Werner dans \cite{werner}.\\

Décrivons rapidement l'espace obtenu. L'adhérence de chaque appartement sera l'espace compact obtenu en rajoutant à l'infini un polygone, dont chaque face sera un complexe de Coxeter avec comme groupe un sous-groupe du groupe de Weyl de l'immeuble. On peut choisir n'importe quel polygone stable par le groupe de Weyl vectoriel, mais pour un choix quelconque, la structure de complexe de Coxeter sera triviale sur de nombreuses faces. La compactification de l'immeuble est obtenue en étendant la compactification d'un appartement de référence, d'une manière qu'on prouvera être unique dans la majorité des cas. Le bord ainsi rajouté à l'immeuble est en fait une réunion d'immeubles affines, dont les groupes de Coxeter sont des sous-groupes du groupe de Coxeter de l'immeuble de départ. Le choix de la décomposition en cônes au départ détermine lesquels de ces sous-groupes interviennent.\\

Par ailleurs, quelques résultats intermédiaires peuvent avoir leur intérêt propre, certains sont valides dans le cadre d'un immeuble quelconque. Il s'agit principalement des résultats de la partie \ref{section:resultats generaux}, qui donnent des critères pour s'assurer qu'une partie d'un immeuble est incluse dans un seul appartement. On prouve notamment que deux galeries tendues sont, quitte à être réduites, incluse dans un même appartement (du système complet d'appartements), puis on généralise au cas d'une galerie tendue et d'une cheminée. Ceci prouve par exemple que deux facettes de quartiers quelconques contiennent des sous-facettes incluses dans un même appartement, généralisant le résultat similaire déjà connu pour les quartiers (voir par exemple la proposition (9.5) de \cite{ronan}). Quelques résultats classiques sur l'enclos d'une partie sont également prouvés, dans le cas d'une partie ne coupant aucune chambre.\\

Dans la partie, \ref{section:donnee initiale}, on énonce et on analyse brièvement les conditions requises sur une décomposition en cônes de l'espace directeur d'un appartement pour définir une compactification de l'immeuble. Dans la partie suivante, on définit la compactification d'un appartement. L'interlude de la partie \ref{section:resultats generaux} permet d'énoncer les quelques résultats nécessaires à la suite qui peuvent avoir un intérêt propre. La partie \ref{section:ensemble} définit l'ensemble qui sera l'immeuble compactifié, la partie \ref{section:topologie} définit la topologie sur cet ensemble et prouve les propriétés attendues, en particulier la compacité. Dans la partie \ref{section:unicite}, on vérifie que la compactification de l'immeuble ainsi définie est unique lorsqu'une compactification d'un appartement est fixée et que l'immeuble provient d'un groupe muni d'une donnée radicielle. Cela fournit un moyen simple de comparer cette compactification avec d'autres, comme celles définies dans \cite{landvogt},\cite{guivarch-remy} et \cite{werner}. Enfin, la partie \ref{section:bord} décrit le bord qu'on vient de rajouter à l'immeuble: il s'agit d'une réunion d'immeubles affines de dimensions inférieures.\\

\section{La donnée initiale}
\label{section:donnee initiale}

\subsection{Conventions, notations}
Un immeuble $\I$ sera vu a priori comme un complexe simplicial vérifiant les axiomes classique (\cite{brown}, \cite{tits74}), ou comme un complexe polysimplicial comme dans \cite{bruhat-tits}. Chaque appartement $A$ est un complexe de Coxeter, dont on note $W(A)$ le groupe de Coxeter. Ce groupe est indépendant de $A$ à isomorphisme près. Lorsqu'on choisit une chambre $c$ de $A$, les réflexions par rapport aux cloisons de $c$ forment une partie $S(c)$ de $W(A)$ telle que $(W(A),S(c))$ est un système de Coxeter. La classe d'isomorphisme de $W(A)$ est caractérisée par un diagramme appelé diagramme de Coxeter. Son ensemble de sommets est en bijection avec $S(c)$, et les sommets correspondant aux réflexions $s$ et $t$ sont reliés si $s$ et $t$ ne commutent pas. On précise alors sur l'arête l'ordre de $st$. Ce diagramme ne dépend pas du choix de $A$ ni de $c$.\\
 On définit sur l'ensemble des facettes le type, c'est une fonction à valeur dans $\mathcal P(S(c))$ qui permet de caractériser les orbites des facettes d'un appartement sous l'action du groupe de Weyl. Les isomorphismes entre appartements préservent le type, par définition.\\
 L'ensemble des chambres de $\I$ est muni d'une distance à valeur dans $W(A)$, appelée $W$-distance. En fait, $\I$ est totalement déterminé par l'ensemble de ses chambres muni de sa $W$-distance, c'est d'ailleurs le point de vue adopté dans \cite{ronan}.\\

Nous nous intéressons principalement aux immeubles affines. Dans ce cas, on identifie $\I$ à sa réalisation géométrique affine, qui est un espace métrique complet dans laquelle les appartements sont des espaces affines euclidiens, les groupes de Coxeter des groupes de transformations orthogonales affines, et les isomorphismes entre appartements des isométries. Voir \cite{brown}, chapitre VI, \cite{bruhat-tits}, \cite{tits86}, \cite{parreau},\cite{euclidean-buildings}. C'est cet espace qu'on se propose de compactifier. Si $A$ est un appartement, le groupe $W(A)$ est engendré par les réflexions orthogonales par rapport à des hyperplans appelés murs (\cite{bourbaki}). Les murs définissent une partition de $A$ dont les parties sont identifiées aux facettes de $A$, les chambres étant les facettes de dimension maximale.\\

 On note $\vec A$ l'espace directeur de $A$, c'est lui aussi un complexe de Coxeter dont le groupe $W(\vec A)$ est l'ensemble des parties vectorielles des éléments de $W(A)$. Il est engendré par les réflexions orthogonales par rapport aux murs de $\vec A$, qui sont les espaces directeurs des murs de $A$. Les murs de $\vec A$ définissent une partition de $\vec A$ en cônes convexes appelés facettes de Weyl ou facettes vectorielles. Ces parties sont identifiées aux simplexes du complexe de Coxeter $\vec A$.\\
 Il existe un système de racine $\phi\subset \vec A^*$ tel que les murs de $\vec A$ sont les noyaux des racines et dont $W(A)$ est le groupe de Coxeter affine associé, au sens de \cite{bourbaki}, 2.5.\\
 Si on fixe une chambre de Weyl $C\subset \vec A$, on note, pour $s\in S(C)$ $\alpha_s$ la racine correspondante (donc $\ker(\alpha_s)=\fix(s)$). Alors $\{\alpha_s\}_{s\in S(C)}$ est une base de $\phi$ et $C=\ens{x\in\vec A}{\forall s\in S(C),\; \alpha_s(x)>0}$. De plus, dans le diagramme de Coxeter de $\vec A$, deux sommets $s$ et $t$ sont reliés si et seulement si $\alpha_s$ n'est pas orthogonale à $\alpha_t$. En considérant ce diagramme comme un simple graphe, on définit les notions de connexité habituelles.\\

Lorsqu'on passe du point de vue complexe simplicial au point de vue géométrique d'un immeuble affine, le vocabulaire change un peu:\\
 Un complexe simplicial est en particulier un ensemble muni d'une relation d'ordre. On notera généralement $\sqsubset$ cette relation, et on dira que $c$ est inclus dans $d$ lorsque $c\sqsubset d$. Mais lorsque $c\sqsubset d$ dans un immeuble affine, alors dans la réalisation géométrique $c$ est inclus dans l'adhérence de $d$, avec $c\subset d$ si et seulement si $c=d$. En général, $c\subset\partial d$.\\
 De plus, dans un complexe simplicial, la facette maximale inférieure à $c$ et à $d$, notée $c\wedge d$ est appelée l'intersection de $c$ et $d$. Dans la réalisation géométrique, l'adhérence de $c\wedge d$ est égale à l'intersection des adhérences de $c$ et $d$, mais $c\wedge d$ ne s'exprime pas de manière directe en fonction de $c$ et $d$ (en fait, $c\wedge d$ est l'intérieur de $\bar c\cap\bar d$ dans $\vect(\bar c\cap\bar d)$).\\

 Si $M$ est un mur et $c$ une chambre dans un système de Coxeter $A$, on notera $\D(M,c)$ le demi-appartement fermé délimité par $M$ et contenant $c$. Sauf précision, un demi-appartement signifiera un demi-appartement fermé. Un demi-appartement peut aussi être défini à l'aide de deux chambres adjacentes $c$ et $d$, ou d'une racine $\alpha$ (on voit a priori les racines comme des formes affines sur un appartement) et d'un entier $k$: on notera $\D(c,d)$ la réunion des chambres fermées plus proches de $c$ que de $d$, et $\D(\alpha,k)=\{x|\alpha(x)+k\geq 0\}$.\\
 Soient $D^+$ et $D^-$ les deux demi-appartements définis par un mur $M$ dans un appartement $A$. On dira que deux parties $\alpha$ et $\beta$ de $A$ sont séparées par $M$ si $\alpha\subset D^+$ et $\beta \subset D^-$, ou l'inverse. On notera alors $\alpha |^M \beta$ (Donc par exemple, $\alpha\subset M\implique \alpha|^M \beta$ quel que soit $\beta$). On dira que ces parties sont séparées strictement par $M$ si en outre $\alpha\not\subset M$ et $\beta\not\subset M$.\\

 L'enclos d'une partie $E$ dans un appartement $A$ est l'intersection de tous les demi-appartements fermés de $A$ contenant $E$. On la note $\cl_A(E)$.\\
 
 Lorsque $A$ et $B$ sont deux appartements ayant au moins une chambre en commun, on peut naturellement identifier $\vec A$ avec $\vec B$. Lorsque la dimension de $A\cap B$ est moindre, on ne peut qu'identifier un sous-espace de $\vec A$ avec un sous-espace de $\vec B$, voici comment on procède:\\
Les sous-espaces $\vect\ens{\vec{xy}\in\vec A}{x,y\in A\cap B}\subset \vec A$ et $\vect\ens{\vec{xy}\in\vec B}{x,y\in A\cap B}\subset \vec B$ sont canoniquement isomorphes. On identifie alors ces deux sous-espaces, et on note l'espace obtenu $\vec A\cap\vec B$. Cet espace vérifie la propriété suivante: si $E$ est un troisième appartement, si $\phi:A\isom E$ et $\psi:B\isom E$ sont deux isomorphismes qui co\"{\i}ncident sur un ensemble $b\subset A\cap B$, alors l'espace directeur $\vec b$ de $\aff(b)$ est inclus dans $\vec A\cap\vec B$, et pour tout $\vec v\in \vec b$, on a $\vec\phi(\vec v)=\vec\psi(\vec v)$.\\
 Cette "intersection des espaces directeurs" ne vérifie pas l'associativité, en fait l'écriture $(\vec A\cap \vec B)\cap \vec C$ n'a même aucun sens, car $\vec A\cap \vec B$ n'est pas uniquement défini au moyen de la structure vectorielle de $\vec A$ et $\vec B$, mais bien de la structure d'espaces affines de $A$ et $B$. (Une notation comme $A\vec\cap B$ serait sans doute plus appropriée.)\\
\rema Dans la suite, $\vec b$ ne signifiera pas en général l'espace directeur de $\aff(b)$.\\

\subsection{Cônes convexes}

 Dans cet article, tous les cônes seront supposé convexes a priori. Un cône vectoriel convexe dans un $\R$-espace vectoriel $\vec E$ est un sous-ensemble de $\vec E$ stable par addition et multiplication par un scalaire strictement positif. Tout cône vectoriel contient $0$ dans son adhérence. Si $E$ est un espace affine dirigé par $\vec E$, un cône convexe affine de $E$ est un sous-ensemble de $E$ de la forme $f=s+\vec f$ où $s\in E$ et $\vec f$ est un cône convexe de $\vec E$.\\
Lorsque $\vec f$ est un cône de $\vec E$ ne contenant pas de droite (on dit aussi "cône pointu"), alors l'écriture $f=s+\vec f$ est unique, ce qui permet de définir $s$ comme étant le \textit{sommet} de $f$, noté $s(f)$.\\
Le cône vectoriel $\vec f$ est quand à lui toujours bien déterminé, on l'appelle la \textit{direction} de $f$. Deux cônes ayant la même direction sont dit \textit{parallèles}. Si $g$ est un cône parallèle à $f$ et $g\subset f$, alors $g$ est un sous-cône parallèle de $f$, et on abrège "sous-cône parallèle" en "scp".\\

\rema Dans \cite{bruhat-tits}, deux parties d'un appartement égales à translation près sont appelées équi\-pol\-lentes au lieu de parallèles.\\

\subsection{Décomposition d'un appartement en cônes}
 On fixe désormais un appartement $A_0$. On notera $W=W(A_0)$ son groupe de Coxeter, et $W^v=W(\fl A_0)$ son groupe de Coxeter vectoriel. On choisit un ensemble $\F$ de parties non vides de $\fl{A_0}$ vérifiant les propriétés suivantes:\\

(H1) $\fl{A_0}=\bigsqcup_{\vec f\in \F} \vec f$. \hspace{1cm}(Le symbole $\sqcup$ signifie "réunion disjointe".)\\

(H2) $\F$ est fini.\\

(H3) $\{0\}\in\F$.\\

(H4) Chaque élément de $\F$ est décrit par un système d'équations et d'inéquations linéaires: pour tout $\vec f\in\F$, il existe $n\in \N$, $\alpha_1,..., \alpha_n \in \fl{A_0}^*$, $r\in \N$ tels que $\vec f=\{x\in \fl{A_0} | \, \alpha_i(x)=0, \, \forall i \in \llbracket 1,r\rrbracket , \,\mathrm{ et }\, \alpha_j(x)>0,\, \forall i\in \llbracket r+1,n\rrbracket \}$. En particulier, chaque élément de $\F$ est un cône convexe, ouvert dans son support.\\

(H5) Le bord d'un cône $\vec f$ de $\F$ est une réunion d'autre cônes de $\F$, qu'on appelle les faces de $\vec f$.\\

(H6) Si $\vec f,\vec g \in \F$ et si $\vec f$ est une face de $\vec g$, alors $\overline{\vec f}=Vect(\vec f) \cap \overline{\vec g}$.\\

Lorsqu'on parle du bord d'un cône, on sous-entend ici le bord dans l'espace vectoriel qu'il engendre. Pour un cône $\vec f$ vérifiant (H4), qui est donc ouvert dans $\vect(\vec f)$, on a $\partial \vec f=\barre{\vec f}\setminus \vec f$.\\
 A partir de ces données, on va définir une compactification de $A_0$. Dès qu'on voudra l'étendre en une compactification de $\I$, il faudra en outre supposer: \\

 (H7) $\F$ est stable par le groupe de Weyl vectoriel $W^v$.\\

\subsection{Conséquences directes des hypothèses sur $\F$}
\label{subsection:consequences directes}

 Chaque élément de $\F$ est un cône convexe, ouvert dans l'espace vectoriel qu'il engendre.
 De plus, $\{0\}$ est le seul cône à contenir $0$, donc aucun élément de $\F$ ne contient de droite, ce qui permet de définir les sommets des cônes affines de direction un élément de $\F$.\\

\rema La condition (H3) a en fait pour unique but de permettre de définir le sommet d'un cône pour faciliter les raisonnements dans la suite, mais elle semble superflue. Étudions brièvement le cas général. Soit $f\in \F$ le cône contenant $0$. Comme $f$ est un cône ouvert dans $\vect(f)$, on a $f=\vect(f)$, c'est-à-dire que $f$ est un espace vectoriel. Comme l'adhérence de tout cône vectoriel contient $0$, par (H5) on voit que $f$ est dans le bord de chaque élément de $\F$. On vérifie alors que chacun de ces éléments est stable par addition par $f$, on peut donc tout quotienter par $f$ pour obtenir un espace vectoriel muni d'une décomposition en cônes vérifiant cette fois toutes les hypothèses (H1) - (H6). Si $\partial (A_0 / f)$ est le bord qu'on va définir dans la section \ref{section:compactification de lappart}, alors la même procédure appliquée à $A_0$ et $\F$ conduirait à rajouter exactement le même bord. En fait, la condition (H3) impose de rajouter à $A_0$ un bord de codimension 1 et non supérieure.\\

 Lorsque $(\alpha_i)_{i\in I\sqcup J}$ est une famille de formes linéaires définissant $\vec f$ comme dans la quatrième hypothèse sur $\F$, on notera juste
 \[\vec f=\{\alpha_i>0,\alpha_j=0,\, i\in I,\, j\in J\}\]
 La famille $(\alpha_i)_{i\in I\sqcup J}$ est nécessairement génératrice de $\fl{A_0}^*$ sans quoi $\vec f$ ou une de ses faces contiendrait un sous-espace vectoriel de $\fl{A_0}$ non réduit à $\{0\}$.\\

 Lorsque $\vec g$ est une face de $\vec f$, alors il existe une famille $(\alpha_i)_{i\in I\sqcup J\sqcup K}$ telle que:
\[  \begin{cases}		 \vec f=\{\alpha_i>0,\alpha_j=0,\, i\in I\sqcup J,\, j\in K\} \\
						 \vec g=\{\alpha_i>0,\alpha_j=0,\, i\in I,\, j\in J\sqcup K\} \end{cases} \]

 \subsection{Exemples}
\label{subsection:exemple}

 Un premier exemple de telle décomposition de $\fl{A_0}$ en cônes est la décomposition en facettes de Weyl, notée $\F^\vide$. Les cônes affines dont les directions sont dans $\F^\vide$ sont les facettes de quartier, et la compactification qu'on obtiendra alors est la compactification polygonale classique, décrite dans \cite{landvogt}.\\
 Un exemple un peu plus général est celui considéré par Annette Werner dans \cite{werner}, où il s'agit grosso modo d'enlever à la décomposition en facettes de Weyl les cloisons d'un certain type. C'est cet exemple que je développe ici.\\

\subsubsection{Décompositions en cônes obtenues à partir d'une partie de $S$}

 On fixe une chambre de Weyl $C_0\subset\fl{A_0}$, soit  $S$ l'ensemble des réflexions par rapport aux cloisons de $C_0$, de sorte que les facettes de $\fl{A_0}$ sont typées par les parties de $S$. Soit $J$ une partie de $S$, l'idée est de rassembler les chambres séparées par une cloison de type $\{j\}$ avec $j\in J$. Pour assurer la convexité, il faut penser à rajouter alors les facettes de dimension plus petite qui se trouvent entre plusieurs chambres rassemblées. Pour s'assurer que (H5) et (H6) seront vérifiées, il faut aussi rassembler les facettes bordant plusieurs chambres rassemblées, lorsqu'elles engendrent le même espace vectoriel. On arrive à la définition suivante:

\begin{defin}
 Si $f$ une facette de Weyl de $\fl{A_0}$ de type $I\subset S$, on note $\bar f^J$ la réunion de $f$ et des facettes de son bord de type inclus dans $I\cup (J\cap I^\perp)$.\\
 \end{defin}

 Si $h$ et $f$ sont deux facettes d'une même chambre fermée, avec $h$ de type $I'$ et $f$ de type $I$, alors $h$ est incluse dans $\bar f^J$ si et seulement si $I'$ est la réunion de $I$ et d'une partie de $J$ disconnectée de $I$. Donc une facette $h$ n'est incluse dans aucun $\bar f^J$, avec $f\not= h$ si et seulement si son type ne contient aucune composante connexe incluse dans $J$. Une telle facette sera dite \textit{admissible}. On dira également que son type est admissible, de sorte qu'une facette est admissible ssi son type est admissible.\\
\begin{defin}
Si $f$ est une facette admissible de $\barre{C_0}$, on pose $J.f=W^v_{J\cap I^\perp}.\bar f^J$.\\
On pose ensuite $\F^J= W^v. \ens{J.f}{ f\text{ facette admissible de $C_0$}}$.\\
\end{defin}

\rema Conformément à la notation déjà introduite, l'ensemble des facettes de Weyl est $\F^\vide$.\\

\begin{prop} L'ensemble $\F^J$ est un ensemble de cônes, qui vérifie les hypothèses (H1)-(H2) et (H4-H7). Chaque facette de Weyl $f$ est incluse dans un unique cône de $\F^J$, qu'on notera $J.f$. Lorsque $J$ ne contient aucune composante connexe de $S$, alors $\F^J$ vérifie également (H3).\end{prop}

\demo\\
Les points (H2), et (H7) sont évidents, (H3) est vrai si et seulement si $\{0\}$ est une facette admissible, ce qui équivaut bien au fait que $J$ ne contient aucune composante connexe de $S$. Il est également clair que $\vec A=\bigcup_{f\in\F^J} f$, et pour prouver (H4) et (H5), il suffit de considérer des éléments de $\F^J$ du type $J.f$, avec $f$ une facette admissible de $C_0$.\\

 On commence par (H4). Soient $\{\alpha_i\}_{i\in S}$ l'ensemble des racines délimitant $C_0$. Soit $I\subset S$ et $f$ la facette de $C_0$ de type $I$, alors
 \[f=\{ \alpha_i=0,\; \alpha_j>0,\; i\in I,\; j\in S\setminus I\},\]
 et \[\bar f^J=\{ \alpha_i=0,\; \alpha_k>0,\; \alpha_l\geq 0,\; i\in I,\; l\in J\cap I^\perp,\; k\in S\setminus(I\cup(J\cap I^\perp)) \}.\]
 Notons $L=J\cap I^\perp$ et $K=S\setminus(I\cup(J\cap I^\perp))$. Montrons que $J.f=\{\alpha_i=0,\; w(\alpha_k)>0,\; i\in I,\; k\in K,\; w\in W^v_{J\cap I^\perp} \}$. En attendant, notons $E$ ce dernier ensemble. C'est un ensemble délimité par des murs, donc une réunion de facettes. De plus, pour tout $w\in W^v_{J\cap I^\perp}$ et tout $i\in I$, $w(\alpha_i)=\alpha_i$, on voit donc que $E$ est stable par $W^v_{J\cap I^\perp}$.\\

Pour montrer que $J.f\subset E$, il suffit donc de montrer que $\bar f^J\subset E$. Soit $x\in \bar f^J$, il vérifie déjà les conditions $\alpha_i(x)=0$, $i\in I$. Soit $k\in K$ et $w\in W^v_{J\cap I^\perp}$. On sait que $\{\alpha_s\}_{s\in S}$ est une base du système de racines de $\fl{A_0}$. Or, $w\in W_{J\cap I^\perp}=W_L$, et $L\cap K=\vide$. Donc $w(\alpha_k)=\alpha_k+\sum_{l\in L} n_l\alpha_l$, avec $n_l\in\N$. Il apparaît ainsi que $w(\alpha_k)(x)>0$. Donc $x\in E$.\\
Montrons l'inclusion inverse. Soit $x\in E$, il existe $w\in W_L$ tel que $\forall l\in L$, $\alpha_l(w.x)\geq 0$. Comme $E$ est stable par $w$, $w.x\in E$, et ainsi $w.x$ vérifie toutes les inégalités prouvant qu'il appartient à $\bar f^J$. Donc $x\in W_L.\bar f^J=J.f$.\\
% on peut trouver une chambre $D$ telle que $x\in \bar d$ et $\forall i\in I$, $\alpha_i(D)>0$. Soit $\Gamma$ une galerie tendue de $C_0$ à $D$, soient $s_1,...,s_n$ la suite des types des cloisons de $\Gamma$, soit $w=s_1...s_n$ de sorte que $w\inv.x\in \bar C_0$. Il ne reste plus qu'à vérifier que $w\inv\in W^v_{J\cap I^\perp}$ et $w\inv.x\in \bar f^J$.\\
% Pour commencer, le mur $\ker(\alpha_{s_1})$ sépare $C_0$ et $D$. Il ne peux donc être de type un élément de $I$, ni de $K$ car $\alpha_k(x)>0$ pour tout $k\in K$. Donc $s_1\in J\cap I^\perp$. Le mur suivant est $s_1.\ker(\alpha_{s_2})=\ker(\alpha_{s_2}\circ s_1)= \ker(s_1.\alpha_{s_2})$. Encore une fois, ce mur ne peut être de type un élément de $I$, ni de $K$, car $s_1.\alpha_k(x)>0$ pour tout $k\in K$. Donc $s_2\in J\cap I^\perp$. Ainsi de suite, on vérifie que pour tout $u\in \llbracket 1,n\rrbracket $, $s_u\in J\cap I^\perp$, ce qui prouve que $w\inv \in W^v_{J\cap I^\perp}$.\\
% Alors pour tout $i\in I$, $w(\alpha_i)=\alpha_i$, donc $\alpha_i(w\inv x)=0$. Comme en outre, pour tout $k\in k$, $\alpha_k(w\inv x)>0$ et $w\inv x\in \bar C_0$, on a bien $w\inv x\in \bar f^J$. Ce qui prouve que $x\in J.f$.\\

 A présent, prouvons que chaque facette de Weyl est incluse dans un unique cône de $\F^J$. Ceci impliquera directement (H1) car les éléments de $\F^J$ sont des réunions de facettes de Weyl.\\
 Soit $f\in\F^\vide$  une facette de Weyl, incluse dans deux éléments de $\F^J$, disons $w_1.J.g_1$ et $w_2.J.g_2$. En translatant tout par un élément de $W^v$, on peut supposer $f\subset \bar C_0$. Soit $I_1$ le type de $g_1$ et $I_2$ celui de $g_2$. Il existe $v_1\in W^v_{J\cap I_1^\perp}$ et $v_2\in W_{J\cap I_2^\perp}$ tels que $f\subset w_1 v_1\bar g_1^J$ et $f\subset w_2 v_2\bar g_2^J$. Comme $f$, $g_1$ et $g_2$ sont des facettes de $\bar C_0$, ceci implique en fait $f\subset \bar g_1^J\cap \bar g_2^J$. Le type de $f$ s'écrit donc $I'=I_1\sqcup J_1=I_2\sqcup J_2$ où $J_1$ et $J_2$ sont des parties de $J$ disconnectées respectivement de $I_1$ et de $I_2$. Soit $K$ une composante connexe du type de $f$, elle est incluse soit dans $I_1$ soit dans $J_1$. De deux choses l'une: soit $K\subset J$ et alors $K$ ne peut être incluse dans $I_1$ car $g_1$ est admissible, soit $K\not\subset J$ et alors elle ne peut être incluse dans $J_1$. On voit donc que $J_1$ est la réunion des composantes connexes de $I'$ incluses dans $J$, $I_1$ est la réunion des autres composantes connexes. Le même raisonnement est valable pour $I_2$ et $J_2$, ce qui prouve $I_1=I_2$ et $J_1=J_2$. D'où $g_1=g_2$ et $J.g_1=J.g_2$. Il reste à regarder $w_1$ et $w_2$. On sait que $w_1 v_1\in \fix(f)=W_{I'}$, et $I'=I_1\cup J_1\subset I\cup (J\cap I^\perp)$. D'où $w_1 v_1\in W^v_{I\cup (J\cap I^\perp)}$, d'où $w_1 J.g_1=J.g_1$. De même, $w_2 J.g_2=J.g_2=J.g_1=w_1 J.g_1$. Ce qui prouve que $f$ est incluse dans un unique cône de $\F^J$.\\

 Prouvons (H5). Soit $f\in \F^\vide$ une facette de Weyl admissible qu'on peut supposer incluse dans $C_0$, soit $I$ son type. Comme $J.f$ est une union de facettes de Weyl, son bord l'est aussi. Soit $g$ une facette de Weyl bordant $J.f$. Quitte à translater $g$ par un élément de $W^v_{J\cap I^\perp}$ qui stabilise $J.f$, on peut supposer $g\subset \bar C_0$, et donc $g\subset \bar f$. Comme $J.f$ est ouvert dans $\vect(J.f)$, on a même $g\subset \bar f\setminus \bar f^J$. Notons $J.g$ le cône de $\F^J$ contenant $g$, il s'agit de prouver que $J.g\subset \partial J.f$. Soit $I'$ le sous ensemble de $S$ obtenu en retirant au type de $g$ toutes ses composantes connexes incluses dans $J$. Soit $h$ la facette de $\bar C_0$ de type $I'$, alors $h$ est admissible et $g\subset \bar h^J\subset J.h$, donc $J.g=J.h$. Si $J_1$ est une composante connexe de $\text{type}(g)$ incluse dans $J$ alors $J_1\cap I$ est une réunion de composantes connexes de $I$ incluses dans $J$, d'où $J_1\cap I=\vide$ car $I$ est admissible. Ceci prouve que $I\subset I'$, et donc $h\subset \bar f$. Ensuite, $J.h=W^v_{J\cap I'^\perp}.\bar h^J\subset \barre{J.f}$ car $W^v_{J\cap I'^\perp}\subset W^v_{J\cap I^\perp}$ et $\bar h^J\subset \bar f$. Et comme $J.h\not= J.f$ car $g\subset J.h$, on a $J.h\cap J.f=\vide$ d'où $J.h\subset \partial J.f$.\\

Il ne reste plus qu'à prouver (H6). Soit $J.h$ une face d'un cône $J.f$. Comme dans le paragraphe précédent, on peut supposer que $h$ et $f$ sont des facettes admissibles de $\bar C_0$. Soit $I$ le type de $f$ et $I'$ le type de $h$. Nous procédons par récurrence sur $\card(I'\setminus I)$, en commençant par montrer l'hérédité.\\

 On suppose donc $\card(I'\setminus I)\geq 2$. Comme $h\not\subset \bar f^J$, il existe $l\in I'\setminus (I\cup (J\cap I^\perp)$. Soit $g$ la facette de $\bar C_0$ de type $I\cup\{l\}$, elle est admissible. Par récurrence on a $\barre{J.g}=\barre{J.f}\cap \vect(J.g)$. De plus, $J.h$ est une face de $J.g$, différente de $J.g$ car $h$ est admissible et ne peut donc être incluse dans $\bar g^J$. Donc par récurrence, $\barre{J.h}=\barre{J.g}\cap \vect(J.h)$. La combinaison des deux égalités donne bien $\barre{J.h}=\barre{J.f}\cap \vect(J.h)$.\\

 Traitons maintenant le cas $\card(I'\setminus I)=1$. Dans l'égalité $\barre {J.h}=\vect(J.h)\cap \barre{J.f}$, l'inclusion $"\subset"$ est évidente.\\
 Soit $x\in \vect(J.h)\cap \barre{J.f}$, cela signifie, d'après la description de $J.f$ obtenue pour prouver (H4): \begin{itemize}
\item pour tout $i\in I'$, $\alpha_i(x)=0$
\item pour tout $w\in W^v_{J\cap I^\perp}$ et $k\in S\setminus(I\cup (J\cap I^\perp))$, $w.\alpha_k(x)\geq 0$. \end{itemize}

Et le but est de prouver:
\begin{itemize}
\item pour tout $i\in I'$, $\alpha_i(x)=0$
\item pour tout $w\in W^v_{J\cap I'^\perp}$ et $k\in S\setminus(I'\cup (J\cap I'^\perp))$, $w.\alpha_k(x)\geq 0$. \end{itemize}
 Parmi ces inégalités, celles qui ne sont pas directement dans la liste des hypothèses sont les $w.\alpha_k(x)\geq 0$ lorsque:
\[k\in 
(S\setminus(I'\cup (J\cap I'^\perp))) \, \setminus \, (S\setminus(I\cup (J\cap I^\perp)))= (I\cup(J\cap I^\perp))\setminus (I'\cup(J\cap I'^\perp))= (J\cap I^\perp)\setminus (I'\cup (J\cap I'^\perp))=(J\cap I^\perp)\setminus (I'\cup I'^\perp)\]
 et lorsque $w\in W^v_{J\cap I'^\perp}$.\\

 Soit $k$ un tel indice et $w\in W^v_{J\cap I'^\perp}$. En particulier, $k\in W^v_{J\cap I^\perp}$.
 Soit $l$ tel que $I'\setminus I=\{l\}$. Comme $I'$ est admissible, $l\not\in J\cap I^\perp$.

 Alors $wk\in W^v_{J\cap I^\perp}$, et $l\in S\setminus(I\cup (J\cap I^\perp))$ d'où d'après les hypothèses $wk.\alpha_l(x)\geq 0$. Mais:
\[wk.\alpha_l(x)=\alpha_l(k w\inv x)=\alpha_l(w\inv x - 2\alpha_k(w\inv x).\alpha_k^\vee)=
 \alpha_l(w\inv x) -2w\alpha_k(x).\langle\alpha_k|\alpha_l\rangle\]
 Comme $l\in I'$ et $w\inv x\in \vect(J.h)$, on a $\alpha_l(w\inv x)=0$. Il reste donc $-2w\alpha_k(x).\langle\alpha_k|\alpha_l\rangle \geq 0$.\\
Mais $k\not=l$ car $l\in I'$ et $k\not\in I'$. De plus, si $k$ était orthogonal à $l$, il serait orthogonal à $I'$, ce qui est exclus. Ainsi, $\langle\alpha_k|\alpha_l\rangle<0$, ce qui entraîne $w\alpha_k(x)\geq 0$.\cqfd\\

\subsubsection{Comparaison avec \cite{werner}}

 On peut dès à présent vérifier de façon élémentaire que les partitions de $\fl{A_0}$ en parties convexes définies dans \cite{werner} sont précisément du type précédent. Ces parties sont définies, pour l'immeuble de Bruhat-Tits d'un groupe $G(K)$ avec $K$ un corps local, à partir d'une représentation linéaire $\rho$ fidèle et de dimension finie de $G(K)$, mais ne dépendent en fait que de la facette de Weyl de $A_0^*$ dans laquelle se trouve le plus haut poids $\lambda_0(\Delta)$ de $\rho$ une fois fixée une base $\Delta$ du système de racine (\cite{werner} théorème 4.5). Soit $J$ le type de cette facette (il est indépendant de $\Delta$). Comme $\rho$ est fidèle, $J$ ne contient pas de composante connexe de $S$. Pour $\Delta$ la base du système de racines correspondant à la chambre $C_0$, on a $J=\ens{s\in S}{\langle \alpha_s|\lambda_0(\Delta)\rangle=0}$. Nous allons vérifier que la partition $\Sigma(\rho)$ de $\fl{A_0}$ définie dans \cite{werner} est $\F^J$.\\
 Soit $\Delta$ une base du système de racines et $Y\subset \Delta$. Il est immédiat d'après la définition donnée dans \cite{werner} (définition 1.1) que $Y$ y est dit admissible si et seulement si il existe $I\subset S$ admissible (au sens défini plus haut) tel que $Y=\{\alpha_i\}_{i\in I}$.\\

 Les parties définies dans \cite{werner} sont notées $F^\Delta_Y$ pour $\Delta$ une base du système de racines et $Y\subset\Delta$ une partie admissible. Ce sont des réunions de facettes de Weyl (\cite{werner}, proposition 4.4), ouvertes dans leur support, qui réalisent une partition de $\fl{A_0}$ stable par le groupe de Weyl $W^v$.  Soit $\Delta$ la base correspondant à $C_0$, et $Y\subset \Delta$ une partie admissible. Soit $f$ la facette de $C_0$ de type $I$, avec $Y=\{\alpha_i\}_{i\in I}$, nous allons voir que $F^\Delta_Y=J.f$. Pour commencer, la proposition \cite{werner} 4.4 montre que $F_Y^\Delta\cap \barre{ C_0}=\bar f^J=J.f\cap \barre{C_0}$.\\
 Montrons que $F^\Delta_Y$ est stable par $W_{J\cap I^\perp}$. Soit $s\in J\cap I^\perp$, alors la facette de $C_0$ de type $I\cup\{s\}$ est incluse dans $\barre{C_0}\cap F^\Delta_Y$. Or $F^\Delta_Y$ est ouverte dans son support, qui est $\bigcap_{i\in I}\ker(\alpha_i)$ et donc qui est stable par $s$. Alors $F^\Delta_Y$ doit contenir la facette $s.f$. D'où $F^\Delta_Y\cap s.F^\Delta_Y\not=\vide$, d'où $s$ stabilise $F^\Delta_Y$. On déduit de ceci que $J.f\subset F^\Delta_Y$.\\
 Pour l'autre inclusion, soit $x\in F^\Delta_Y$. En choisissant $y\in f$ de manière générique, on peut s'assurer que $[x,y]$ ne rencontre que des facettes de dimension $\dim(f)-1$. Soit $g$ la première facette différente de $f$ rencontrée par ce segment en partant de $y$. Comme $[x,y]\subset F_Y^\Delta$, par convexité de $F_Y^\Delta$, on obtient que $g\subset \bar f\cap F_Y^\Delta=\bar f^J$. Comme $g$ est de codimension $1$ dans $f$, le type de $g$ est $I\cup\{s\}$, avec $s\in J\cap I^\perp$. Alors la facette $sf$ est incluse dans $J.f$, et elle contient "la suite" du segment $[x,y]$, c'est à dire un intervalle ouvert $]u,v[$ tel que $]u,v[\cup ([x,y]\cap g) \cup ([x,y] \cap f)$ est connexe. On a $J.f=J.(sf)$, et $F_{sY}^{s\Delta}=s.F_Y^\Delta=F_Y^\Delta$ (car ces deux cônes contiennent la facette $sf$). On peut appliquer la proposition \cite{werner} 4.4, dans la chambre $s.C_0$, on obtient que $F_Y^\Delta\cap s\bar C_0=J.f\cap s\bar C_0$ est la réunion des facettes de $s.f$ de type inclus dans $I\cup(J\cap I^\perp)$. Alors la prochaine facette rencontrée par le segment $[x,y]$ est de type $I\cup\{t\}$, avec $t\in J\cap I^\perp$, elle est bien incluse dans $J.f$, la facette de dimension maximale suivante est $st.f$ qui est aussi dans $J.f$ car $st\in W_{J\cap I^\perp}$. Ainsi de suite, on vérifie que tout le segment $[x,y]$ est inclus dans $J.f$ (puisque $[x,y]$ ne rencontre qu'un nombre fini de facettes).\\

\subsubsection{Dessins, autres exemples}
\label{subsection:dessins}
 Voici les dessins de quelques décompositions en cônes de l'appartement vectoriel de type $A_2$. On note $s$ et $t$ les réflexions par rapport aux cloisons de la chambre de base.\\
\begin{figure}[h!]
\includegraphics[width=5cm, height=5cm]{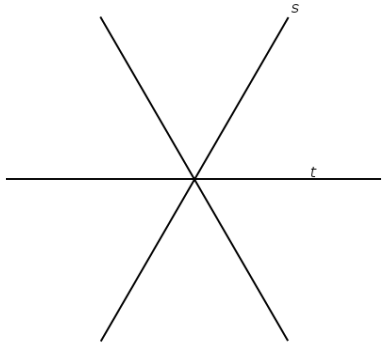}
\centering
\caption{Décomposition en facettes de Weyl. Je la laisse en arrière-plan dans les exemples suivants.}
\end{figure}

\begin{figure}[h!]
\includegraphics[width=5cm, height=5cm]{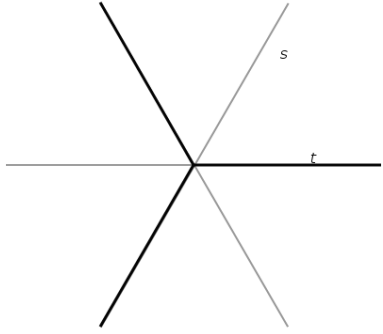}
\centering
\caption{Décomposition de type $\F^{ \{s\} }$: on retire les cloisons de type $s$.}
\end{figure}

\begin{figure}[h!]
\includegraphics[width=5cm, height=5cm]{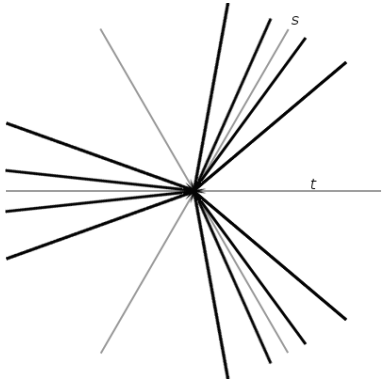}
\centering
\caption{Un autre exemple de décomposition vérifiant (H1)-(H7).}
\end{figure}

\newpage
\section{Compactification de $A_0$}
\label{section:compactification de lappart}

\subsection{L'ensemble $\overline{A_0}$ }

\begin{defins}\begin{itemize}

\item Un $\F$-cône affine dans $A_0$ est un cône dont la direction est dans $\F$. On note $\F_{A_0}$ l'ensemble des $\F$-cônes affines de $A_0$.
\item Deux $\F$-cônes affines sont équivalents lorsqu'ils sont parallèles et que leur intersection est non vide. L'intersection contient alors un autre $\F$-cône affine, parallèle aux deux premiers. On note $f\sim_{A_0} g$, ou juste $f\sim g$ lorsqu'il n'y a pas d'ambiguïté.

\end{itemize}\end{defins}

\begin{prop} La relation "être équivalent" est une relation d'équivalence sur $\F_{A_0}$.\end{prop}
Remarquons également que les relations "être un sous-cône" et "être un sous-cône parallèle" sont des relations d'ordre.\\

\begin{defin} On pose $\bar{A_0}= \F_{A_0}/\sim_{A_0}$. Pour $f\in \F_{A_0}$, on note $[f]_{A_0}$, ou juste $[f]$ lorsqu'aucune confusion n'est possible, la classe de $f$. Et si $x=[u+\vec f]_{A_0}$ avec $u\in A_0$, $\vec f\in \F$, on appelle $\vec f$ le cône directeur ou la direction de $x$ (il est uniquement déterminé).\\
 Pour $\vec f\in \F$, on note $A_{0 \vec f}$ l'ensemble des points de $\overline{A_0}$ de cône directeur $\vec f$, c'est la "façade" de type $\vec f$ de l'appartement $A_0$. La projection $A_0\rightarrow A_{0 \vec f}$ est notée $p_{A_0,\vec f}$ ou juste $p_{\vec f}$. Lorsqu'un cône $f$ de direction $\vec f$ est fixé, on pourra noter pour simplifier $A_f=A_{\vec f}$.
 \end{defin}

\begin{prop}
 Soit $\vec f\in \F_{A_0}$. La façade $A_{0 \vec f}$ est un espace affine isomorphe à $A_0/ \vect(\vec f)$, et son espace vectoriel directeur est $\fl{A_{0 \vec f}}\simeq \fl{A_0}/ \vect{\vec f}$.
\end{prop}

 Lorsqu'un cône vectoriel $\vec f\in \F$ s'écrit $\vec f=\{\alpha_i>0,\alpha_j=0,\, i\in I,\,j\in J\}$ avec $(\alpha_i)_{i\in I\sqcup J}$ une famille de formes linéaires, alors les $\alpha_j$, $j\in J$ s'identifient à des formes linéaires sur $\fl{A_{0 \vec f}}$, elles forment même une famille génératrice de $\fl{A_{0 \vec f}}^*$. Quand aux $\alpha_i$, $i\in I$, on sait qu'elles envoient tout représentant de tout point de $A_{0 \vec f}$ sur un voisinage de $+\infini$ dans $\R$, on dira donc qu'elles prennent la valeur $+\infini$ sur $A_{0 \vec f}$.\\

\subsection{Topologie sur $\overline{A_0}$}
\subsubsection{Définitions}

\begin{defin} Soit $\mathcal U$ l'ensemble des voisinages de $0$ dans $\fl{A_0}$ stables par le groupe de Weyl $W^v$. Pour $U$ un tel voisinage, et $g\in \F_{A_0}$, on pose:
\[\mathcal V_{A_0}(g,U)= \{ x\in \bar{A_0} \,|\, \exists h\in x  \tq  h\subset g+U\}\]
C'est l'ensemble des points de $\bar{A_0}$ ayant un représentant inclus dans $g+U$.\\
Lorsqu'il n'y aura pas d'ambiguïté, je noterai juste $\V(g,U)$.\end{defin}

\rema Comme le produit scalaire de $\fl{A_0}$ est $W^v$-invariant, les boules de centre $0$ sont des exemples d'éléments de $\U$.\\

\begin{prop} Il existe une unique topologie sur $\bar{A_0}$ telle que pour tout $x\in \bar{A_0}$, l'ensemble des $\mathcal V_{A_0}(g,U)$ pour $U\in \mathcal U_{A_0}$ et $g\in x$ soit une base de voisinage de $x$.\end{prop}

\demo\\
 Il suffit de vérifier que pour tout $x\in \bar{A_0}$, $U,V \in \mathcal U_{A_0}$, $g,g' \in x$, il existe $W\in \mathcal U_{A_0}$ et $h\in x$ tels que $\V_{A_0}(h,W)\subset \V_{A_0}(g,V)\cap \V_{A_0}(g',U)$. Mais comme $g\sim_{A_0} g'$, $g\cap g'$ contient un autre représentant $h$ de $x$. Soit $W=U\cap V$, alors $h+W\subset g+U \cap g'+V$, d'où $\V_{A_0}(h,W)\subset \V_{A_0}(g,V)\cap \V_{A_0}(g',U)$.\cqfd\\

 \rema La fonction $\V_{A_0}$ ainsi définie est croissante en les deux variables: si $U\subset V$ et $f\subset g$ alors $\V(f,U)\subset \V(f,V)$ et $\V(f,U)\subset \V(g,U)$.\\
 On remarque aussi qu'on peut remplacer $\U$ par n'importe quelle base de voisinages de $0$ dans $\fl{A_0}$, on obtiendra alors encore une base de la même topologie sur $\overline{A_0}$. On peut notamment choisir une base dénombrable de voisinages de $0$, ensuite si on se limite par exemple aux cônes dont les coordonnées du sommet sont rationnelles, on obtient une base dénombrable de la topologie de $\overline{A_0}$: $\overline{A_0}$ est à base dénombrable d'ouverts, les caractérisations par des suites sont possibles.\\
 
 Afin de pouvoir raisonner à l'aide de suites dans la suite, on va donner une caractérisation des suites convergentes. Commençons par le lemme suivant:
 \begin{lemme}
  Soit $\vec f\in \F$ et soit $\{\alpha_i=0,\, \alpha_j>0,\, i\in I,\,j\in J\}$ un système d'équations et d'inéquations déterminant $\vec f$. Soit $U\in \U$. Soit $(x_n)_n$ une suite dans $A_0$ telle que $\alpha_i(x_n)\rightarrow_{n\rightarrow \infini} 0$ pour $i\in I$, et $\lim_{n\rightarrow \infini}\alpha_j(x_n)=+\infini$ pour $j\in J$.\\
  Alors à partir d'un certain rang, les $x_n$ sont dans $U+\vec f$.
 \end{lemme}
  \pv\\
  La seule difficulté provient du fait que $\{\alpha_i\}_{i\in I\cup J}$ est une famille génératrice mais pas forcément libre de $A_0^*$. Soit $I'\subset I$ maximal tel que $\{\alpha_i\}_{i\in I'}$ soit libre. Soit ensuite $J'\subset J$ tel que $\{\alpha_i\}_{i\in I'\cup J'}$ soit une base de $A_O^*$.
  Il existe $\epsilon$ tel que $U':=\{x|\; |\alpha_i(x)|<\epsilon\, \forall i\in I'\cup J'\}\subset U$.\\
  Pour $n$ assez grand, on a $|\alpha_i(x_n)|<\epsilon$ pour $i\in I'$. Soit alors $y_n\in A_0$ tel que $\alpha_i(y_n)=\alpha_i(x_n)$ pour $i\in I'$ et $\alpha_j(y_n)=0$ pour $j\in J'$. Le point $y_n$ est alors dans $U$, et il suffit de vérifier que $x_n-y_n\in \vec f$.\\
  Or pour tout $i\in I'$, $\alpha_i(x_n-y_n)=0$. Comme tout $\alpha_i,i\in I$ est combinaison linéaire des $\alpha_j,j\in I'$, on obtient que $\alpha_i(x_n-y_n)=0$ pour tout $i\in I$. Enfin, $U'$ est borné, donc pour tout $j\in J$, $\alpha_j(y_n)$ aussi. Comme $\alpha_j(x_n)$ tend vers $\infini$, on vérifie bien que $\alpha_j(x_n)-\alpha_j(y_n)$ est positif dès que $n$ est assez grand.\cqfd\\

  Il est maintenant facile de prouver la caractérisation des suites convergentes suivante:\\
  
  \begin{prop}\label{prop:convergence dans A}
   Soit $(x_n)_{n\in \N}$ une suite de $\overline{A_0}$ dont tous les éléments ont la même direction $\vec f$. Soit $l=[a+\vec g]\in \overline{A_0}$. Alors $(x_n)_n$ converge vers $l$ si et seulement si les deux conditions suivantes sont réalisées:\\
   \begin{itemize}
   \item $\vec f$ est une face de $\vec g$.
   \item Pour toute famille $\{\alpha_i\}_{i\in I\sqcup J\sqcup K}\subset A_0^*$ de formes linéaires telle que $\vec g=\{x|\alpha_i(x)=0,\alpha_j(x)>0\, ,i\in I,\, j\in J\sqcup K\}$ et  $\vec f=\{x|\alpha_i(x)=0,\alpha_j(x)>0\, ,i\in I\sqcup J,\, j\in K\}$ on a $\lim_{n\rightarrow \infini}\alpha_i(x_n)=\alpha_i(a)$ pour $i\in I$ et $\lim_{n\rightarrow \infini}\alpha_j(x_n)=\infini$ pour $j\in J$. (Il existe une telle famille $\{\alpha_i\}$ car $\vec f$ est une face de $\vec g$).
   \end{itemize}\end{prop}
 
 \demo\\
  On commence par supposer que $(x_n)$ tend vers $l$. Si $\vec f$ n'est pas une face de $\vec g$, alors il existe une forme linéaire $\alpha$ telle que $\alpha(\vec f)\subset \R^+$ et $\alpha(\vec g)\subset \R^-$, et un vecteur $\vec v\in \vec f$ tel que $\alpha(\vec v)>0$. Donc $\alpha(\vec f)$ contient un voisinage de l'infini (dans $\R\cup\{\pm \infini\}$), et il en va de même de tous les $\alpha(h)$ pour $h$ un représentant d'un $x_n$. Il est donc impossible qu'un tel $h$ soit inclus dans un ensemble de la forme $U+a+\vec g$ dès que $U$ est borné, il est alors impossible que $x_n$ soit dans le voisinage $\V(a+\vec g,U)$ de $l$ correspondant. Et $\U$ contient bien des éléments bornés, prendre par exemple des boules.\\
  Soit ensuite un ensemble de formes linéaires $\{\alpha_i\}_{i\in I\sqcup J\sqcup K}\subset A_0^*$ tel que $\vec g=\{x|\alpha_i(x)=0,\alpha_j(x)>0\, ,i\in I,\, j\in J\sqcup K\}$ et  $\vec f=\{x|\alpha_i(x)=0,\alpha_j(x)>0\, ,i\in I\sqcup J,\, j\in K\}$.
   Soit $i\in I$, montrons que $\alpha_i(x_n)$ tend vers $\alpha_i(a)$. Soit $\epsilon>0$, il existe $U\in \U$ tel que $\alpha_i(U)\subset ]-\epsilon,\epsilon[$. Et pour $n$ assez grand, $x_n\in \V(a+\vec g,U)$ donc $x_n$ a un représentant inclus dans $U+a+\vec g$. Mais $\alpha_i(U+a+\vec g)\subset ]\alpha_i(a)-\epsilon,\alpha_i(a)+\epsilon[$, car $\alpha_i(\vec g)=\{0\}$. On en déduit que $|\alpha_i(a)-\alpha_i(x_n)|<\epsilon$. On a ainsi prouvé que $\lim_{n\rightarrow \infini}\alpha_i(x_n)=\alpha_i(a)$\\
  Soit $j\in J$, vérifions que $\alpha_j(x_n)$ tend vers l'infini. Soit $m\in \R$. Comme $\alpha_j(\vec g)$ est un voisinage de l'infini, il existe $\vec v\in \vec g$ tel que $\alpha_j(a+\vec v)>m+1$. Soit $U\in \U$ tel que $\alpha_j(U)\subset ]-1,1[$. Comme $a+\vec v +\vec g$ est un autre représentant de $l$, à partir d'un certain rang, $x_n$ doit être dans $\V(a+\vec v +\vec g,U)$. Et cela implique $\alpha_j(x_n)>m$.\\
  
  A présent, supposons que les deux conditions sont vérifiées par la suite $(x_n)$ et montrons qu'elle converge alors vers $l$.\\
  Comme $\alpha_k(\vec f)$ contient un voisinage de $\infini$ pour tout $k\in K$, on peut choisir $y_n$ un point d'un représentant de $x_n$ pour tout $n$, de sorte que les suites $\alpha_k(y_n)$ tendent vers $\infini$. Les suites $\alpha_i(y_n)$ pour $i\in I$ et $\alpha_j(y_n)$ pour $j\in J$ tendent quand à elle obligatoirement vers les $\alpha_i(a)$ et $\infini$, respectivement. Soit $U\in \U$, qu'on peut supposer ouvert, et soit $a'+\vec g$ un représentant de $l$, on peut supposer $a'\in a+\vec g$. On a pour tout $i\in I$, $\alpha_i(a)=\alpha_i(a')$ car $a$ et $a'$ diffèrent d'un élément de $\vec g$. Alors le lemme indique que $y_n$ est dans $U+a'+\vec g$ à partir d'un certain rang. Mais comme $y_n$ est un point d'un représentant de $x_n$, $y_n+\vec f$ est un représentant de $x_n$. Et $y_n\in U+a'+\vec g$, avec $U+a'+\vec g$ un ouvert, et $\vec f\subset \overline{\vec g}$ implique que $y_n+\vec f\subset U+a'+\vec g$. D'où $x_n\in\V(a'+\vec g,U)$.\\
  Ceci étant vrai pour tout $U\in \U$ ouvert, et pour tout $a'\in a+\vec g$, cela prouve la proposition.\cqfd\\

\subsubsection{Séparation}

\begin{prop} L'espace topologique $\bar{A_0}$ est séparé.\end{prop}

\demo\\
Soient $x$ et $y$ dans $\bar{A_0}$, $x\not= y$, et soient $f\in x$, $g\in y$ des représentants. Il faut trouver $U,V\in \mathcal U_{A_0}$, $f'\sim_{A_0}f$ et $g'\sim_{A_0} g$ tels que $\V_{A_0}(f',U)\cap \V_{A_0}(g',V)=\emptyset$. D'après la définition de $\V_{A_0}$, il suffit de s'assurer que $f'+U\cap g'+V=\emptyset$.\\

 Pour commencer, supposons que $\vec f=\vec g$, c'est-à-dire que $f$ et $g$ sont parallèles. Alors $f\cap g=\emptyset$ car sinon on aurait $f\sim_{A_0} g$ puis $x=y$. Soient $\{\alpha_i\}_{i\in I}$ et $\{\beta_j\}_{j\in J}$ des formes linéaires sur $A_0$ telles que $\vec f=\{ z| \forall i,\, \alpha_i(z)=0 $ et $\forall j, \,\beta_j(z)>0\}$.  Soit $\epsilon >0$ tel que $\exists i\in I$, $|\alpha_i(x)-\alpha_i(y)|>2\epsilon$, soit $B$ une boule centrée en $0$ telle que $\alpha_i(B)\subset ]-\epsilon,\epsilon[$. Alors $U=V=B$, $f'=f$ et $g'=g$ conviennent.\\
 
  Étudions maintenant le cas où $\vec f\not= \vec g$. Alors $\vec f\cap\vec g=\vide$. Les adhérences de ces deux cônes ne sont pas égales, supposons par exemple $\vec f\not\subset \bar{\vec g}$, soit $\vec u\in \vec f -\bar{\vec g}$. Pour tout $\lambda >0$, $f+\lambda \vec u\subset f$. Soit $\alpha$ une forme linéaire séparant $\vec f$ et $\vec g$: supposons par exemple $\alpha (\vec f)\subset \mathbb R^-$ et $\alpha (\vec g)\subset \mathbb R^+$. Alors $\alpha (\vec u)<0$, et pour $\lambda$ assez grand, $f+\lambda \vec u$ et $g$ sont séparés par un hyperplan $\alpha^{-1}(k)$, $k\in \mathbb R$. Quitte à choisir $\lambda$ encore un peu plus grand, on peut supposer que $\forall z\in f+\lambda \vec u$, $\alpha(z)\leq k-1$ et $\forall z\in g$, $\alpha(z)\geq k$. Soit alors $U\in \mathcal U_{A_0}$ tel que $\alpha(U)\subset ]-1/2,1/2[$, on a $f+\lambda \vec u +U\cap g+U=\vide$, ce qui achève la preuve.\cqfd\\

  \subsubsection{L'inclusion canonique}
  
  \begin{defin} Soit \fonc{i}{A_0}{\bar{A_0}}{x}{[\{x\}]_{A_0}}.
   Cette fonction est bien définie car $\{0\}\in \F$. Il est immédiat qu'elle est injective. On l'appelle l'injection canonique de $A_0$ dans $\bar{A_0}$.\end{defin}
  
  \begin{prop} L'injection canonique est continue, ouverte, d'image dense dans $\bar{A_0}$. De plus, pour $U\in \mathcal U$, et $f\in \F_{A_0}$, $i\inv (\V(f,U))=f+U$.\end{prop}
  \demo\\
  L'assertion $i\inv (\V(f,U))=f+U$ découle directement des définitions. On voit alors que si $U$ est ouvert, $i\inv (\V(f,U))$ est ouvert. Comme une base de la topologie de $\overline{A_0}$ est constituée des $\V(f,U)$ pour $U\in \U$ ouvert, ceci prouve que $i$ est continue.\\
  L'image est dense: si $f\in\F_{A_0}$, $U\in \U$ alors pour tout $x\in f+U$, $\{x\}\in i(A_0)\cap \V(f,U)$\\
  Enfin $i$ est ouverte car les ensembles de la forme $x+U$, $x\in A_0$, $U\in \U$ forment une base de voisinages de $A_0$ et ont pour image les $\V(\{x\},U)$.\cqfd\\
  
  On identifie donc au moyen de $i$ $A_0$ à un ouvert dense de $\bar A_0$. Une fois cette identification faite, on peut remarquer l'égalité suivante: $\forall U\in\U$, $f\in\F_{A_0}$, \[\overline{\V_{A_0}(f,U)}=\overline{f+U}\]

  \subsubsection{Compacité}
  \begin{prop} L'espace topologique $\overline{A_0}$ est compact.\end{prop}

  \demo\\
  Il reste à montrer que toute suite admet une valeur d'adhérence. Soit donc $(x_n)_{n\in \N}$ une suite dans $\overline{A_0}$. Comme $\F$ est fini, on peut supposer que tous les $x_n$ ont une même direction $\vec f$. Fixons une origine $0\in A_0$ et identifions $\fl{A_0}$ et $A_0$. Alors chaque $x_n$ a un représentant inclus dans un cône de $\F$. Comme $\F$ est fini, on peut supposer qu'il existe $F\in \F$ tel que chaque $x_n$ a un représentant dans $F$. On a forcément $\vec f\subset \bar F$.\\
  
  Choisissons un système d'inéquations pour $F$: 
\[F=\{x\in A_0\, | \alpha_i(x)>0 \, \forall i\in I, \, \alpha_i(x)=0\, \forall i\in J\}\]
  les $\{\alpha_i\}_{i\in I\cup J}$ étant une famille génératrice de $A_0^*$. Un système d'inéquations de $\vec f$ est alors de la forme:
  \[\vec f=\{x\in A_0\, | \alpha_i(x)>0 \, \forall i\in I_1, \, \alpha_i(x)=0\, \forall i\in J\cup I_2\}\]
  
  où $I=I_1\sqcup I_2$ est une partition de $I$. Pour $i\in J$, on a $\alpha_i(x_n)=0$ pour tout $n$ puisque les $x_n$ ont un représentant inclus dans $F$. Et pour $i\in I_1$, on a $\alpha_i(x_n)=\infini$.\\
  Quitte à extraire une sous-suite de $(x_n)$, on peut supposer que pour tout $i\in I_2$, soit $\alpha_i(x_n)$ converge vers une limite $\lambda_i$, soit $\alpha_i(x_n)$ tend vers $\infini$. Soit $C=\{i\in I_2|\lim_{n\infini}\alpha_i(x_n)=\lambda_i\}$ et $N=\{i\in I_2| \lim_{n\infini}\alpha_i(x_n)=\infini\}$ ($C$ pour "convergent" et $N$ pour "non convergent").\\
  
  Posons \[\vec g=\{x\in A_0\, | \alpha_i(x)>0 \, \forall i\in I_1\sqcup N, \, \alpha_i(x)=0\, \forall i\in J\sqcup C\}.\]
   Nous allons maintenant montrer qu'il existe $u\in A_0$ tel que $\forall i\in C$, $\alpha_i(u)=\lambda_i$, puis que $g:=u+\vec g\in \F_{A_0}$ puis enfin que la suite $\{x_n\}$ converge vers $[g]_{A_0}$. Nous aurons besoin du lemme suivant:\\
   
     \begin{lemme} \label{lemme:facette non vide} Soit $E$ un espace vectoriel de dimension finie, $(\alpha_i)_{i\in I\sqcup J}$ une famille génératrice finie dans $E^*$, $(x_n)_n$ une suite dans $E$ vérifiant $\forall i\in I$, $(\alpha_i(x_n))_n$ est bornée, et $\forall j\in J$, $\alpha_j(x_n)$ tend vers $\infini$. Alors il existe $z\in E$ tel que $\forall i\in I$, $\alpha_i(z)=0$ et $\forall j\in J$, $\alpha_j(z)>0$.\end{lemme}
  
  \pv\\
   On montre le lemme par récurrence sur $|J|$ le cardinal de $J$. On peut supposer que les coordonnées des $x_n$ selon $J$ sont toujours non nulles.\\
   Si $|J|=0$ alors $z=0$ convient.\\
   Si $|J|=1$, soit $j$ l'élément de $J$. Soit $z_n=\frac{x_n}{\alpha_j(x_n)}$. Alors pour $i\in I$, $\alpha_i(z_n)\rightarrow 0$, et $\alpha_j(z_n)=1$ pour tout $n$. Tous les $\alpha_i(z_n)$ convergent donc, et comme la famille des $\alpha_i$ est génératrice de $E^*$, la suite $z_n$ admet une limite dans $E$. Et cette limite vérifie clairement les conditions requises.\\
   Supposons maintenant $|J|=k>1$. Quitte à prendre une sous suite de $(x_n)_n$, il existe une énumération $J=\{\alpha_{j_1},...,\alpha_{j_k}\}$ telle que:\\
   \begin{itemize}
   \item $\forall n\in \N$, $\alpha_{j_1}(x_n)\leq...\leq\alpha_{j_k}(x_n)$.
   \item $\forall u\in \llbracket 1,k\rrbracket $, la suite $(\frac{\alpha_{j_u}(x_n)}{\alpha_{j_k}(x_n)})_n$ converge.\end{itemize}
   On pose alors $z_n=\frac{x_n}{\alpha_{j_k}(x_n)}$, la suite $(z_n)$ converge vers une limite $z_\infini\in E$. Cette limite vérifie $\alpha_i(z_\infini)=0$ pour tout $i\in I$, et elle a au moins une coordonnée ($\alpha_{j_k}$) strictement positive. Soit $J_1=\{j\in J|\alpha_j(z_\infini)>0\}$.\\
   En retirant aux $x_n$ un multiple convenable de $z_\infini$, et en prenant éventuellement encore une sous suite, on peut rendre bornées un ensemble non vide $\{\alpha_j\}_{j\in K}$ de nouvelles coordonnées, avec $K\subset J_1$,  sans changer les $\alpha_i(x_n) ,i\in I$, et en s'assurant que les $\alpha_j , j\in J_1-K$ qui ne sont pas devenue bornées tendent toujours vers l'infini. Alors par hypothèse de récurrence, il existe $z_1\in E$ qui annule les $\alpha_i$, $i\in I\sqcup K$ et tel que $\alpha_j(z_1)>0$ pour $j\in J-K$. Alors $z=z_1+z_\infini$ convient, ce qui prouve le lemme.\cqfd\\
  
  \textit{montrons l'existence d'un tel $u$}:\\
  
  Posons $\lambda_i=0$ pour $i\in J$. Si la famille $\{\alpha_i\}_{i\in J\sqcup C}$ est libre, l'existence de $u$ ne pose aucun problème. Soit $K\subset J\sqcup C$ tel que $\{\alpha_i\}_{i\in K}$ soit une sous famille libre maximale de $\{\alpha_i\}_{i\in J\sqcup C}$. On peut alors trouver $u$ tel que $\alpha_i(u)=\lambda_i$ pour $i\in K$. Si $i\in (J\sqcup C) \setminus K$ alors $\alpha_i$ est combinaison linéaire des $\alpha_j$, $j\in K$. Disons $\alpha_i=\sum_{j\in K}a_j\alpha_j$. En appliquant cette égalité aux $x_n$ et en passant à la limite, on obtient $\lambda_i=\sum_{j\in K}a_j\lambda_j$. Donc $\alpha_i(u)=\sum_{j\in K}a_j\alpha_j(u)=\sum_{j\in K}a_j\lambda_j=\lambda_i$. Ainsi, $u$ convient bien.\\
  
  \textit{montrons que $g\in \F$}:\\
  
  Pour montrer que $g=u+\vec g\in \F_{A_0}$, il faut montrer que $\vec g\in \F$. On a $\vec g\subset \bar F$, et $\vec g$ est l'intérieur de Vect$(\vec g)\cap \overline{F}$. Donc $\vec g$ est une des faces de $F$, d'après les hypothèses faites sur $\F$. En fait, la seule difficulté est de montrer que $\vec g\not=\vide$. On va pour ce faire construire directement un point de $\vec g$ à partir de la suite $(x_n)$.\\
  Pour tout $n\in \N$, on peut choisir $t_n\in A_0$ un point d'un représentant de $x_n$ de sorte que $t_n\in F$ et $\alpha_i(t_n) \rightarrow_{n\rightarrow \infini} \infini$ pour tout $i\in I_1$. On a alors $\alpha_i(t_n)\rightarrow \infini$ pour $i$ dans $I_1\sqcup N$, $\alpha_i(t_n)\rightarrow \lambda_i$ pour $i\in C$ et $\alpha_i(t_n)=0$ pour $i\in J$. Alors le lemme prouve l'existence d'un point de $\vec g$.\\

   \textit{Montrons que la suite $(x_n)_n$ converge vers $[u+\vec g]_{A_0}$:}\\
   
   Résumons: on dispose d'une famille génératrice $(\alpha_i)$ de $A_0^*$, indexée par $I_1\sqcup C\sqcup N\sqcup J$, avec:\\
   \begin{itemize}
   \item $\alpha_i(x_n)=0$ pour $i\in J$.
   \item $\alpha_i(x_n)$ tend vers $\lambda_i$ pour $i\in C$.
   \item $\alpha_i(x_n)=\infini$ pour $i\in I_1$.
   \item $\alpha_i(x_n)$ tend vers $\infini$ pour $i\in N$.
  \end{itemize}
  De plus $\vec g=\{x|\alpha_i(x)>0,\alpha_j(x)=0\, \forall i\in I_1\sqcup N, j\in J\sqcup C\}$ et $\vec f=\{x|\alpha_i(x)>0,\alpha_j(x)=0\, , \forall i\in I_1,j\in j\sqcup N\sqcup C\}$, et donc $\vec f\subset \overline{\vec g}$.\\
  Nous sommes donc exactement dans la situation de la proposition \ref{prop:convergence dans A}, et la suite $(x_n)_{n\in\N}$ tend vers $[u+\vec g]_{A_0}$.\cqfd\\

  \subsection{Structure du bord de $\bar A_0$}

  \subsubsection{Prolongement des automorphismes}
  
   Pour que les automorphismes de $A_0$ se prolongent en des homéomorphismes de $\overline{A_0}$, il est nécessaire de supposer que $\F$ est stable par le groupe de Weyl vectoriel $W^v$ (hypothèse (H7)), ce que nous ferons dorénavant. Un tel prolongement est forcément unique puisque $A_0$ est dense dans $\overline{A_0}$.
  
  \begin{prop} Soit $w\in W$, et soit \fonc{\tilde w}{\overline{A_0}}{\overline{A_0}} { \left[f\right]_{A_0} }{ [w(f)]_{A_0} }.
   Alors $\tilde w$ est un homéomorphisme bien défini de $\overline{A_0}$ qui prolonge $w$, c'est donc l'unique prolongement continu de $w$ à $\overline{A_0}$.\\
  \end{prop}
   
   \demo\\
   Si $f=x+\vec f\in\F_{A_0}$, alors $w(f)=w(x)+\vec w(\vec f)\in \F_{A_0}$ car $\vec w\in W^v$ préserve $\F$. Donc $W$ préserve $\F_{A_0}$. On voit également que $W$ préserve la relation d'équivalence $\sim_{A_0}$ donc $\tilde w$ est bien défini. Comme de plus $w$ est bijectif sur $\F_{A_0}$, on obtient que $\tilde w$ est bijectif.\\
   Soit $x=[f]_{A_0}\in \overline{A_0}$, et $\V(f,U)$ un voisinage de $x$, montrons que $\tilde w\inv(\V(f,U))$ contient un voisinage de $\tilde w\inv(x)=[w\inv(f)]$. On calcule:   $\tilde w\inv(\V(f,U))=\{[w\inv (g)]|g\subset f+U\}=\{[g]|w(g)\subset f+U\}=\{[g]|g\subset w\inv(f+U)\}$. Mais $w\inv(f+U)=w\inv(f)+U$ donc 
   $\tilde w\inv(\V(f,U))=\V(w\inv(f),U\}$, qui est un voisinage de $\tilde w\inv(x)$. Donc $\tilde w$ est continue.\\
   Comme $\overline{A_0}$ est compact (donc en particulier séparé), $\tilde w$ est automatiquement fermée, c'est donc un homéomorphisme.\cqfd\\
   
   Bien entendu, l'action de $W^v$ sur $\fl{A_0}$ s'étend elle aussi de la même manière en une action par homéomorphismes sur $\bigsqcup_{\vec f\in \F}\fl{A_{0,\vec f}}$. L'action affine et l'action vectorielle sont compatibles au sens suivant: si $w\in W$ envoie $A_{0 \vec f}$ sur $A_{0 \vec g}$, si $\vec w\in W^v$ est sa partie vectorielle, alors $\vec w$ envoie $\fl{A_{0 \vec f}}$ sur $\fl{A_{0 \vec g}}$ et $\vec w_{|\fl{A_{0 \vec f}}}$ est bien la partie vectorielle de $w_{|A_{0 \vec f}}$.\\

   \subsubsection{Coeurs de cônes}
   Pour faire le lien entre les cônes de $\F$ et les facettes de Weyl vectorielles, la notion du coeur d'un cône développée dans ce paragraphe est utile. Il s'agit d'attacher à un cône de $\F$ une facette vectorielle, ou une partie de facette vectorielle, qui le caractérise.\\
   
  \textit{rappel:} Dans un complexe de Coxeter abstrait, si $f$ est une facette, alors $f^*$ est le complexe de Coxeter formé de toutes les facettes $g$ supérieures à $f$, c'est-à-dire telles que $f\sqsubset g$.\\
 Pour définir la notion similaire dans le cadre de la réalisation géométrique d'un immeuble abstrait, nous prenons les conventions suivante: si $\vec\delta \subset \vec A$ est un cône inclus dans une facette de Weyl, on appelle \textit{étoile} de $\vec\delta$ dans $A$ et $\vec\delta^*_A$ (parfois juste $\vec\delta^*$) l'union de toutes les chambres fermées contenant $\vec\delta$. C'est un cône convexe fermé d'intérieur non vide. L'étoile de $\vec\delta$ a la propriété que tout mur coupant son intérieur doit contenir $\vec\delta$. De plus, si $\phi:A\rightarrow B$ est un isomorphisme alors $\vec\phi(\vec\delta)^*_B=\vec\phi(\vec\delta^*_A)$. \\
 Si $f=x+\vec f$ avec $\vec f$ un cône inclus dans une facette de Weyl, on pose $f^*_A=x+\vec f^*_A$. \\

   \begin{defin}
   Soit $\vec f\in \F$, soit $W^v_{\vec f}=\stab_{W^v}(\vec f)$ le stabilisateur de $\vec f$. On pose $\delta(\vec f):=\vec f^{W^v_{\vec f}}=\{x\in \vec f \; |\; \forall w\in W^v_{\vec f},\, w(x)=x\}=\vec f\cap \fix_{A_0}(W^v_{\vec f})$, c'est le coeur de $\vec f$.
   \end{defin}
   
   Voici quelques propriétés simples du coeur d'un cône. Remarquons qu'on obtient des informations non triviales sur les cônes de $\F$, en particulier sur leurs stabilisateurs dans $W^v$.\\
   
   \begin{prop}\label{prop:coeur vect} Soit $\vec f\in \F$. Alors:
   \begin{enumerate}
   \item $\forall w\in W^v$, $\delta(w. \vec f)=w.\delta(\vec f)$.
   \item $\delta(\vec f)$ est un cône vectoriel convexe non vide.
   \item Si $\vec f\not=\vec g\in\F$, alors $\delta(\vec f)\cap\delta(\vec g)=\vide$.
   \item $\delta(\vec f)$ est inclus dans l'intersection des murs coupant $\vec f$, et $W^v_{\vec f}$ contient les réflexions par rapport à ces murs.
    \item $W^v_{\vec f}=\fix_{W^v}(\delta(\vec f))$, c'est à dire qu'un élément de $W^v$ stabilise $\vec f$ ssi  il fixe $\delta(\vec f)$.
   \item $\delta(\vec f)$ est inclus dans une facette de Weyl.
  
   \item $\vec f\subset \delta(\vec f)^*$.
   \item (amélioration de 4.)  $\delta(\vec f)$ est égal à l'intersection de $\vec f$ et des murs coupant $\vec f$, et $W^v_{\vec f}$ est le sous groupe de Coxeter engendré par les réflexions selon ces murs.

	\item Si $\vec g\in\F$ et si $\vec f\subset \overline{\vec g}$, alors il existe une facette de Weyl $\vec h$ telle que $\delta(\vec f)\cup\delta(\vec g)\subset \overline{\vec h}$.
 %  \item Si $\vec f\subset \partial\vec g$, alors $\delta(\vec f)\subset\partial\delta(\vec g)$.
      
   \end{enumerate}
   \end{prop}
   
   \demo\\
   \begin{enumerate}
   \item Comme $W^v_{w\vec f}=w W^v_{\vec f} w\inv$, les points fixes de $W^v_{w.\vec f}$ sont $w.Fix_{A_0}(W^v_{\vec f})$. D'où $\delta(w.\vec f)=w.\vec f \cap w.\fix_{A_0}(W^v_{\vec f})=w.(\vec f\cap \fix_{A_0}(W^v_{\vec f})=w.\delta(\vec f)$.

   \item L'ensemble des point fixes d'une application linéaire est un cône vectoriel convexe. Et $\vec f$ aussi. Donc $\delta(\vec f)$ est une intersection de cônes vectoriels convexe, c'en est donc un aussi. Soit $z\in \vec f$. Comme $W^v_{\vec f}$ est fini, on peut définir $g$ comme le barycentre  de $\{ w(z)|w\in W^v_{\vec f}\}$. Par convexité de $\vec f$, $g\in \vec f$, et c'est un point fixe pour $W^v_{\vec f}$. Donc $\delta(\vec f)\not=\vide$.
   
   \item Ceci est clair car $\delta(\vec f)\subset \vec f$.
   
   \item Soit $\vec M$ un mur coupant $\vec f$, soit $\sigma$ la réflexion selon $\vec M$. Comme $\F$ est stable par $W^v$, $\sigma(\vec f)\in \F$. Mais $\sigma$ fixe au moins un point de $\vec f$, donc $\vec f\cap \sigma(\vec f)\not=\vide$, donc $\sigma(\vec f)=\vec f$ et $\sigma\in W^v_{\vec f}$. Alors $\delta(\vec f)\subset \fix_{\vec f}(\sigma)=\vec M\cap \vec f$.
   
     \item  Soit $w\in W^v$. Si $w\in W^v_{\vec f}$ alors $w$ fixe $\delta(\vec f)$ par définition. Réciproquement, si $w$ fixe $\delta(\vec f)$, alors $w$ fixe au moins un point de $\vec f$, donc $w(\vec f)=\vec f$ (on a déjà fait ce raisonnement).
   
   \item Il s'agit d'abord de prouver que pour tout mur vectoriel $\vec M$, $\delta(\vec f)$ est soit inclus dans $\vec M$, soit inclus dans un des deux demi-appartements ouverts délimités par $\vec M$. Mais c'est une conséquence de la convexité de $\delta (\vec f)$ et du quatrième point. Ensuite, d'après le point précédent, il suffit de vérifier que $\fix_{W^v}(\delta(\vec f))=\fix_{W^v}(\vec g)$, où $\vec g$ est la facette de Weyl contenant $\delta(\vec f)$. Ceci découle du fait que $W^v$ préserve l'ensemble des facettes de Weyl, et que si $w\in W^v$ stabilise une facette, alors il la fixe.
   
   \item Soit $C$ une chambre fermée coupant $\vec f$. Montrons que $C$ contient $\delta(\vec f)$. Supposons par l'absurde qu'il existe $y\in \delta(\vec f)\setminus C$. Soit $x\in C\cap \vec f$, alors $[x,y]\subset \vec f$. Soit $z$ tel que $C\cap [x,y]=[x,z]$. On a $z\not=y$ car $y\not\in C$. Alors $z$ est dans un mur $\vec M$ qui borde $C$ et qui ne contient par $[x,y]$. Mais $z\in \vec f$, donc $\vec M$ coupe $\vec f$, donc d'après le quatrième point, $\delta(\vec f)\subset \vec M$. En particulier, $y\in \vec M$, d'où $[x,y]\in \vec M$, ce qui est une contradiction.
   
   \item Notons $\vec g$ la facette de Weyl contenant $\delta(\vec f)$. Alors $W^v_{\vec f}=\fix_{W^v}(\vec g)$, c'est un sous groupe de Coxeter de $W^v$, engendré par les réflexions selon les murs contenant $\vec g$. Mais ces murs sont justement les murs coupant $\vec f$ (utiliser 4.), d'où la seconde partie de 8. On sait alors de plus que $\overline{\vec g}$ est l'ensemble des points fixes de $W^v_{\vec f}$, et aussi l'intersection des murs contenant $\vec g$. On obtient alors que $\delta(\vec f)=\overline{\vec g}\cap \vec f$ puis que $\delta(\vec f)$ et l'intersection de $\vec f$ et des murs coupant $\vec f$.

\item Si $\vec f\subset \overline{\vec g}$, alors $\delta(\vec f)\subset \overline{\vec g}\subset \delta(\vec g)^*$, d'où le résultat.\cqfd\\
   
  % \item $\vec f\subset \overline{\vec g}\implique W^v_{\vec g}\subset W^v_{\vec f}
  % \implique \fix_{\vec A}(W^v_{\vec f})\subset \fix_{\vec A}(W^v_{\vec g})$. Comme $\delta(\vec f)=\vec f\cap \fix_{\vec A}(W^v_{\vec f})$, et analoguement pour $g$, en utilisant l'hypothèse $\vec f\subset\partial \vec g$, on obtient bien $\delta(\vec f)\subset \partial{\delta(\vec g)}$.\cqfd\\
     
   \end{enumerate}

 Lorsque $\vec f\in\F$, $\vec f$ n'est pas forcement inclus dans une facette de Weyl. On peut néanmoins poser $\vec f^*=\delta(f)^*$. Ceci contient $\vec f$, par le point 7. de la proposition. On notera alors, si $f\in\F_{A_0}$, $f^*=s(f)+\vec f^*$.\\

\begin{figure}[h]
\includegraphics[width=5cm, height=5cm]{A2quelconque.eps}
\includegraphics[width=5cm, height=5cm]{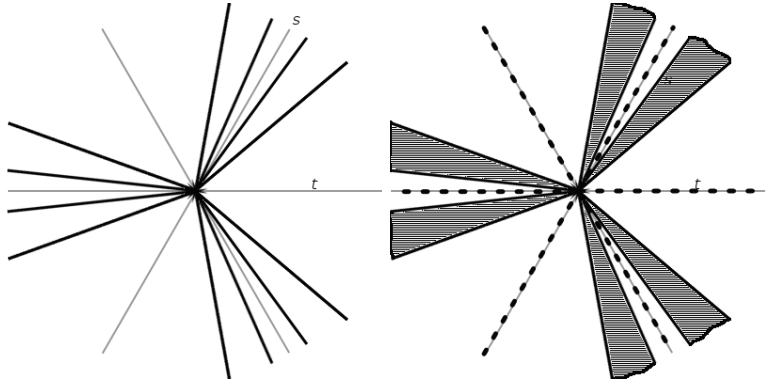}
\centering
\caption{Les coeurs d'un autre exemple de décomposition en cônes.}
\label{figure:coeurs}
\end{figure}

\textit{Exemples}\begin{itemize}
\item Dans le cas où $\F=\F^\vide$ est l'ensemble des facettes de Weyl de $\vec A_0$, on a $\delta(\vec f)=\vec f$, pour tout $\vec f\in \F$.
\item Dans le cas où $\F=\F^J$ comme dans \ref{subsection:exemple}, si $\vec f$ est une facette admissible de $\barre{\vec C_0}$ de type $I$, alors par définition, $W^v_{J\cap I^\perp}\subset \stab(J.\vec f)$. D'autre part, il est clair que $W^v_I=\stab(\vec f)\subset \stab J.\vec f$. Il est facile de vérifier qu'en fait $\stab(J.\vec f)=W^v_{I\cup(J\cap I^\perp}$. Dès lors, le coeur de $J.\vec f$ est la facette de $\barre{ \vec C_0}$ de type $I\cup(J\cap I^\perp)$. De manière plus générale, si $\vec f$ est une facette admissible, alors $\delta(J.\vec f)$ est la sous-facette de Weyl de $\vec f$ de type $I\cup (J\cap I^\perp)$.\\
\item La figure \ref{figure:coeurs} représente les coeurs des cônes de la décomposition du troisième exemple de \ref{subsection:dessins}. Les coeurs de dimension 2 sont hachurés, ceux de dimension 1 sont en pointillé. Pour faciliter la lecture, la figure de gauche rappelle cette décomposition en cônes.
\end{itemize}

   Pour finir, définissons le coeur d'un cône affine:\\
   \begin{defin} Soit $f=x+\vec f$ un cône affine. On définit alors $\delta(f)=x+\delta(\vec f)$, c'est le coeur de $f$.\end{defin}

   Le coeur d'un cône affine est toujours caractéristique du cône:
   \begin{prop} L'ensemble des coeurs de cône affine possibles est $\ens{x+\vec \delta}{x\in A_0,\; \vec\delta\in \delta(\F)}$. Et chacun de ces coeurs est le coeur d'exactement un cône affine.\end{prop}
   \demo Il est évident que chacun des $x+\vec \delta$, pour $x\in A_0$ et $\vec\delta$ le coeur d'un cône vectoriel est le coeur d'un cône affine. Étudions l'unicité. Supposons $\delta(f)=\delta(g)$ avec $f,g\in\F_{A_0}$. Comme $\delta(f)$ et $\delta(g)$ sont des sous-cônes de $f$ et $g$, ils ne contiennent pas de sous-espace vectoriel non trivial et ont donc chacun un unique sommet, tout comme $f$ et $g$. On peut alors conclure, d'après la définition du coeur d'un cône, que le sommet de $f$ est aussi celui de $\delta(f)$, de même pour $g$. Comme $\delta(f)=\delta(g)$, on obtient que tous ces sommets sont égaux. D'où $\delta(\vec f)=\delta(\vec g)$, puis $\vec f=\vec g$ d'après l'étude du cas vectoriel, d'où enfin $f=g$.\cqfd\\
   
   \begin{defin} Soit $\delta\in\delta(\F_{A_0})$ un coeur de cône dans $A_0$. On appelle cône affine de $A_0$ engendré par $\delta$ l'unique $f\in\F_{A_0}$ tel que $\delta(f)=\delta$.\end{defin}
   
   Voici quelques propriétés immédiates du coeur d'un cône affine:
   \begin{prop}Soit $f\in \F_{A_0}$. Alors:\\
\begin{itemize} 
   \item Si $\delta$ est un coeur tel que $\vec\delta\subset \barre{\delta(\vec f)}$ et $\delta\subset f$, alors le cône engendré par $\delta$ est inclus dans $f$. Si $\vec\delta \subset \overline{\delta(\vec f)}$ et $\delta \subset \bar f$, alors le cône engendré par $\delta$ est inclus dans $\bar f$.
   \item Pour tout $w\in W$, $w(f)=f$ si et seulement si $w$ fixe $\delta(f)$.
   \end{itemize}
   \end{prop}

      % On voit donc que $w$ stabilise $f$ équivaut à $w(x)=x$ et $\vec w$ stabilise $\vec f$. Et d'autre part $w$ fixe $\delta(f)$ équivaut à $w$ fixe $x$ et $\vec w$ fixe $\delta(\vec f)$. Le résultat découle alors de la proposition précédente qui traite le cas vectoriel.\cqfd\\
   
   \subsubsection{Structure de complexe de Coxeter vectoriel sur une façade}
   \label{subsubsection:cx cx vectoriel de la facade}
    Fixons un cône $\vec f\in \F$. Dans ce paragraphe, on notera $F:=A_{0 \vec f}$ la façade de $A_0$ de type $f$. C'est un espace affine dont l'espace directeur est $\vec F= \fl{A_0}/\vect(\vec f)$. On munit $\vec F$ du produit scalaire obtenu en l'identifiant avec $(\vec f)^\perp$. Nous commençons par étudier la structure de  $\vec F$. Notons $p=p_{A_0,\vec f}$ la projection sur $F$, et $\vec p:\fl{A_0}\rightarrow \vec F$ sa partie vectorielle.\\

     Soit $\mathcal M^v$ l'ensemble (fini) des murs vectoriels contenant $\vec f$. Si $\vec M\in\mathcal M^v$, alors $p(\vec M)$ est un hyperplan de $\vec F$. Nous dirons que les $p(\vec M)$, pour $\vec M\in \mathcal M^v$ sont les \textit{murs} de $\vec F$. La fonction \fonc{p_{\mathcal M^v}}{\M^v}{\{\text{murs de }\vec F\}}{\vec M}{p(\vec M)} est injective, nous permettant d'identifier $\M^v$ à l'ensemble des murs de $\vec F$. Cet ensemble de murs fait de $\vec F$ un complexe de chambres, nous voulons prouver qu'il s'agit en fait d'un complexe de Coxeter. Il nous reste donc à trouver un groupe de Coxeter agissant simplement transitivement sur l'ensemble des chambres de $\vec F$.\\
     
     Rappelons que le groupe $W^v$ agit sur $\cup_{\vec f\in \F}\vec A_{0 \vec f}$ et que le stabilisateur de $\vec f$ est $W^v_{\vec f}=\{w\in W^v|w(\vec f)=\vec f\}$. Pour $w\in W^v_{\vec f}$, on note $p(w)$ l'application induite sur $\vec F$. Le groupe $p(W^v_{\vec f})$ préserve $\mathcal M^v$ et le produit scalaire de $\vec F$, c'est le candidat naturel comme groupe de Coxeter de $\vec F$. Montrons que c'est effectivement un groupe de Coxeter et qu'il agit simplement transitivement sur l'ensemble des chambres de $\vec F$. On observe la suite exacte courte:

     \[0\rightarrow \fix_{W^v_{\vec f}}(\vec f^\perp) \rightarrow W^v_{\vec f} \,\stackrel{p}{\rightarrow}\, p(W^v_{\vec f})\,\rightarrow 0\]
     
     D'après le paragraphe précédent, $W^v_{\vec f}$ est le fixateur dans $W^v$ de $\delta(\vec f)$, ou encore le fixateur de la facette de Weyl contenant $\delta(\vec f)$. Choisissons une chambre $C$ contenant $\delta(\vec f)$ dans son adhérence, soit $S$ le système générateur de $W^v$ formé des réflexions par rapport aux murs de $C$. Alors il existe $I\subset S$ tel que $W^v_{\vec f}=W_I=<s| s\in I>$. De plus, $I=\{s\in S|s(\delta(\vec f))=\delta(\vec f)\}$.\\

     Comme $W^v$ préserve le produit scalaire sur $\fl{A_0}$, les éléments de $W^v_{\vec f}$ stabilisent $\vect(\vec f)$ et $(\vec f)^\perp$. Dans une base adaptée, leurs matrices sont du type $Mat( w)= \begin{bmatrix}\psi(w) & 0\\
     						0 & \phi(w)\\ \end{bmatrix}$, où $\psi(w)$ est la matrice de $w$ restreint à $\vect(\vec f)$ et $\phi(w)$ est la matrice de $w$ restreint à $(\vec f)^\perp$, c'est aussi la matrice de $p(w)$.\\
     Soit $s\in I$, $s$ est diagonalisable et préserve les espaces $\vect (\vec f)$ et $(\vec f)^\perp$. On peut donc choisir une base adaptée comme précédemment en imposant en plus que $\psi(s)$ et $\phi(s)$ soient diagonales. Le spectre de $s$ est composé d'un seul $-1$ et les autres valeurs propres sont $1$, donc soit $\psi(s)$ est une réflexion et $\phi(s)=id$, soit $\phi(s)$ est une réflexion et $\psi(s)=id$.\\
 Notons $I_1$ l'ensemble des réflexions de $I$ se trouvant dans le premier cas, c'est-à-dire qui agissent sur $\vec f$ et fixent $\vec f^\perp$. Notons $I_2$ l'ensemble des réflexions se trouvant dans le second cas, qui fixent $\vec f$ et donnent une réflexion sur $\vec f^\perp$.\\
 Alors $W_{I_2}\subset \fix_{W^v}(\vec f)$, $W_{I_1}\cap W_{I_2}=\{id\}$, et $p_{|W_{I_2}}$ est injective. De plus, tout élément de $W_{I_1}$ commute avec tout élément de $W_{I_2}$ et $W^v_{\vec f}=W_I=W_{I_1\sqcup I_2}$ est engendré par $W_{I_1}$ et $W_{I_2}$.\\
 Par conséquent, $W^v_{\vec f}\simeq W_{I_1}\times W_{I_2}$, avec $W_{I_2}\simeq \fix_{W^v}(\vec f)\simeq p(W^v_{\vec f})$ ($I_1$ et $I_2$ sont deux parties disjointes du diagramme de Coxeter de $\vec A_0$). En particulier, $p(W^v_{\vec f})$ est un groupe de Coxeter.\\
     
     A présent, on vérifie que le complexe de chambres défini par $\M^v$ sur $\vec F$ est isomorphe au complexe de Coxeter $W_{I_2}.C$. On dispose d'une fonction $W_{I_2}$-équivariante de $W_{I_2}.C$ dans les chambres de $\vec F$: il s'agit de la fonction $p_{\mathcal C}$ qui a une chambre $D$ dans $W_{I_2}.C$ associe la chambre de $\vec F$ contenant $p(D)$.\\
     Cette fonction est injective car les murs de $W_{I_2}.C$ contiennent tous les points fixes de $W_{I_2}$ donc en particulier $\vec f$, donc ce sont tous des éléments de $\M$. Deux chambres distinctes dans $W_{I_2}.C$ sont donc séparées par un mur de $\M^v$, elles donnent donc deux chambres distinctes dans $\vec F$.\\
     Enfin vérifions que $p_{\mathcal C}$ est surjective. Toutes les réflexions par rapport à un mur de $\M^v$ sont des éléments de $W_{I_2}$. Comme ces dernières agissent transitivement sur les chambre de $\vec F$ et comme $p_{\mathcal C}$ est $W_{I_2}$-invariante, on obtient bien la surjectivité.\\
     
     Résumons:\\
     \begin{prop}\label{prop:cx cx de facade}
     Soit $\vec f\in \F$, alors l'espace vectoriel directeur de la façade $A_{0\vec f}$ est muni d'une structure de complexe de Coxeter où:\begin{itemize}
     \item les murs sont les hyperplans du type $p(\vec M)$ où $\vec M$ est un mur de $\vec A_0$ contenant $\vec f$.
     \item le groupe de Coxeter est $p(W^v_{\vec f})\simeq \fix_{W^v}(\vec f)$ (isomorphe à $W_{I_2}$ avec les notations précédentes).\end{itemize}
     Si $C$ est une chambre fermée contenant $\delta(\vec f)$, alors ce complexe de Coxeter est isomorphe à $\fix_{W^v}(\vec f).C$, via l'application $p_{\mathcal C}:w.C\mapsto $la chambre de $\vec F$ contenant $p(w.C)$.\\
     \end{prop}
     
\rema Ce complexe de Coxeter n'est pas forcément essentiel. Le cas extrême est atteint lorsque $\vec f$ est une droite, incluse dans une chambre. Alors $I=\vide$, donc $A_{0 \vec f}$ est trivial comme complexe de Coxeter, et pourtant il est de dimension $\dim(\I)-1$.\\
En fait, $A_{0\vec f}$ est essentiel $\equi\: \fix_{A_0}(W_{I_2})=\vect(\vec f)$ $\equi\: \vect(\vec f)=\vect(\cl(\vec f))$.\\
En particulier, $A_{0\vec f}$ est essentiel si $\vec f$ est une réunion de facettes de Weyl.\\

    Dans la décomposition $W^v_{\vec f}\simeq W_{I_1}\times W_{I_2}$, le premier facteur est donc le groupe de Coxeter de $\vec F$. Concernant le deuxième facteur, on montre facilement le résultat suivant:
\begin{prop}\label{prop:cx cx de facade 2}
 Soit $\vec E$ le supplémentaire orthogonal de $\vect(\delta(\vec f))$ dans $\vect(\vec f)$, de sorte qu'en identifiant $\vec F$ avec $\vec f^\perp$, on a $\vec A_0=\vect(\delta(\vec f))\oplus \vec E\oplus \vec F$, et la somme est orthogonale.\\
 Alors les murs de $\vec A_0$ contenant $\delta(\vec f)$ contiennent soit $\vec E$, soit $\vec F$. On sait déjà que les premiers induisent les murs de $\vec F$, les second définissent un complexe de Coxeter sur $\vec E$, dont le groupe de Coxeter est $W_{I_1}$. Ce complexe de Coxeter est essentiel.\\
\end{prop}

 En terme de diagrammes de Dynkin, disons que si $D$ est le diagramme de $A_0$, $I$ (qu'on identifie à un sous-diagramme de $D$) est obtenu en enlevant les noeuds de $D$ correspondant au type de la facette de Weyl contenant $\delta(\vec f)$. Le diagramme obtenu est séparé en deux parties non reliées, le digramme $I_1$ de $\vec E$ et le diagramme $I_2$ de $\vec F$.\\
 Notons que si le diagramme de $A_{\delta(\vec f)}$ est connexe, alors $\vec E$ ou $\vec F$ doit être trivial, c'est-à-dire sans mur, ce qui correspond aux situations suivantes:\\
\begin{itemize}
\item $\vec E$ trivial $\equi\: I_1=\vide \equi \stab_{W^v}(\vec f)=\stab_{W^v}(\delta(\vec f))\equi \vec f\subset \cl(\delta(\vec f))$. Comme $\delta(\vec f)=\cl(\delta(\vec f))\cap \vec f$, ceci équivaut en fait à $\vec f=\delta(f)$ donc à $\vec E=\{0\}$.
\item $\vec F$ trivial $\equi\: I_2=\vide\equi \fix_{W^v}(\vec f)=\{e\}\equi \cl(\vec f)$ contient une chambre.\end{itemize}

   \subsubsection{Structure de complexe de Coxeter affine sur une façade}
   
   A présent, soit $\M=\{M \;\text{mur de}\; A_0\; |\; \vec M\in \M^v\}$. Alors $\M$ est un ensemble discret d'hyperplans.\\
 Le sous-groupe de $W$ qui stabilise la façade $F=A_{0 \vec f}$ est $W_{\vec f} = \{ w\in W| \vec w\in W^v_{\vec f}\}=\{w\in W | \vec w(\vec f)=\vec f\}$. Ce sous-groupe contient les symétries par rapports aux murs de $\M$, il préserve la structure euclidienne de $F$ et cet ensemble d'hyperplans. En notant $p(w)$ l'élément de $Isom(F)$ induit par $w$, on voit donc que $p(W_{\vec f})$ agit sur $F$ en préservant $\M$, et qu'il contient les réflexions orthogonales par rapport aux hyperplans de $\M$.\\
   Ceci suffit à prouver que $\M$ définit un complexe de Coxeter affine sur $F$, et que $p(W_{\vec f})$ contient le groupe de Coxeter correspondant.\\
   \begin{defin} L'ensemble des facettes de $\bar A_0$ est la réunion des ensembles de facettes de chacune de ses façades.\end{defin}

      \subsection{Compactification de chaque appartement}
   
    Soit $A$ un appartement de $\I$ soit $\phi:A_0\rightarrow A$ un isomorphisme et $\vec \phi:\fl{A_0}\rightarrow \vec A$ sa partie linéaire.\\
     On notera encore $\F$ l'ensemble des $\vec \phi(\vec f)$, $\vec f\in \F$. Comme deux isomorphismes entre $A_0$ et $A$ diffèrent d'un élément du groupe de Weyl $W$, le choix de $\phi$ n'intervient pas dans cette définition. De la même manière, on notera encore $\U$ l'ensemble des $\vec\phi(U)$, $U\in \U$, ou si l'on préfère, l'ensemble des voisinage de $0$ dans $\vec A$ stables par $W(\vec A)$.\\
 Comme dans la partie précédente, on définit $\F_A$ l'ensemble des cônes affines de $A$, $\bar A$ le compactifié de $A$, $\ens{A_{\vec f}}{\vec f\in \F}$ l'ensemble des façades de $\bar A$, $\V_A$ la fonction qui à un cône de $\F_A$ et un élément de $\U$ associe un voisinage dans $\bar A$, le coeur d'un cône, etc...\\

Un isomorphisme $\psi:A\isom B$ entre deux appartements induit un isomorphisme entre chacun de ces objets pour $A$ et l'objet correspondant pour $B$. Il s'étend notamment en un unique homéomorphisme de $\bar{A}$ sur $\bar B$, noté encore $\psi$. Explicitement, on a $\phi( [f]_{A_0})=[\phi(f)]_A$.\\
 Pour toute façade $A_{\vec f}$ de $A$, $\psi$ induit une isométrie entre $A_{\vec f}$ et $B_{\vec\psi(\vec f)}$, cet homéomorphisme induit une bijection entre les murs de $A_{\vec f}$ et ceux de $B_{\vec\psi(\vec f)}$, c'est donc un isomorphisme de complexes de Coxeter.\\
 Dans la suite, on appellera isomorphismes d'appartements compactifiés les applications telles $\phi$ qui sont des homéomorphismes entre deux appartements compactifiés qui induisent sur chaque façade un isomorphisme de complexe de Coxeter.

   \section{Quelques résultats généraux sur les immeubles}
\label{section:resultats generaux}

    Le but principal de cette section est de donner des critères permettant de s'assurer qu'une partie de $\I$ est incluse dans un appartement. Les résultats de \ref{subsection:parties closes} sont vrais dans un immeuble quelconque, nous cesserons donc pour cette partie de supposer $\I$ affine.\\
 A défaut d'une réalisation affine, il sera toujours possible d'utiliser la réalisation comme cône de Tits d'un complexe de Coxeter $\Sigma$ quelconque. Il s'agit d'un cône convexe $\mathcal C$ dans un espace vectoriel de dimension finie $V$ dont le sommet est $0$. Les murs de $\Sigma$  correspondent à des hyperplans vectoriels de $V$ (qui rencontrent l'intérieur de $\mathcal C$), les facettes sont des cônes convexes de sommet $0$, les facettes de dimension $i$ sont des cônes ouverts dans un sous-espace vectoriel de dimension $i$ (voir \cite{bourbaki} 4.6).\\

  La notation $\langle f \rangle$ pour une facette $f$ dans un complexe de Coxeter $\Sigma$ représentera l'intersection des murs de $\Sigma$ qui contiennent $f$. Si $\Sigma$ est un complexe de Coxeter affine, identifié à sa représentation affine, alors $\langle f \rangle=\aff(f)$, et si $\Sigma$ est un complexe de Coxeter quelconque identifié à son cône de Tits, alors $\langle f \rangle=\vect(f)$.\\

   On rappelle d'abord quelques résultats classiques, sur lesquels s'appuie toute la suite:\\
 
   \begin{prop} Soit $\I$ un immeuble, de système de Coxeter $(W,S)$. Son système complet d'appartement est l'ensemble des parties de $\I$ isomorphes au complexe de Coxeter de $(W,S)$.\end{prop}
   
   A partir de cette proposition, on montre le résultat fondamental suivant (voir \cite{ronan}, théorème 3.6 page 31):
   \begin{prout}  \label{th:fonda} Soit $Z$ un ensemble de chambre de $\I$, isométrique pour la $W$-distance de $\I$ à une partie de $\Sigma(W,S)$. Alors il existe un appartement du système complet d'appartements de $\I$ contenant $Z$.\end{prout}

   On note immédiatement deux conditions équivalentes à celle donnée par le théorème:
   \begin{cor}\label{cor:du th fonda}
   Soit $Z$ un ensemble de chambres dans $\I$. Les conditions suivantes sont équivalentes:\\
   \begin{itemize}
   \item Pour tout $c,d,e\in Z$, $\delta(c,e)=\delta(c,d)\delta(d,e)$, où $\delta$ est la $W$-distance de $\I$.
   \item Trois chambres quelconques de $Z$ sont incluses dans un même appartement.
   \item Il existe un appartement contenant $Z$.\end{itemize}
   \end{cor}

   Voici deux situations particulières où s'applique ce théorème:
   \begin{cor}\label{cor:inclusion dans un appartement}
   \begin{itemize}
   \item Si $Z$ est une galerie tendue (pas forcément finie), alors $Z$ est contenu dans un appartement du système complet d'appartements de $\I$.
   \item Si $D$ est un demi-appartement délimité par un mur $M$, et si $c$ est une chambre de $\I$ dont une cloison est incluse dans $M$, alors $D\cup c$ est inclus dans un appartement du système complet d'appartements de $\I$.\end{itemize}
   \end{cor}

Enfin, voici une version assez générale du résultat parfois appelé « lemme fondamental des immeubles »:
\begin{lemme}\label{lemme:fonda des immeubles}
 Soient $\Sigma$ et $\Sigma'$ des complexes de chambres. On suppose que dans $\Sigma'$, toute cloison est dans au plus deux chambres. Soit $c$ une chambre de $\Sigma$ et $d$ une chambre de $\Sigma'$. Alors il existe au plus un morphisme de complexes de chambres injectif de $\Sigma$ dans $\Sigma'$ qui envoie $c$ sur $d$.\\
En particulier, si $\Sigma\subset\Sigma '$ et $c=d$, un tel morphisme est forcément l'inclusion.

\end{lemme}

      \subsection{Parties closes dans un appartement}\label{subsection:parties closes}
   
   Voici quelques rappels sur la notion de partie close, principalement issus de \cite{bruhat-tits} pour le cas affine. Pour le cas général considéré ici, on peut se reporter à \cite{remy} 5.4.3. On donne une preuve de deux résultats classiques, à savoir la clôture de l'intersection de deux appartements et l'existence d'un isomorphisme entre deux appartements fixant leur intersection dans le cas général où l'intersection des appartements ne contient pas forcément de chambre.\\

%\subsubsection{Généralités}

   \begin{defin} \label{def:clos1}(rappel) Une partie $Z$ dans un appartement $A$ est dite close si c'est une intersection de demi-appartements de $A$. Si $Z$ est une partie quelconque de $A$, on note $\text{Cl}_A(Z)$ (ou juste $\cl(Z)$) et on appelle enclos de $Z$ la plus petite partie close contenant $Z$. C'est aussi l'intersection de tous les demi-appartements de $A$ contenant $Z$.
   \end{defin}
   
   \begin{prop} Si $Z$ est la fermeture d'un ensemble $E$ de chambres, alors $Z$ est clos si et seulement si $E$ contient toutes les galeries minimales entre deux éléments de $E$.\end{prop}
   
   \demo\\
   Le sens $\implique$ est clair. Pour l'autre sens, il s'agit de montrer que si $c$ est une chambre de $A$ qui n'est pas dans $E$, alors il existe un mur séparant $c$ de toutes les chambres de $E$.\\
   Soit $d_1,...,d_k,c$ une galerie minimale de $E$ à $c$. Alors $d_1\in E$ et $d_2\not\in E$. Soit $M$ le mur entre $d_1$ et $d_2$. Soit $e\in E$, si $M$ séparait $d_1$ et $e$ alors il existerait une galerie minimale entre $e$ et $d_1$ passant par $d_2$, ce qui impliquerait que $d_2\in E$, ce qui est impossible. Donc $M$ sépare $c$ et $E$.\cqfd\\
   
Cette proposition permet une définition de la clôture dans $\I$, pour des ensembles de chambres, cohérente avec la précédente:
\begin{defin} \label{def:clos2} Un ensemble de chambres de $\I$ est dit clos s'il contient toutes les galeries minimales entre deux de ses chambres. L'enclos d'un ensemble de chambres est le plus petit ensemble de chambres clos le contenant.\end{defin}
La proposition précédente prouve que si $Z$ est un ensemble de chambres inclus dans un appartement $A$, alors $Z$ est clos (au sens \ref{def:clos2}) si et seulement si $\bar Z$ est clos dans $A$ (au sens \ref{def:clos1}).\\

 L'enclos d'une partie, même convexe, d'un appartement n'est pas toujours évidente. Par exemple dans un appartement de type $\tilde B_2$ on trouve facilement des parties convexes de toute dimension non nulle, dont l'enclos contient des chambres qui n'ont aucun point en commun avec la partie de départ.\\

\subsubsection{Projection}

Pour étudier la clôture d'ensembles de facettes qui ne sont pas forcement l'adhérence d'un ensemble de chambre, la notion de projection est très utile. Rappelons-en les principales propriétés.

\begin{defin} Soit $\Sigma$ un complexe de chambres, soient $c$ et $d$ deux simplexes de $\Sigma$. On appelle projection de $d$ sur $c$ et on note $\proj_c(d)$ l'intersection des chambres finales de toutes les galeries minimales de $d$ à $c$.\end{defin}

\begin{prop} \label{prop:projections} Soit $\Sigma$ un complexe de Coxeter, et $c$, $d$ deux simplexes de $\Sigma$. Alors:
\begin{itemize}
\item $\proj_c(d)\subset \cl(c\cup d)$
\item Soit $g$ la chambre terminale d'une galerie minimale de $d$ à $c$, et soit $\M$ l'ensemble des murs contenant $c$ et $d$. Alors $\proj_c(d)$ est l'intersection de $g$ et des murs de $\M$.
\item L'étoile de $\proj_c(d)$ est formée des chambres terminales des galeries tendues de $d$ à $c$ et des intersections de ces chambres. Son groupe de Coxeter, c'est-à-dire $\fix_{W(\Sigma)}(\proj_c(d))$ est le groupe engendré par les réflexions selon les murs de $\M$.
\item  $\dim(proj_c(d)) \geq \dim (d)$ avec égalité si et seulement si $c\subset \langle d \rangle_\Sigma$.
\end{itemize}
\end{prop}

\rema La dimension d'une facette $f$ dans un complexe de Coxeter $\Sigma$ est son nombre de sommet. Ceci coïncide avec la dimension de $\vect(f)$ dans le cône de Tits, et avec la dimension de $\aff(f)$ plus $1$ dans une réalisation affine si $\Sigma$ est de type affine et irréductible. Dans le cas non irréductible, la représentation affine de $\Sigma$ contient des polysimplexes qui engendrent un espace de dimension moindre que leur nombre de sommets moins $1$.\\
 Cependant, la proposition est encore vraie dans ces cas si on remplace la dimension d'une facette par la dimension de l'espace affine qu'elle engendre.\\

\demo\\
Pour le premier point, soit $M^+$ un demi-appartement délimité par le mur $M$ et contenant $c$ et $d$. Soit $\Gamma$ une galerie minimale de $d$ à $c$, si $\Gamma$ n'est pas incluse dans $M^+$, alors en la pliant le long de $M$ on obtient une autre galerie minimale de $d$ à $c$ qui elle est incluse dans $M^+$. La chambre finale de cette galerie contient $\proj_c(d)$, d'où $\proj_c(d)\subset M^+$. Ceci prouve que $\proj_c(d)\subset \cl(c\cup d)$.\\

Passons au second point. Soit $\Gamma=g_0,...,g_k$ une galerie tendue de $d$ à $c$ avec $g_k=g$ et $k=d(c,d)$. Soit $M\in \M$, notons $\sigma_M$ la réflexion selon $M$. Alors $\sigma_M(\Gamma)$ est aussi une galerie minimale de $d$ à $c$, donc $\proj_c(d)\subset g\wedge\sigma_M(g)=g\wedge M$. Nous prouvons ainsi que $\proj_c(d)\subset g\wedge \bigwedge_{M\in \M} M$.\\
Pour montrer l'inclusion réciproque, il s'agit de prouver que toute galerie tendue de $d$ à $c$ se termine par une chambre contenant $g\wedge \bigwedge_{M\in M} M$. Soit $\Theta=h_0,...,h_k$ une telle galerie. Il suffit de prouver que les murs séparant $g$ de $h_k$ sont dans $\M$. Si $N$ est un tel mur, il contient $g\wedge h_k$ donc en particulier $c$. De plus, il ne peut couper ni $\Gamma$ ni $\Theta$ sans quoi on pourrait réduire ces galeries par un pli le long de $N$. Donc $N$ sépare $g_0$ de $h_0$, donc contient $g_0\wedge h_0$, et en particulier $d$. Ceci prouve que $N\in\M$, et conclut la démonstration du second point. \\

 Notons $\mathcal T$ l'ensemble des chambres terminales de galeries tendues de $d$ à $c$. Notons aussi $W'$ le groupe engendré par les réflexions selon un mur de $\M$. Nous venons de voir que tout mur séparant deux chambres de $\mathcal T$ est dans $\M$. Nous avons vu juste avant que $W'$ stabilise $\mathcal T$. Ceci entraîne que $\mathcal T=W'.g$, et que $\mathcal T$ contient toutes les galeries tendues entre deux de ses éléments. De plus, l'action de $W'$ sur $\mathcal T$ est simplement transitive puisque celle de $W(\Sigma)$ sur les chambres de $\Sigma$ l'est.
Soit $S'$ l'ensemble des réflexions selon les murs de $\M$ bordant $\mathcal T$. Si $w\in W'$, soit $\Lambda$ une galerie tendue de $g$ à $wg$, elle est incluse dans $\mathcal T$. Le premier mur qu'elle rencontre correspond à une réflexion $s\in S'$, et alors $l(sw)<l(w)$. On prouve ainsi par récurrence que $W'$ est engendré par $S'$, c'est donc le sous-groupe parabolique de $W(\Sigma)$ correspondant à la chambre $g$ et aux réflexions $S'$. C'est donc le fixateur de la facette de $g$ fixée par $S'$. Mais $\fix_\Sigma(S')=\fix_\Sigma(W')=\bigcap_{w\in W'}wg$, c'est l'intersection des chambres de $\mathcal T$, et par définition, c'est $\proj_c(d)$. Au final, $W'=\fix_{W(\Sigma)}(\proj_c(d))$, et $\mathcal T=W'.g$ est l'ensemble des chambres de $(\proj_c(d))^*$. Ceci prouve le troisième point.\\

 Pour le dernier point, on considère la réalisation géométrique $\mathcal C\subset V$ de $\Sigma$ comme cône de Tits.\\
 Soit $\M_g$ l'ensemble des murs de $\M$ qui bordent $g$, ils correspondent à une famille d'hyperplans linéairement indépendants (au sens où les formes linéaires les définissant sont indépendantes). Nous avons vu que $\proj_c(d)$ est l'intersection de $g$ et des murs de $\M_g$, d'où $\codim(\proj_c(d))=\card(\M_g)=\codim(\bigcap_{M\in\M_g}M)$. Comme tous ces murs contiennent $d$, on a $d\subset\bigcap_{M\in\M_g}M$, d'où $\dim d\leq \dim ( \bigcap_{M\in\M_g}M)=\dim(\proj_c(d)$.\\
 Le cas d'égalité est atteint lorsque tous les murs contenant $d$ contiennent également $\proj_c(d)$. On a alors $c\sqsubset \proj_c(d) \subset \langle d\rangle_\Sigma$. Réciproquement, si $c \subset \langle d\rangle_\Sigma$, alors tous les murs contenant $d$ contiennent aussi $c$ et donc sont dans $\M$. Au vu du point précédent, $\M$ est l'ensemble des murs contenant $\proj_c(d)$. On a donc $\langle d\rangle_\Sigma=\langle \proj_c(d)\rangle_\Sigma$ d'où $\dim d=\dim(\proj_c(d))$.\cqfd\\

\subsubsection{Les parties closes sont des complexes de chambre}

\begin{prop} \label{prop:clos implique cx de chb} Une partie close dans un complexe de Coxeter est un complexe de chambre.\end{prop}

\demo\\
Soit $E$ une partie close dans un complexe de Coxeter $\Sigma$.\\
Soit $c$ un simplexe de $E$, choisissons $d$ un simplexe de dimension maximale dans $E$. La proposition \ref{prop:projections} montre que $\proj_c(d)$ est inclus dans $E$ et est de dimension supérieure à celle de $d$. Comme $d$ est de dimension maximale, les dimensions sont en fait égale, ceci prouve que $c$ est inclus dans un simplexe de dimension maximale.\\

Il reste à prouver qu'entre deux simplexes de $E$ de dimension maximale existe une galerie dans $E$. Soient $c$ et $d$ deux tels simplexes, on procède par récurrence sur la distance, dans $\Sigma$, entre $c$ et $d$.\\
Si $d(c,d)=0$, cela signifie que $c$ et $d$ sont inclus dans une même chambre $C$ de $\Sigma$. Soit alors $c\vee d$ la plus petite facette de $C$ contenant $c$ et $d$, elle est incluse dans $\cl(c,d)$ donc dans $E$. Comme $c$ et $d$ sont de dimension maximales dans $E$, on a $\dim(c\vee d)=\dim(c)=\dim(d)$, d'où $c=c\vee d=d$. La galerie de longueur $0$ formée uniquement de la chambre $c$ relie $c$ à $d$.\\
Si $d(c,d)>0$, soit $\Gamma=g_0... g_k$, $k\geq 1$, une galerie tendue dans $\Sigma$ de $d$ à $c$. Comme $c \not\sqsubset g_{k-1}$, le simplexe $\sigma:=c\wedge g_{k-1}$ est de codimension $1$ dans $c$ (c'est une cloison de $E$). Maintenant $\proj_\sigma(d)$ est une chambre de $E$ adjacente à $c$ et strictement plus proche dans $\Sigma$ de $d$ que $c$. Ce qui conclut la récurrence.\cqfd\\

Le cas particulier d'une intersection de murs est intéressant, puisqu'on prouve alors qu'il s'agit d'un complexe de chambre "au plus mince":
\begin{lemme} \label{lemme:mur mince} Une intersection de murs dans un complexe de Coxeter est un complexe de chambres dans lequel une cloison est incluse dans au plus deux chambres.\end{lemme}
\demo\\
Soit $\M$ un ensemble de murs dans un complexe de Coxeter $\Sigma$, et $E$ leur intersection. En particulier, c'est une partie close de $\Sigma$ donc un complexe de chambres.\\
Soit $\sigma$ une cloison de $E$. On considère la réalisation géométrique $\mathcal C$ comme cône de Tits du complexe de Coxeter $\sigma^*$ formé des simplexes de $\Sigma$ contenant $\sigma$.
 %C'est un cône convexe dans un espace vectoriel $V$ partitionné par des cônes vectoriels convexes correspondant aux simplexes de $\sigma^*$. Les murs sont l'intersection d'hyperplans de $V$ avec $\mathcal C$.
  La facette minimale de $\sigma^*$ est $\sigma$, elle correspond à $\{0\}$ dans le cône de Tits. Les chambres de $E$ contenant $\sigma$ correspondent à des demi-droites issues de $0$. Nous voulons montrer qu'il ne peut y avoir qu'au plus deux telles demi-droites.\\
 Les demi-appartements de $\Sigma$ contenant $E$ soit contiennent $\sigma^*$ soit contiennent $\sigma$ dans le mur qui les bordent et induisent alors un demi-appartement de $\sigma^*$. Ceci prouve que $E\cap\sigma^*$ est égale à l'intersection des demi-appartements de $\sigma^*$ induits par les demi-appartements de $\Sigma$ contenant $E\cap\sigma^*$ et dont le bord contient $\sigma$. Donc $E\cap \sigma^*$ est une partie close de $\sigma^*$. Dans $V$, c'est l'intersection de $\mathcal C$ avec des demi-espaces, c'est donc un convexe. Donc si $d$ est une chambre de $E$ contenant $\sigma$, alors $E\cap\sigma^*\subset \vect(d)$, en identifiant $E\cap\sigma^*$ et $d$ à leurs images dans $V$. Ainsi les chambres de $E$ contenant $\sigma$ sont, dans $V$, des demi-droites issues de $0$ incluse dans la droite $\vect(d)$: il n'y en a que deux possibles.\cqfd\\

 Ce dernier résultat permet de prouver la caractérisation géométrique des parties closes suivante:

   \begin{prop}\label{prop:convexe et ferme implique clos} 
Soit $\Sigma$ un complexe de Coxeter, $\mathcal C\subset V$ sa réalisation comme cône de Tits dans l'espace vectoriel $V$. Soit $E$ un ensemble de facettes dont l'image dans $\mathcal C$ est fermée et convexe. Alors $E$ est une partie close.
%Si $\I$ est un immeuble affine, $E$ une partie convexe et fermée dans un appartement $A$ de $\I$, et si $E$ est une union de facettes de $A$, alors $E$ est close.
\end{prop}
   \rema\begin{itemize}
\item Réciproquement, une partie close de $\Sigma$ est clairement convexe et fermée dans $\mathcal C$, et égale à une union de facettes.
\item Ce résultat est encore vrai dans le cadre d'un complexe de Coxeter affine en remplaçant la réalisation comme cône de Tits par la réalisation comme espace affine (et en remplaçant  dans la preuve $\vect(f)$ par $\aff(f)$ ).\\
\end{itemize}

   \pv\\
On identifie $E$ à son image dans $\mathcal C$.
   Soit $f\subset E$ une facette de dimension maximale. Soit $x\in f$. Pour tout $y\in E$, le segment $[x,y]$ est contenu dans $E$. Comme $f$ est une facette maximale de $E$, ce segment doit rester, au voisinage de $x$, dans $f$. Donc $[x,y]\subset\vect(f)$. Ceci prouve que $E\subset\vect(f)$.\\
   A présent, soit $g\subset \vect(f)$ une facette qui n'est pas incluse dans $E$, montrons qu'il existe un mur $M$ tel que $E|^M g$.
 Soit $x\in f$, $y\in g$. Comme $E$ est fermée, il existe $z\in[x,y]$ tel que $[x,y]\cap E=[x,z]$. Quitte à déplacer le point $x$ dans $f$, on peut supposer que $z$ est dans une facette $h\subset \vect(f)$ de dimension dim$(f)-1$. Il existe un mur $M$ tel que $h=M\cap\vect(f)$, vérifions que $M$ convient. S'il existait un point $a\in E$ tel que $x|^M a$, alors le segment $[a,z]$ serait inclus dans $E$. Mais il existe au plus deux facettes de dimension $\dim(f)$ incluses dans $\vect(f)$ et contenant $z$ dans leur adhérence, par le lemme \ref{lemme:mur mince}. Or $[x,z]$ et $[a,z]$ coupent deux telles facettes distinctes. Ceci prouve qu'il y en a deux, et qu'elles sont dans $E$. Mais ceci contredit la définition de $z$ car l'adhérence de leur réunion contient un voisinage de $z$ dans $\vect(f)$ inclus dans $E$.\cqfd\\

\subsubsection{Intersection de deux appartements}

   \begin{prop}\label{prop:intersection close} L'intersection de deux appartements $A$ et $B$ est une partie close de $A$ et de $B$.\\
%De plus, si $f$ est une facette maximale de $A\cap B$, alors $A\cap B$ est inclus dans l'intersection des murs de $A$ ou de $B$ contenant $f$.
\end{prop}

   \demo\\
    Soit $f$ une facette maximale dans $A\cap B$, c'est à dire une facette qui n'est incluse dans l'adhérence d'aucune autre facette incluse dans $A\cap B$. Montrons que $A\cap B\subset \langle f\rangle_A$.\\
    Soit $g$ une facette dans $A\cap B$. Soient $\Gamma_A=a_1,...,a_k$ et $\Gamma_B=b_1,...,b_k$ deux galeries minimales de $f$ à $g$, l'une dans $A$ l'autre dans $B$ (deux galeries minimales entre deux facettes ont forcément la même longueur). Soit $\Delta=(a_1=d_1,...,d_l=b_1)$ une galerie minimale de $a_1$ à $b_1$ et $Z$ un appartement la contenant. On a $a_1\wedge b_1=f$ puisque $f$ est une facette maximale dans $A\cap B$, et toutes les chambres de $\Delta$ contiennent $f$.\\
    Soit $\rho=\rho_{A,a_1}$, $\rho(\Gamma_B)$ est une galerie minimale entre $f$ et $g$ dans $A$. Si $\phi$ est l'isomorphisme de $Z$ sur $A$ fixant $a_1$, les murs séparant $a_1$ et $\rho(b_1)$ sont les images par $\phi$ des murs de $Z$ séparant $a_1$ et $b_1$. Si $M$ est un tel mur, $M$ ne sépare pas $f$ et $g$ puisque $f\subset M$. Donc $g\subset M$, sans quoi une des galeries $\Gamma_A$ ou $\rho(\Gamma_B)$ traverserait $M$ et ne serait donc pas minimale. Il reste donc à prouver que l'intersection de ces murs est $\langle f\rangle_A$, ou de manière équivalente, que l'intersection des murs de $Z$ traversés par $\Delta$ est $\langle f\rangle_Z$.\\
    
    Soit $S$ le système générateur du groupe de Coxeter $W$ de $Z$ formé par les réflexions selon les cloisons de $a_1$. Soit $I\subset S$ le type de $f$, $W_I=<I>$ le fixateur de $f$ dans $W$. Soient $M_1,..,M_l$ les murs traversés par $\Delta$, $\sigma_1,...,\sigma_l$ les réflexions correspondantes, et $U$ le groupe qu'elles engendrent. Il suffit de montrer que $U=W_I$. Pour commencer, $\sigma_1\in I$ puisque $M_1$ est un mur de $a_1$ et que $\sigma_1$ fixe $f$. Notons $s_1=\sigma_1$. Ensuite, on a $\sigma_2=\sigma_1\circ s_2\circ\sigma_1$, avec $s_2\in S$. On constate alors que $s_2$ fixe $f$. En continuant ainsi jusqu'à $\sigma_l$, on constate que $U=<s_1,s_2,s_3,...,s_l>$ où les $s_i$ sont des éléments de $S$ qui fixent $f$, c'est à dire des éléments de $I$. Il ne reste plus qu'à prouver que chaque élément de $I$ apparaît dans $\{s_1,...,s_l\}$. Soit $s\in I$, supposons par l'absurde que $s\not\in\{s_1,...,s_l\}$. Soit $h$ la facette de type $I-\{s\}$, c'est une facette strictement plus grande que $f$, et elle est fixée par toutes les $s_i$. Elle est donc fixée par toutes les $\sigma_i$, mais ceci entraîne que $h\subset a_1\wedge b_1$, et cela contredit l'hypothèse de maximalité sur $f$.\\
    
    Nous avons ainsi prouvé que $A\cap B\subset \langle f\rangle_A$. On peut conclure par récurrence sur la longueur $l$ de $\Delta$.\\
    Si $l=0$, alors $f$ est une chambre. Il est alors bien connu que dans ce cas, $A\cap B$ est l'adhérence d'un ensemble de chambres, et que toute galerie minimale entre deux chambres de $A\cap B$ est incluse dans $A\cap B$. Donc $A\cap B$ est clos.\\
    Si $l>0$, utilisons le corollaire \ref{cor:inclusion dans un appartement} pour trouver un appartement $A_1$ contenant un demi-appartement de $A$ délimité par $d_1\wedge d_2$ ainsi que $d_1$ et $d_2$. Cet appartement contient le mur de $A$ définit par $d_1\wedge d_2$, il contient donc en particulier $\langle f\rangle_A$, donc $A\cap B$. Donc $A\cap B=(A\cap A_1)\cap (B\cap A_1)$. Par hypothèse de récurrence, $A\cap A_1$ et $A_1\cap B$ sont deux parties closes de $A_1$, donc $(A\cap A_1)\cap (B\cap A_1)$ est close dans $A_1$. Mais il existe un isomorphisme de $A_1$ sur $A$ qui fixe $A\cap A_1$, donc en particulier $(A\cap A_1)\cap (B\cap A_1)$, ce qui conclut la preuve.\cqfd\\

\subsubsection{Isomorphisme entre deux appartements}

%Nous pouvons enfin montrer qu'entre deux appartements $A$ et $B$ existe un isomorphisme fixant l'intersection. 

\begin{prop} \label{prop:isom fixant l intersection}Soient $A$ et $B$ deux appartements, alors il existe un isomorphisme $\phi:A\isom B$ fixant $A\cap B$. De plus, si $\psi:A\isom B$ est un isomorphisme fixant une facette de dimension maximale de $A\cap B$, alors $\psi$ fixe $A\cap B$.\end{prop} 
 
\demo\\
 Par \ref{prop:intersection close} et \ref{prop:clos implique cx de chb}, $A\cap B$ est un complexe de chambre. Soit $d$ une chambre de $A\cap B$. Soient $C$ et $D$ des chambres de $A$ et $B$, respectivement, contenant $d$. Soit $Z$ un appartement contenant $C\cup D$. En composant un isomorphisme de $A$ sur $Z$ fixant $C$ puis un isomorphisme de $Z$ sur $B$ fixant $D$, on obtient un isomorphisme $\phi:A\isom B$ fixant $d$. Montrons que $\phi$ fixe $A\cap B$.\\
Le dernier point de la proposition \ref{prop:projections} montre que $A\cap B\subset \langle d\rangle_A$ et $A\cap B\subset \langle d\rangle_B$. Or $\phi(\langle d\rangle_A)=\langle d\rangle_B$, d'où $\phi(A\cap B)\subset \langle d\rangle_B$. Comme $\langle d\rangle_B$ est un complexe de chambre "au plus mince", par \ref{lemme:mur mince}, et que $\phi$ induit un morphisme injectif entre $A\cap B$ et  $\langle d\rangle_B$ le "lemme fondamental des immeubles" \ref{lemme:fonda des immeubles} montre que $\phi$ fixe $A\cap B$.\cqfd\\

   \subsection{Cheminées}
 La notion de cheminée permettra d'utiliser pleinement le théorème \ref{th:deux galeries dans appart} (ci-dessous). Nous en donnons ici une définition ainsi que les propriétés élémentaires. On suppose dorénavant que $\I$ est affine.\\

\begin{lemme}\label{lemme:direction de cheminee} Soit $A$ un appartement, $c$ une facette de $A$, $\vec f$ une facette de $\vec A$. Alors les points de $\cl(c+\vec f)$ restent à distance bornée de $c+\vec f$.\\
De plus, si $d$ est une autre facette de $A$, $\vec g$ une autre facette de $\vec A$ telles que $\cl(d+\vec g)\subset \cl(c+\vec f)$, alors $\vec g\subset \barre{\vec f}$.\end{lemme}

\pv\\
 Soit $\{\alpha_i\}_{i\in I\sqcup J}$ des formes affines sur $A$ telles que $\{\vec\alpha_i\}_{i\in I\sqcup J}$ soit une base du système de racines de $\vec A$ avec $\vec f=\{ \vec\alpha_i>0,\,: \vec\alpha_j=0,\: i\in I,\; j\in J\}$ et telles que les murs de $A$ de direction $\ker \alpha_i$ sont les $\alpha_i\inv(\{n\})$ pour $n\in \Z$. Alors:
\[ c+\vec f\subset \{ \alpha_i>inf_c(\alpha_i),\: \alpha_j\in \alpha_j(c),\: i\in I,\; j\in J\}
\subset \{ \alpha_i\geq n_i,\,: \alpha_j\in [n_j,m_j],\: i\in I,\; j\in J\} \]
en appelant, pour $i\in I\sqcup J$, $n_i$ la partie entière de $\inf_c(\alpha_i)$ et pour $j\in J$, $m_j$ le plus petit entier supérieur à $\sup_c(\alpha_j)$.\\
Le troisième ensemble est une intersection de demi-appartements donc une partie close, donc contient $\cl(c+\vec f)$. Il est de plus égal à $\barre{b+\vec f}$ avec $b=\{\alpha_i=n_i,\: \alpha_j\in [n_j,m_j] \}$. Comme $\{\vec\alpha_i\}_{i\in I\sqcup J}$ est une base de $\vec A^*$ et qu'elle reste bornée sur $b$, on déduit que $b$ est borné. Mais $c$ est également borné, on voit alors facilement que les points de $\barre{b+\vec f}$ sont à distance bornée de $c+\vec f$. Et ceci reste en particulier vrai pour les points de $\cl(c+\vec f)$.\\

Maintenant si $\vec g\not\subset \barre{\vec f}$, cela signifie qu'il existe $\vec v\in \vec g\setminus \barre{\vec f}$. Alors $d(\vec v,\vec f)>0$. Pour $x$ un point de $d$, les points de la forme $x+\lambda \vec v$ sont dans $\cl(d+\vec g)$ et peuvent être rendus arbitrairement éloignés de $c+\vec f$, car $c$ est borné. Ceci contredit que $\cl(d+\vec g)\subset \cl(c+\vec f)$.\cqfd\\

   \begin{defin}
 \begin{itemize}

   \item   Une cheminée dans un appartement $A$ est une partie de $A$ de la forme $R=\cl(c+\vec f)$, où $c$ est une facette de $A$ et $\vec f$ une facette de $\vec A$. La facette $\vec f$ est uniquement déterminée par $R$ d'après le lemme \ref{lemme:direction de cheminee}, on l'appelle la direction de $R$.\\

   \item Soient $R$ et $R'$ deux cheminées. On dira que $R'$ est une sous-cheminée de $R$ lorsque $R'\subset R$. On dira que $R'$ est une sous-cheminée pleine, abrégé en scp, lorsque de plus $R$ et $R'$ ont la même direction et $\aff(R)=\aff(R')$.\\

   \end{itemize}
   \end{defin}
   
On abrège de la même manière sous-cheminée pleine et sous-cône parallèle, mais les deux notions ont à peu près le même sens, l'une dans le cadre des cônes affines et l'autre dans le cadre des cheminées.\\
%Il est immédiat que les relations "être une sous-cheminée" et "être une sous-cheminée pleine" sont des relations d'ordre. On remarque également que si $R'$ est une sous-cheminée de $R$, $R''$ une sous-cheminée de $R'$ qui est également une scp de $R$, alors $R'$ est une scp de $R$.\\

   \begin{prop} \label{prop:sous cheminee}Soit $R=\cl(c+\vec f)$ une cheminée dans $A$, soit $d$ une facette de $A$ incluse dans $R$. Alors $\cl(d+\vec f)\subset R$, c'est donc une sous-cheminée de $R$.\end{prop}
   \demo\\
    Soit $D$ un demi-appartement de $A$ contenant $R$, montrons qu'il contient $d+\vec f$. Ceci n'est pas tout à fait évident car $d$, et donc $d+\vec f$ n'est pas forcément incluse dans $c+\vec f$. Soit $\alpha\in A^*$, $k\in\Z$ tels que $\D=\D(\alpha,k)=\{x\in A|\alpha(x)\geq k\}$. Nous savons déjà que $d\subset R\subset \D$. Soit $\vec v\in \vec f$, $x\in d$, montrons que $x+\vec v\in \D$. Dans le cas contraire, nous aurions $\alpha(x)\geq k$ et $\alpha(x+\vec v)<k$, d'où $\vec\alpha(\vec v)<0$. Choisissons alors un point quelconque $y\in c$, on aura $\alpha(c)\geq k$ car $c\subset \D$, mais $\alpha(c+\lambda\vec v)<0$ pour $\lambda$ un réel assez grand. ceci contredit le fait que $c+\vec f\subset \D$. Donc $x+\vec v\in \D$.\cqfd\\
    
   Le lemme suivant permet dans certains cas de se ramener à des cheminées dont la base est une chambre.
\begin{lemme}\label{lemme: cheminee de base une chambre}
Soit $R=\cl(x+\vec f)$ une cheminée. Soit $x$ une chambre dont l'adhérence contient $c$, et $R_1=\cl(x+\vec f)$. Alors toute sous-cheminée pleine de $R_1$ contient une sous-cheminée pleine de $R$.
\end{lemme}
\pv
Soit $R_1'=\cl(y+\vec f)$ une scp de $R_1$. Le cône $y+\vec f^*$ contient une scp $R'$ de $R$, et pour montrer que $R'\subset R_1'$, il suffit de prouver qu'aucun mur dont la direction contient $\vec f$ ne sépare strictement $R$ de $R_1'$. Soit $M$ un tel mur. Comme $R$ et $R_1'$ sont deux parties incluse dans $R_1=\cl(x+\vec f)$, si $M$ sépare strictement $R$ de $R_1'$, alors $M$ doit séparer strictement deux parties de $x+\vec f$. Mais c'est impossible pour un mur dont la direction contient $\vec f$.\cqfd\\

   \begin{defin} Soit $R=\cl(c+\vec f)$ une cheminée. Une demi-droite caractéristique de $R$ est une demi-droite $\delta$ dont l'extrémité est dans $c$ et dont la direction est incluse dans $\vec f$ (en particulier, $\delta\subset R$).\end{defin}
   \rema Si $\delta$ est une demi-droite, alors $\cl(\delta)$ est une cheminée dont $\delta$ est caractéristique.\\

   \begin{prop} \label{prop:demi droite caracteristique} Soit $R=\cl(c+\vec f)$ une cheminée dans un appartement $A$. Soit $\delta$ une demi-droite caractéristique de $R$. Alors:\\
   \begin{itemize}
   \item $R=\cl(\delta)$
	\item $\aff(\cl(\delta))$ est l'intersection des murs contenant $\delta$
   \item Si $\delta'$ est une demi-droite incluse dans $\delta$, alors $\cl(\delta')$ est une scp de $R$.\end{itemize}
   \end{prop}
   
   \demo\\
    Pour le premier point, la proposition \ref{prop:sous cheminee} prouve que $\cl(\delta)\subset R$. Pour montrer l'autre inclusion, soit $\D=\D(\alpha,k)$ un demi-appartement contenant $\delta$, montrons qu'il contient $R$. Soit $x$ le sommet de $\delta$, $\vec \delta$ sa direction, de sorte que $\delta=x+\vec\delta$. Le fait que $x$ est dans $c$ prouve que $\D$ contient $c$. Et le fait que $\D$ contient une demi droite dirigée par $\vec\delta$ prouve que $\alpha(\vec\delta)\subset [0,\infini[$. Mais $\vec\delta$ est inclus dans $\vec f$, et $\alpha$ est une racine, d'où $\alpha(\vec f)\subset [0,\infini[$. On obtient alors $c+\vec f\subset \D$, d'où $R\subset \D$.\\
    
    Pour prouver le second point, soit $\M_\delta$ l'ensemble des murs contenant $\delta$. Soit $M\in\M_\delta$, alors $\delta$ est dans les deux demi-appartements délimités par $M$, donc $\cl(\delta)$ est dans l'intersection de ces deux demi-appartements, c'est-à-dire dans $M$. Ceci prouve que $\aff(\cl(\delta))\subset \bigcap_{M\in \M_\delta} M$.\\
    Pour l'inclusion réciproque, comme la répartition des murs de $A$ est discrète, il existe $x\in\delta$ qui n'est inclus dans aucun autre mur que ceux de $\M_\delta$. Soit $\epsilon$ la distance minimale entre $x$ et un mur n'appartenant pas à $\M_\delta$, on a $\epsilon>0$ encore par la répartition discrète des murs. Alors la boule dans $\bigcap_{M\in \M_\delta} M$ de centre $x$ et de rayon $\epsilon$ est incluse dans $\cl(\delta)$, ce qui prouve que 
   $\bigcap_{M\in \M_\delta} M\subset \aff(\cl(\delta))$.\\

	Enfin pour le troisième point, il faut montrer que $\aff(\cl(\delta))=\aff(\cl(\delta'))$. Ceci provient directement du point précédent.\cqfd\\
   
   \rema Guy Rousseau définit les cheminées comme l'enclos dans $A$ d'une demi-droite, et d'un germe de segment de même origine. Dans le cas d'un immeuble où les murs sont répartis de manière discrète, c'est une définition équivalente au vu de la proposition qui précède.\\
   
   Ainsi, une cheminée peut-être caractérisée par une demi-droite. Nous allons maintenant voir qu'une demi-droite est incluse dans une galerie tendue, ce qui permettra de ramener les raisonnements sur les cheminées à des raisonnement sur les galeries tendues.\\
   
   \begin{prop} \label{prop:galerie sur droite} Soit $\delta$ une demi-droite et $s$ son origine, ou $\delta$ un segment et $s$ une extrémité de $\delta$. Alors pour toute chambre $c_0$ contenant $s$ dans son adhérence, il existe une galerie tendue commençant par $c_0$ contenant $\delta$ dans son adhérence. De plus, l'adhérence de toute chambre ou cloison de cette galerie rencontre $\delta$.\end{prop}
   \pv\\
    Soit $\mathcal M$ l'ensemble fini des murs contenant $\delta$. Soit $E$ l'intersection de tous les demi-appartements $\D(M,c_0)$ pour $M\in\mathcal M$. La demi-droite $\delta$ est incluse dans $E$. Elle est découpée en segments d'intérieur non vide (pour la topologie induite sur $\delta$), et chacun de ces segments est inclus dans une unique chambre fermée incluse dans $E$. En effet, si un tel segment est inclus dans deux chambres, alors il existe un mur $M$ séparant ces chambres et contenant ce segment. Comme le segment est d'intérieur non vide, $M$ contient $\delta$, donc $M\in\mathcal M$, donc au plus une des deux chambres est dans $E$. Soient $(I_k)_{k\in\N}$ ces segments, rangés de manière que le sommet de $\delta$ soit dans $I_0$ et que $\forall k$, $I_k\cap I_{k+1}\not= \vide$. Soient $(c_k)_{k\in\N}$ les chambres correspondantes, $c_0$ est bien la chambre donnée dans l'énoncé. Pour $k\in \N$, soit $\Gamma_k$ une galerie minimale entre $c_k$ et $c_{k+1}$. C'est une galerie incluse dans $E$ car $E$ est clos. On définit enfin $\Gamma$ comme la concaténation des $\Gamma_k$ auxquels on a retiré leurs dernières chambres, c'est à dire $c_{k+1}$, pour éviter qu'elles n'apparaissent deux fois de suite. Cette galerie contient bien $\delta$, il reste à vérifier qu'elle est tendue.\\
    Remarquons pour commencer que puisque $\Gamma\subset E$, les murs de $\mathcal M$ ne séparent pas les chambres de $\Gamma$. Supposons qu'il existe un mur $M$ coupant deux fois $\Gamma$, entre $d_1$ et $d_2$ puis entre $e_1$ et $e_2$. Il existe $k$ et $l$ tels que $d_1\cup d_2 \subset\Gamma_k$ et $e_1\cup e_2\subset \Gamma_l$. Toutes les chambres de $\Gamma_k$ contiennent dans leur adhérence le point de $I_k\cap I_{k+1}$, donc $M$ contient aussi ce point. De même $M$ contient le point de $I_l\cap I_{l+1}$. Mais $M$ ne peut contenir deux points distincts de $\Delta$ sans contenir $\Delta$, donc $k=l$. Alors $M$ coupe $\Gamma_k$ deux fois, mais c'est impossible puisque $\Gamma_k$ est une galerie minimale.\cqfd\\

   \subsection{Inclusion d'une partie d'appartement et d'une chambre dans un appartement}

On suppose dorénavant $\I$ munit de son système complet d'appartements.\\

    Considérons une galerie tendue $\Gamma=c_1,...,c_n$ de longueur $n$. Soit $A_1$ un appartement contenant $c_1$, le but est d'étudier quelles parties de $A_1$ sont incluses dans un appartement contenant également $c_n$. On peut supposer que $c_2$ n'est pas dans $A_1$. Soit $M_1$ le mur de $A_1$ contenant la cloison $c_1\wedge c_2$. Alors d'après le corollaire \ref{cor:inclusion dans un appartement}, chacun des deux demi-espaces de $A_1$ délimités par $M_1$ est inclus dans un appartement contenant également $c_2$. Soit $A_2$ un tel appartement contenant un demi-appartement de $A_1$ et $c_2$. On définit $M_2$ comme étant le mur de $A_2$ contenant la cloison $c_2\wedge c_3$. Alors chacun des demi-appartements de $A_2$ peut être inclus dans un appartement contenant $c_3$. Mais ces demi-appartements contiennent en quelque sorte des quarts d'appartement de $A_1$. En continuant ainsi jusqu'au bout de la galerie $\Gamma$, on prouve la:\\
    
    \begin{prop} \label{prop:inclusion d'une partie et d'une chambre dans un appart} Soit $A$ un appartement et $c$ une chambre. Soit $n$ la longueur minimale entre une chambre de $A$ et $c$. Alors il existe un découpage $\mathcal P \subset \mathcal P(A)$ de $A$ en $2^n$ parties fermées, (certaines pouvant être vides) tel que \begin{itemize}
    \item $\mathcal P$ est obtenu en coupant $A$ en deux le long d'un mur, puis en coupant chacun des demi-appartements obtenus encore en deux le long d'autres murs, et en répétant cette opération encore $n-2$ fois
    \item chaque élément de $\mathcal P$ est inclus dans un appartement contenant également $c$.\end{itemize}
    \end{prop}
    \rema Cette proposition est vraie dans un immeuble quelconque.\\

    Voici deux résultats qui seront très utiles par la suite directement obtenus grâce à cette proposition.
    \begin{cor}\label{cor:inclusion d'une galerie ou d'une facette et d'une chambre dans un appart}
    \begin{itemize}
    \item Soit $\Gamma=c_0,c_1,...$ une galerie tendue infinie. Soit $d$ une chambre de $\I$. Alors il existe $n_0\in \N$ et un appartement $A$ contenant $d$ et tous les $c_i$ pour $i\geq n_0$.
    \item Soit $f$ une facette de quartier, et $c$ une chambre. Alors il existe un scp de $f$ qui est inclus dans un appartement contenant également $c$.\end{itemize}\end{cor}
    Pour prouver ce corollaire, il suffit de remarquer que si $A$ est un appartement contenant $\Gamma$, si $M$ est un mur de $A$ alors un des deux demi-appartements de $A$ définis par $M$ contient tous les $c_i$ à partir d'un certain rang, car une galerie tendue ne peut être coupée deux fois par un même mur.\\
     Et de la même manière, si $A$ contient $f$, alors un des demi-appartements définis par $M$ contient un scp de $f$.\cqfd\\
    
    Nous voulons maintenant préciser la proposition \ref{prop:inclusion d'une partie et d'une chambre dans un appart}, en déterminant certaines zones de $A$ qui seront obligatoirement incluses dans un élément de la partition $\mathcal P$. Pour cela, il s'agit d'identifier des parties de $A$ qui ne seront jamais coupées par un mur "séparant $A$ de $c$". Commençons par donner un sens précis à ceci:\\

\begin{prop}\label{prop:inclusion d'une partie de A et d'une galerie dans un appart}
Soit $A\in\A$, soit  $\Gamma=(c_i)_{i\in I}$ une galerie tendue telle que $c_0\in A$, avec $I$ un intervalle de $\N$ contenant $0$. Soit $Z$ un appartement contenant $\Gamma$, et $\phi:Z\isom A$ un isomorphisme fixant $c_0$. Soit $\beta$ une partie de $A$, connexe, contenant $c_0$, telle qu'aucun $\phi(M)$ pour $M$ un mur de $Z$ traversé par $\Gamma$ ne sépare strictement deux points de $\beta$.\\
Alors il existe un appartement contenant $\cl(\beta)\cup\Gamma$.\end{prop}

\demo\\
On construit par récurrence une suite d'appartements $(Y_i)_{i\in I}$ telle que pour tout $i\in I$, on ait $\beta\cup c_0\cup c_1\cup...\cup c_i \subset Y_i$.\\
 Pour $i=0$, $Y_0=A$ convient. Supposons à présent $Y_i$ construit, avec $i\in I$ tel que $i+1\in I$. Soit $M_i$ le mur de $Y_i$ contenant $c_i\wedge c_{i+1}$. Si nous prouvons que $M_i$ ne coupe pas $\beta$ alors ceci entraînera que $\beta\cup c_1\cup...\cup c_i\subset \D(M_i,c_i)$, et le corollaire \ref{cor:inclusion dans un appartement} prouvera l'existence d'un $Y_{i+1}$ convenable.\\
Supposons par l'absurde que $M_i$ rencontre $\beta$. Soit $\chi_i:Y_i\isom A$ l'isomorphisme fixant $\beta$ (il est unique car $\beta$ contient la chambre $c_0$). Alors $\chi_i(M_i)$ rencontre $\beta$. Mais $\chi_i(M_i)=\phi\circ\psi_i(M_i)$, ou $\psi_i:Y_i\isom Z$ est l'isomorphisme fixant $c_0$. Et $\psi_i(M_i)$ est un mur de $Z$ coupant $\Gamma$. Ceci contredit les hypothèses, et prouve donc l'existence de $Y_{i+1}$.\\

  Chaque $Y_i$ contient en fait $\cl(\beta)$, car $A\cap Y_i$ est une partie close. Et $\cl(\beta)$ est égal à l'adhérence de l'ensemble de ses chambres car il contient une chambre. A présent, pour tout triple $(c,d,e)$ de chambres incluse dans $\cl(\beta)\cup\Gamma$, il existe $Y_i$ contenant $c\cup d\cup e$. Ceci prouve, en utilisant le corollaire \ref{cor:du th fonda} l'existence d'un appartement $Y$ contenant $\cl(\beta)\cup \Gamma$.\cqfd\\

    Passons à présent à quelques situations particulières où s'applique ce résultat.\\

\begin{defin}\label{defin:base}
 Soit $A$ un appartement et $\vec\delta$ un cône inclus dans une facette de $\vec A$. Soit $\delta=a+\vec\delta$ un cône de direction $\vec \delta$. On appelle base de $\delta$, et on note $b(\delta)$ l'intersection de $\delta$ et des demi-appartements contenant un voisinage de $a$ pour la topologie induite sur $\delta$.
\end{defin}

Cette définition est indépendante de l'appartement $A$ contenant $\delta$ choisi. De plus, comme les murs de $A$ sont répartis de manière discrète, $b(\delta)$ contient un voisinage de $a$ pour la topologie induite sur $\delta$. En particulier, la base de $\delta$ contient une base affine de $\aff_B(\delta)$, pour tout appartement $B$ contenant $\delta$. Dis autrement, $\aff_B(b(\delta))=\aff_B(\delta)$. Enfin, $b(\delta)$ est inclus dans une facette de $A$.\\
\rema Cette définition sera utilisée en \ref{section:topologie} pour définir la topologie sur l'immeuble compactifié.\\

\begin{prop}\label{prop:cone et chambre dans un appart}
Soit $A$  un appartement, $\vec\delta$ une facette de $\vec A$ non triviale, $a\in A$, on note $\delta:=a+\vec \delta$. Soit $b$ la base de $\delta$.
Soit $B$ un autre appartement contenant $b$, et $\psi:A\isom B$ un isomorphisme fixant $b$.
Soit $x$ une facette de $B$ incluse dans $a-\vec\psi(\vec\delta^*_A)$ et telle qu'aucun mur dont la direction contient $\vec\delta$ ne sépare strictement $a$ de $x$.\\
 Alors pour toute chambre $\vec C\subset \vec A$ incluse dans $\vec\delta^*_A$, il existe un appartement $Z$ contenant $(a+\vec C)\cup \{x\}$.\\
En particulier, pour tout point $t$ de $\delta^*_A$, il existe un appartement contenant $t$, $a$ et $x$.\end{prop}

\rema \begin{itemize}
\item  $\psi(\delta)$ et $\vec\psi(\vec\delta)$ sont indépendants du choix de $\psi$.
\item Comme $b\subset A\cap B$ contient une base affine de $\delta$, on a en fait $\vec\delta\subset \vec A\cap\vec B$. Cependant, pour les démonstrations de cette sous-partie, il me parait plus clair de ne pas faire cette identification, c'est-à-dire de distinguer $\vec\delta$ de $\vec\psi(\vec\delta)$.
\item Dans la suite, on pourra noter $-\delta=a-\vec\delta$ le cône opposé à $\delta$, et de manière générale si $f$ est un cône ne contenant pas de droite, $-f=s(f)-\vec f$.\\
 \end{itemize}

\demo\\
Soit $c_0$ une chambre de $A$, coupant $a+\vec C$ et dont l'adhérence contient $b$.
 Soit $\Gamma=c_0,...,c_k$ une galerie tendue entre $c_0$ et $x$. Si $Z$ est un appartement contenant $\Gamma$, et $\psi_Z:A\isom Z$ un isomorphisme fixant $b$, alors $x$ est une facette de $Z$, les seuls murs de $Z$ dont la direction contient $\vec\psi_Z(\vec\delta)$ et qui séparent $a$ de $x$ contiennent $a$ ou $x$, et enfin $x\in -\psi_Z(\delta)^*$. On est ainsi ramené au cas où l'intersection des deux appartements considérés contient la chambre $c_0$.\\

 Soit $\beta=(a+\vec C)\cup c_0$, c'est un connexe contenant $c_0$, il faut maintenant vérifier que pour tout mur $M$ de $\Gamma$, $\psi_Z\inv(M)$ ne coupe pas $\beta$.\\
Soit $M$ un tel mur. Alors $\psi_Z\inv(M)$ sépare strictement $c_0$ et $x$, donc coupe $-\delta^*$. Si en plus $\psi_Z\inv(M)$ coupe $\beta$ alors il coupe $a+\vec C$ donc en particulier l'intérieur de $\delta^*$. Ceci implique $\vec\delta\subset \fl{\psi_Z\inv(M)}$ puis $\vec\delta\subset\vec M$. Mais alors $x\subset M$ ou $a\in M$. Le premier cas est exclus car $M$ sépare strictement $x$ de $c_0$. Quand au second cas, il implique $\delta\subset M$, ce qui empêche que $M$ coupe $a+\vec C$.\cqfd\\

 Dans la variante suivante de ce résultat, on veut obtenir un appartement contenant $\delta^*\cup x$. Pour cela
on peut renforcer l'hypothèse sur les murs dont la direction contient $\vec\delta$, en imposant par exemple que les seuls murs séparant $b$ de $x$ contiennent $x$.\\
 Cependant, nous nous contenterons de l'hypothèse $x\subset -\psi( \delta)$, qui présente l'avantage qu'en poussant un tout petit peu le raisonnement, on arrive à trouver un appartement contenant non seulement $x$, mais en fait $\delta^*\cup-\psi( \delta)$. Ceci donne le résultat:\\

 \begin{prop}\label{prop:coeur et etoile dans un appart} Soient $A$, $\delta$, $a$, $b$, $B$, $\psi$ comme précédemment.\\
Alors il existe un appartement contenant $\delta^*_A\cup-\psi( \delta)$
% Soit $x$ une facette de $B$ tel que les seuls murs de $B$ dont la direction contient $\vec\psi(\vec\delta)$  qui séparent $b$ de $x$ contiennent $b$ et $x$. On suppose en outre que $x\in -\psi(\delta)^*$ (ce qui signifie  $x\in a-\vec\psi(\vec\delta)^*$).\\
%Alors un existe un appartement contenant $\delta^*_A\cup x$.
\end{prop}

\demo\\
 On commence, comme pour la preuve précédente, par se ramener au cas où $A\cap B$ contient une chambre. Soient $c$ et $d$ deux chambres contenant $b$, incluses respectivement dans $A$ et $B$. Soit $\Lambda=d,d_1,...,d_k,c$ une galerie entre $d$ et $c$. En utilisant $k+1$ fois le corollaire \ref{cor:inclusion dans un appartement}, on trouve une suite $B_0=B,B_1,...,B_{k+1}$ d'appartements tels que $d_i\subset B_i$ (en notant $d_0=d$ et $d_{k+1}=c$) et tels que $B_i\cap B_{i+1}$ contient un demi-appartement délimité par le mur $M_i$ de $B_i$ contenant $d_i\wedge d_{i+1}$. Or ce mur contient $b$, on vérifie alors par une simple récurrence que pour tout $i$, $M_i$ contient $-\psi(\delta)$.\\
Soit $Z=B_{k+1}$: cet appartement contient $-\psi(\delta)$ et une chambre $c$ de $A\cap Z$. Soit $\phi:Z\isom A$ fixant $c$, alors $-\psi(\delta)=-\phi\inv(\delta)$.\\

Soit $l$ une demi-droite dans $Z$ passant par $a$, dont la direction est intérieure à $-\vec\phi\inv(\vec\delta)$, et dont l'origine est dans $b$. Elle est donc incluse dans $-\psi(\delta)\cup b$. Soit $\Gamma=(c_i)_{i\in\N}$ une galerie tendue le long de $l$, donnée par la proposition \ref{prop:galerie sur droite} telle que $c_0=c$. Soit $\beta=\cl(c_0\cup \delta^*_A)\subset A$. Pour pouvoir appliquer la proposition \ref{prop:inclusion d'une partie de A et d'une galerie dans un appart}, il nous reste à prouver que pour tout mur $M$ traversé par $\Gamma$, $\phi(M)$ ne coupe pas $\beta$.\\

Soit $M$ un tel mur, il coupe $l$ sans la contenir, donc $\vec l\not\subset \vec M$, d'où $\vec\phi\inv(\vec\delta)\cap \vec M=\vide$ et $\vec M$ ne peut couper $\vec\phi\inv(\vec\delta^*)$.\\
Notons $\vec M^+$ et $\vec M^-$ les deux demi-appartements de $\vec Z$ délimités par $\vec M$ en choisissant $\vec\phi\inv(\vec\delta^*)\subset \vec M^+$. Soient $M^+$ et $M^-$ les demi-appartements de $Z$ correspondants, nous avons alors $c_0\subset M^+$. Passons du côté de $A$: nous savons que $c_0\subset \phi(M^+)$ et $\vec\delta^*\subset \vec\phi(\vec M^+)$. Donc $\beta\subset\cl(c_0+\vec\delta^*)\subset \phi(M^+)$. Ainsi, $\phi(M)$ ne coupe pas $\beta$.\\
 D'après la proposition \ref{prop:inclusion d'une partie de A et d'une galerie dans un appart}, il existe un appartement $Y$ contenant $\beta\cup \Gamma$. Alors $l\subset Y\cap Z$ d'où $\cl(l)\subset Y\cap Z$. Mais $\cl(l)\supset -\phi\inv(\delta)$: c'est une conséquence de la proposition \ref{prop:demi droite caracteristique} (en fait $\cl(l)$ est égal à la cheminée $\cl_Z(b-\vec\phi\inv(\vec\delta) )$). Comme $-\phi\inv(\delta) = -\psi(\delta)$, la proposition est prouvée.\cqfd\\

    \subsection{Rétractions par rapport à deux chambres adjacentes}
    
    Nous étudions ici comment change une rétraction lorsqu'on remplace la chambre de base par une chambre adjacente.\\% Le résultat est simple et sera souvent utile dans la suite.\\
    
    \begin{prop}\label{prop:retractions adjacentes}
     Soient $c_1$ et $c_2$ deux chambres adjacentes dans un appartement $A$. Soit $f$ la cloison $c_1\wedge c_2$ et $M$ le mur de $A$ contenant $f$, $\sigma$ la réflexion selon $M$. Soient encore $\rho_1=\rho_{c_1,A}$ et $\rho_2=\rho_{c_2,A}$. Alors pour toute chambre $d\subset \I$,
     \begin{itemize}
     \item $\rho_2(d)=\rho_1(d)$ ou bien $\rho_2(d)=\sigma\circ\rho_1(d)$.
     \item Si $\rho_1(d)$ est du même côté de $M$ que $c_1$ alors $\rho_2(d)=\rho_1(d)$.
     \item Si $\rho_1(d)$ est du même côté de $M$ que $c_2$ alors les deux possibilités peuvent survenir, et on a $\rho_2(d)=\rho_1(d)$ si et seulement si $c_1$ est dans l'enclos de $d$ et $c_2$.
     \end{itemize}
    \end{prop}
    
    on peut aussi énoncer un résultat un peu moins précis mais plus clair:\\
    \begin{cor} \label{cor:retractions adjacentes} On reprend les hypothèses de la proposition. Alors $\rho_2(d)=\rho_1(d)$ si et seulement si il existe un appartement contenant $c_1$, $c_2$ et $d$.\end{cor}
    
    \rema Le sens $\Leftarrow$ du corollaire est de toute façon une conséquence directe de la définition d'une rétraction.\\
    
    \demo\\
    Soit $B$ un appartement contenant $c_1$ et $d$, soit $\phi:B\rightarrow A$ l'isomorphisme fixant $c_1$. Alors $\rho_1(d)=\phi(d)$.\\
    Si $\rho_1(d)$ est du même côté de $M$ que $c_1$ alors $d$ est du même côté de $\phi\inv(M)$ que $c_1$. En vertu du corollaire \ref{cor:inclusion dans un appartement}, il existe un appartement $Z$ contenant $d$, $c_1$ et $c_2$, ce qui prouve que $\rho_1(d)=\rho_2(d)$.\\
    
    Dans l'autre situation, appelons $c_2'$ la chambre adjacente à $c_1$ le long de $f$ dans $B$. Cette chambre $c_2'$ est donc dans l'enclos de $c_1$ et de $d$. Si $c_2'=c_2$ alors $c_1$, $c_2$ et $d$ sont dans un même appartement, d'où $\rho_1(d)=\rho_2(d)$ et $c_2$ est bien dans l'enclos de $c_1$ et de $d$. Il ne reste donc qu'à vérifier que si $c_2\not=c_2'$ alors $\rho_2(d)=\sigma\circ\rho_1(d)$ et $c_2'$ n'est pas dans l'enclos de $c_1$ et $d$. Le deuxième point est prouvé car nous disposons d'un appartement, $B$ qui contient $c_1$ et $d$, mais pas $c_2$ sans quoi il y aurait trois chambres contenant $f$ dans $B$.\\
    Pour prouver le premier point, constatons que la formule $\rho_2(x)=\sigma\circ\rho_1(x)$ est vraie lorsque $d=c_2'$. Soit $B^+$ le demi appartement délimité par $\phi\inv(M)$ et contenant $c_2'$, donc également $d$. Par le corollaire \ref{cor:inclusion dans un appartement}, il existe un autre appartement $B'$ contenant $B^+$ et $c_2$. Si $\phi':B'\fleche A$ est l'isomorphisme correspondant, alors $\rho_2(d)=\phi'(d)$ et l'égalité $\rho_2(c_2)=\sigma\circ\rho_1(c_2)$ devient $\phi'(c_2)=\sigma\circ\phi(c_2)$. Restreignons $\sigma\circ\phi$ et $\phi'$ à $B^+$: nous obtenons deux morphismes de complexe de chambre injectifs de $B^+$ dans $A$, qui co\"incident en $c_2$. Comme $B^+$ est convexe, ceci entraîne que $\phi'$ et $\sigma\circ\phi$ co\"incident sur $B^+$.\cqfd

    \subsection{Inclusion de deux galeries tendues dans un appartement}
   
\subsubsection{Le théorème}

   \begin{defin} Soit $\Gamma=(c_i)_{i\in I}$ une galerie. Une sous-galerie de $\Gamma$ est une galerie obtenue en enlevant un nombre fini de chambres à $\Gamma$. Dans le cas où $I=\N$ et où $\Gamma$ est tendue, une sous-galerie de $\Gamma$ est donc de la forme $(c_i)_{i\geq i_0}$ pour un $i_0\in \N$.\end{defin}

   \begin{prout}\label{th:deux galeries dans appart}
    Soient $\Gamma=(c_i)_{i\in \N}$ et $\Delta=(d_i)_{i\in \N}$ deux galerie tendues. Alors il existe une sous-galerie de $\Gamma$ et une sous-galerie de $\Delta$ qui sont incluses dans un même appartement.
   \end{prout}
   
   \demo\\
   \textit{Notations:}\\
    Soit $A$ un appartement contenant $\Gamma$. Pour $i\in \N$, soit $\rho_i=\rho_{c_i,A}$, $M_i$ le mur de $A$ entre $c_i$ et $c_{i+1}$, $\sigma_i$ la réflexion selon $M_i$, de sorte que pour tout $i$ et pour toute chambre $c$ de $\I$, on a en vertu de la proposition \ref{prop:retractions adjacentes} $\rho_{i+1}(c)=\rho_i(c)$ ou $\rho_{i+1}(c)=\sigma_i\circ\rho_i(c)$. Remarquons que par continuité de $\rho_{i+1}$, la même formule s'applique pour le calcul de $\rho_{i+1}$ pour deux chambres adjacentes, sauf dans le cas où ces deux chambres sont projetées par $\rho_i$ sur des chambre adjacentes à $M_i$.\\
   Choisissons un point spécial $0$ pour identifier $A$ et $\vec A$. Quitte à réduire $\Gamma$, cette galerie est incluse dans une chambre vectorielle $C$. Grâce à la proposition $\ref{prop:inclusion d'une partie et d'une chambre dans un appart}$, on voit que, quitte à la réduire, la rétraction d'une galerie tendue est une galerie tendue. En particulier, $\rho_i(\Delta)$ est incluse à translation près dans une chambre vectorielle. Choisissons $D_i$ une chambre vectorielle à distance maximale de $C$ contenant un translaté de $\rho_i(\Delta)$. Notons enfin $\delta_i$ la distance entre la chambre $C$ et $D_i$. Comme le complexe de Coxeter vectoriel est fini, les $\delta_i$ sont majorées par une certaine constante $\delta_{max}$.\\
   
   \textit{plan de la preuve:}\\
   Pendant cette démonstration, nous allons étudier en quoi $\rho_i$ peut différer de $\rho_{i+1}$. Nous verrons qu'il y a trois cas à distinguer: soit $\rho_i$ et $\rho_{i+1}$ coïncident sur presque toutes les chambres de $\Delta$, soit elles diffèrent sur un nombre infini de ces chambres et $\delta_{i+1}>\delta_i$, soit elles diffèrent sur un nombre infini de ces chambres et $\delta_{i+1}=\delta_i$. Nous allons commencer par traiter le cas où seul le premier cas intervient (préliminaire 1). Ensuite, nous allons vérifier que les trois cas cités sont effectivement les seuls possibles, nous en profiterons pour analyser un peu plus en détail le deuxième cas (préliminaire 2). Alors nous verrons qu'on peut réduire $\Gamma$ pour se ramener à une situation où le deuxième cas n'apparaît pas, puis qu'on peut la réduire une seconde fois pour que le troisième cas n'apparaisse plus non plus. \\

   \textit{Préliminaire 0: le cas très facile:}\\
   Supposons que pour tout $i\in \N$ et pour toute chambre $d$ de $\Delta$, $\rho_i(d)=\rho_{i+1}(d)=\rho(d)$. Cela implique par le corollaire \ref{cor:retractions adjacentes} que $d$, $c_i$ et $c_{i+1}$ sont inclus dans un même appartement pour tout $i$ et $d$. Par ailleurs, quitte à réduire $\Delta$, on peut supposer que $\Delta$ est incluse dans un appartement contenant aussi $c_1$. L'égalité des $W$-distances $\delta(a,c)=\delta(a,b)\delta(b,c)$ est alors réalisée dans le cas où $\{a,b,c\}=\{c_1,d,d'\}$, avec $d,d'\in\Delta$, et dans le cas où $\{a,b,c\}=\{c_i,c_{i+1},d\}$, avec $i\in\N$ et $d\in \Delta$ (et bien sûr aussi dans le cas $\{a,b,c\}\subset \Gamma$ ou $\{a,b,c\}\subset \Delta$). Le cas $\{a,b,c\}=\{c_i,c_j,d\}$, $i,j\in\N$ et $d\in\Delta$, s'obtient aisément par récurrence sur $|i-j|$. Le dernier cas est $\{a,b,c\}=\{d,d',c_i\}$ avec $d,d'\in\Delta$, il se traite à partir des cas précédents: $\delta(d,d')=\delta(d,c_1)\delta(c_1,d') =\delta(d,c_i)\delta(c_i,c_1)\delta(c_1,c_i)\delta(c_i,d') =\delta(d,c_i)\delta(c_i,d')$.\\
   Il ne reste plus qu'à appliquer le théorème \ref{th:fonda}: $\Gamma$ et $\Delta$ sont incluses dans un même appartement.\\
   
   \textit{Préliminaire 1: le cas facile:}\\
   Supposons que pour tout $i$, $\rho_i(d)=\rho_{i+1}(d)=\rho(d)$ pour toutes les chambres $d$ de $\Delta$ sauf un nombre fini (ce nombre pouvant a priori varier avec $i$). Nous allons construire par récurrence une sous-galerie $\Delta'$ de $\Delta$ de sorte que $\Gamma$ et $\Delta'$ soient dans la situation très facile du préliminaire 0. Plus précisément, nous allons construire pour tout $i\in\N$ une sous-galerie $\Delta_i$ de $\Delta$ telle que $\forall k< i$, $\forall d\in \Delta_i$, $\rho_k(d)=\rho_{k+1}(d)$. Pour des raisons techniques, nous allons imposer en outre que $\rho_i(\Delta_i)=\rho_1(\Delta_i)$ est une galerie tendue. Ensuite, nous vérifierons que la suite $(\Delta_i)$ est stationnaire à partir d'un certain rang, alors la "limite" de cette suite conviendra.\\
   Commençons par choisir pour $\Delta_1$ une sous galerie de $\Delta$ telle que $\rho_1(\Delta_1)$ est tendue. Supposons ensuite construite $\Delta_i$, et construisons $\Delta_{i+1}$. Supposons qu'il existe $d_k\in\Delta_i$ tel que $\rho_{i+1}(d_k)=\sigma_i\circ\rho_i(d_k)$. Par la proposition \ref{prop:retractions adjacentes}, $\rho_i(d_k)$ est dans le demi-appartement $\D(M_i,c_{i+1})$. Comme par hypothèse seul un nombre fini de chambres $e$ de $\Delta$ vérifient $\rho_{i+1}(e)=\sigma_i\circ\rho_i(e)$, il existe $d_l\in\Delta_i$ tel que $\rho_{i+1}(d_l)=\rho_i(d_l)$. Nous pouvons choisir $d_k$ et $d_l$ adjacentes, c'est à dire $l=k\pm 1$. Alors $\rho_i(d_k),\rho_i(d_l)=\rho_{i+1}(d_l),\rho_{i+1}(d_k)$ est une galerie de longueur 2 dans $A$, donc $d(\rho_i(d_k),\rho_{i+1}(d_k))\leq 2$. Mais $\rho_{i+1}(d_k)=\sigma_i(\rho_i(d_k))$ donc  $d(\rho_i(d_k),\rho_{i+1}(d_k))$ est impair. En effet $\sigma_i$ induit une involution sur l'ensemble des murs séparant $\rho_i(d_k)$ et $\rho_{i+1}(d_k)$, et cette involution n'a qu'un point fixe, $M_i$, car un mur stabilisé par $\sigma_i$ ne peut séparer un point et son image par $\sigma_i$. Il y a donc un nombre impair de murs séparant $\rho_i(d_k)$ et $\rho_{i+1}(d_k)$.\\
 Ainsi $d(\rho_i(d_k),\rho_{i+1}(d_k))=1$, les deux chambres $\rho_i(d_k)$ et $\rho_{i+1}(d_k)$ sont de part et d'autre de $M_i$, et $\rho_{i+1}(d_k)=\rho_{i+1}(d_l)=\rho_i(d_l)$. Comme $\rho_i(\Delta_i)$ est tendue, elle ne peut couper $M_i$ une autre fois: $k$ et $l$ sont donc les seuls indices adjacents tels que $\rho_{i+1}(d_k)=\sigma_i\circ\rho_i(d_k)$ et $\rho_{i+1}(d_l)=\rho_i(d_l)$. Nous en déduisons que la première formule est vraie pour tous les indices inférieurs à $k$, et la seconde pour tous les indices supérieurs à $k$ (en particulier, $l=k+1$). Nous avons donc une description précise de la situation: une partie finie de $\rho_i(\Delta)$ est dans le demi-appartement $\D(M_i,c_{i+1})$, et c'est précisément cette partie qui va être déplacée lors du passage à $\rho_{i+1}$. On pose $\Delta_{i+1}=(d_j)_{j\geq k}$, et cela conclut la construction de la suite $(\Delta_i)$.\\
   Il reste à vérifier qu'au cours de cette construction, il n'y a qu'un nombre fini d'étapes où on a réellement retiré des chambres à $\Delta_i$ pour obtenir $\Delta_{i+1}$. Soit $I$ l'ensemble des indices $i$ tel que $\Delta_{i+1}$ est différent de $\Delta_i$. Nous allons montrer qu'il ne peut y avoir deux indices $i$ et $j$ dans $I$ tels que $\vec M_i=\vec M_j$, ce qui prouvera que $I$ est fini puisque l'ensemble des murs vectoriels est fini. Supposons par l'absurde l'existence de tels $i$ et $j$, avec $i<j$. Lors de la construction de $\Delta_{i+1}$, on a fait en sorte que $\rho_{i+1}(\Delta_{i+1})$ soit inclus dans $\D(M_i,c_i)$. Comme $\Gamma$ est tendue et $\vec M_i=\vec M_j$, $\D(M_i,c_i)\subset\D(M_j,c_j)$. D'après la construction, la galerie $\rho_j(\Delta_j)$ est une sous galerie de $\rho_i(\Delta_i)$, elle est donc incluse dans $\D(M_j,c_j)$. Mais cela implique alors que $\rho_{j+1|\Delta_j}=\rho_{j|\Delta_j}$ et cela contredit le fait que $j\in I$. Ce qui prouve le cas facile.\\

   \textit{Préliminaire 2: que se passe-t-il lorsque $\rho_i\not=\rho_{i+1}$ pour un nombre infini de chambres de $\Delta$?}\\
   Il s'agit à présent d'étudier ce qui se passe lorsque $\rho_{i+1}\not=\rho_i$ pour un nombre infini de chambres de $\Delta$, c'est à dire pour toutes les chambres de $\Delta$ à partir d'un certain rang, puisque $\rho_i(\Delta)$ est tendue à partir d'un certain rang. Supposons que $\delta_{i+1}\leq \delta_i$ (c'est le seul cas qu'il est utile d'étudier, comme nous allons le voir).\\
    Tout d'abord, $\rho_i(\Delta)$ est à partir d'un certain rang dans le demi appartement $\D(M_i,c_{i+1})$.\\
     Par ailleurs, $\rho_{i+1}(\Delta)$ est alors à translation près dans la chambre $\vec\sigma_i(D_i)$, donc $d(C,D_i)\geq d(C,\vec\sigma_i(D_i))$. Cela signifie que $\vec M_i$ sépare $C$ et $D_i$. Mais $M_i$ est un translaté de $\vec M_i$ dans la direction de $C$. (Pour résumer les placements des différentes parties étudiées par rapport aux deux murs $M_i$ et $\vec M_i$, je note: $\Gamma,c_{i+1},\rho_i(\Delta) |^{M_i} c_i |^{\vec M_i} D_i$, où $|^M$ symbolise le mur $M$.) Le point qui nous intéresse est celui-ci: un translaté de $\rho_i(\Delta)$ est dans $\D(M_i,c_{i+1})$ donc dans $\D(\vec M_i,C)$ et un autre translaté est dans $D_i$ donc dans $-\D(\vec M_i,C)$. Cela implique que $\rho_i(\Delta)$ reste à distance bornée de $\vec M_i$, au sens où il existe un majorant $M$ à la distance entre une chambre de $\rho_i(\Delta)$ et $M_i$.\\
     Alors $2M$ est un majorant de la distance entre $\rho_i(d)$ et $\rho_{i+1}(d)$ pour les chambres $d\in\Delta$. Ceci implique que $\rho_{i+1}(\Delta)$ est incluse à translation près dans la même chambre vectorielle que $\rho_i(\Delta)$, c'est-à-dire dans $D_i$. Mais comme $\delta_{i+1}\leq\delta_i$, on trouve que $\delta_{i+1}=\delta_i$ et on peut choisir $D_{i+1}=D_i$.\\

    \textit{le coeur de la démonstration:}\\
    Lorsque $\rho_i=\rho_{i+1}$ pour presque toutes les chambres de $\Delta$, il est clair que $\delta_i=\delta_{i+1}$, et qu'on peut choisir $D_{i+1}=D_i$. Lorsque $\rho_i\not=\rho_{i+1}$ pour un nombre infini de ces chambres, nous venons de voir que $\delta_{i+1}\leq\delta_i$ implique que $\delta_{i+1}=\delta_i$ et qu'on peut encore supposer $D_{i+1}=D_i$.
    Donc la fonction $i\mapsto\delta_i$ est croissante. Comme elle est bornée par $\delta_{max}$, elle est constante à partir d'un certain rang. Quitte à raccourcir $\Gamma$, on peut la supposer constante. Nous sommes alors pour tout $i$ soit dans le cas étudié dans le préliminaire 2, soit dans le cas $\rho_i=\rho_{i+1}$ pour presque toutes les chambres de $\Delta$. En particulier on peut supposer que pour tout $i$, $D_i=D_1$. Il s'agit bien sûr de se ramener à présent au cas du préliminaire 1.\\
    
    Soit $F$ la facette vectorielle égale à l'intersection de tous les murs vectoriels à distance bornée desquels reste $\rho_1(\Delta)$. Soit $\mathcal M$ l'ensemble des murs vectoriels contenant $F$, soit $W^v_F=\{w\in W^v|wF=F\}$. Nous allons montrer que pour tout $i$, $\mathcal M$ est égal à l'ensemble des murs à distance bornée desquels reste $\rho_i(\Delta)$.\\
     Commençons par $i=1$. Soient $M_1,..,M_k$ les murs à distance bornée desquels reste $\rho_1(\Delta)$, il s'agit en fait de vérifier que $\mathcal M=\{M_1,...,M_k\}$. Soient $\alpha_1,..,\alpha_k$ des formes linéaires tels que $M_j=\ker \alpha_j$. Soit $M=\ker \alpha $ un mur dans $\mathcal M$. Alors $\alpha$ est combinaison linéaire des $\alpha_j$. Mais chaque $\alpha_j$ est bornée sur $\rho_1(\Delta)$, donc $\alpha$ est également bornée sur $\rho_1(\delta)$, c'est à dire que $\rho_1(\Delta)$ reste à distance bornée de $M$.\\
     Supposons à présent le résultat vrai au rang $i$. Si $\rho_i$ et $\rho_{i+1}$ ne diffèrent que sur un nombre fini de chambres de $\Delta$, le résultat reste évidemment vrai au rang $i+1$. Sinon, nous sommes dans le cas étudié dans le préliminaire 2. Donc $\rho_i(\Delta)$ reste à distance bornée de $M_i$, et par hypothèse de récurrence, $\vec M_i\in \mathcal M$, et $\sigma_i\in W^v_F$. L'ensemble des murs à distance bornée de $\rho_{i+1}(\Delta)$ est donc $\sigma_i(\mathcal M)$, mais ceci vaut $\mathcal M$ puisque $\sigma_i$ fixe $F$.\\
    
 %   Comme $\rho_1(\Delta)$ est une galerie tendue à partir d'un certain rang, et comme elle reste à distance bornée de chaque $M\in\mathcal M$, on peut, quitte à réduire $\Delta$, se placer dans le cas où aucun mur parallèle à un mur de $\mathcal M$ ne coupe $\rho_1(\Delta)$. Une récurrence facile montre que cette propriété est alors encore vérifiée par tous les $\rho_i(\Delta)$.\\
    On utilise maintenant la compactification de $A$ définie par $\F=\{$facettes vectorielles de Weyl$\}$ (c'est à dire la compactification polygonale classique). Soit $p$ la projection de $A$ sur l'espace affine $A_F\simeq A/\vect(\vec F)$. On prend $p(0)$ comme origine dans $A_F$ pour l'identifier à son espace vectoriel directeur.  L'ensemble $\mathcal M$ s'identifie (via $p$) à l'ensemble des murs vectoriel de $A_F$, $W^v_F=\stab_{W^v}(F)$ est le groupe de Coxeter vectoriel. L'ensemble des murs affines de $A_F$ est $\{p(M)|M$ mur affine de $A$ et $\vec M\in \mathcal M\}$, et si $M\in\mathcal M$ on notera encore $M$ au lieu de $p_F(M)$. Enfin, si $w\in W$ vérifie $\vec w\in W^v_F$ alors $w$ agit sur $A_F$ en préservant les murs, notons $\tilde w\in GA(A_F)$ l'automorphisme affine induit par $w$ sur $A_F$. (Attention: la fonction $w\mapsto \tilde w$ n'est pas injective et $\tilde w$ n'est pas forcement dans le groupe de Coxeter de $A_F$.)\\
    Chaque galerie $\rho_i(\Delta)$ est tendue à partir d'un certain rang, et reste à distance bornée des murs de $\mathcal M$, elle ne coupe donc qu'un nombre fini de murs parallèles à un mur de $\mathcal M$. Donc à partir d'un certain rang, sa projection sur $A_F$ reste dans une unique chambre de $A_F$, notons cette chambre $g_i$.
      Les chambres vectorielles $C$ et $D_1$ sont également projetées dans des chambres vectorielles de $A_F$ (de manière pas nécessairement surjective a priori), qu'on notera encore $C$ et $D_1$ pour simplifier. Nous avons vu lors du préliminaire 2 que $\vec M_i$ sépare $C$ et $D_1$, si $\rho_{i+1}\not=\rho_i$ pour un nombre infini de chambres de $\Delta$. Comme $\vec M_i\in \mathcal M$, $M_i$ est un mur de $A_F$, et la relation "$M_i$ sépare $D_1$ et $C$" est encore vraie dans $A_F$. On a de plus dans ce cas que $g_{i+1}=\tilde\sigma_i(g_i)$. Mais $M_i$ est un translaté de $\vec M_i$ dans la direction de $C$, donc $M_i$ sépare $g_{i}$ et $D_1$. Ceci entraîne que $g_{i+1}$ est strictement plus proche de $D_1$ que $g_{i}$, au sens où pour toute chambre $e$ incluse dans $D_1$, $d(e,g_{i+1})<d(e,g_i)$.\\
     Et s'il existe $i$ tel que $g_i\in D_1$, alors la suite des $g_k$ est constante pour $k\geq i$. Car si $g_k\subset D_1$, il n'est plus possible que $M_k$ sépare $D_1$ de $g_k$, donc plus possible que $g_{k+1}=\tilde\sigma_k(g_k)$.\\
     Nous avons donc prouvé que la suite $(g_i)$ est finalement stationnaire. Nous pouvons raccourcir $\Gamma$ de sorte qu'elle soit complètement stationnaire. Cela implique que pour tout $i$, $\rho_i\not=\rho_{i+1}$ seulement pour un nombre fini de chambres de $\Delta$. Mais c'est précisément le cas, facile, du préliminaire 1: la preuve est faite.\cqfd\\

   \subsubsection{Conséquences}
   
   Nous voici en mesure d'énoncer les résultats les plus forts de cette partie. Nous allons voir par exemple qu'étant donnés deux cônes affines $f$ et $g$, il existe un appartement $A$ contenant non seulement des sous-cônes parallèles de $f$ et $g$, mais encore un "épaississement" de ces sous-cônes parallèles. Précisément, $A$ contiendra des cheminées de direction $\vec f$ et $\vec g$, contenant des scp de $f$ et $g$ et dont les bases sont des chambres (\ref{prop:cheminees dans appart}).\\

 Pour commencer, le théorème permet de prouver que deux demi-droites peuvent être réduites pour être incluses dans un même appartement:

\begin{prop} Soient $\delta_1$ et $\delta_2$ deux demi-droites. Alors il existe $\delta_1'$ et $\delta_2'$ deux demi-droites incluses respectivement dans $\delta_1$ et $\delta_2$ et un appartement $Z$ tel que $\delta_1'\cup\delta_2'\subset Z$. Dès lors, $Z$ contient $\cl(\delta_1')\cup\cl(\delta_2')$.
\end{prop}
\textit{preuve:} Cela découle de la proposition \ref{prop:galerie sur droite} et du théorème \ref{th:deux galeries dans appart}\cqfd\\

Ceci permet de traduire le théorème en terme de cheminées:

      \begin{prop} \label{prop:cheminees dans appart} Soient $R_1$ et $R_2$ deux cheminées, dans les appartement respectifs $A_1$ et $A_2$. Alors il existe $R_1'\; scp\; R_1$ et $R_2' \;scp\; R_2$ deux sous cheminées pleines qui sont incluses dans un même appartement.
   \end{prop}
   
   \demo\\
    Soient $\delta_1$ et $\delta_2$ des demi-droites caractéristiques de $R_1$ et $R_2$. Soit $B$ contenant une demi-droite $\delta_1'$ incluse dans $\delta_1$ et une demi-droite $\delta_2'$ incluse dans $\delta_2$. Alors $B$ contient $\cl(\delta_1')$ qui est une scp de $R_1$ et $\cl(\delta_2')$ qui est une scp de $R_2$.\cqfd\\
   
 Le même procédé permet quelques variations:\\

	\begin{prop} \label{prop: facettes dans appart}
 \begin{enumerate}
\item Soit $R$ une cheminée et $c$ une chambre, alors il existe un appartement contenant $c$ et une scp de $R$
\item Soient deux cônes $f=x+\vec f$ et $g=y+\vec g$ où $\vec f\subset \vec A$ et $\vec g\subset \vec B$ sont inclus dans des facettes de Weyl. Alors il existe un appartement contenant des sous-cônes parallèles $f'$ et $g'$ de $f$ et $g$.
\item Soit $f=x+\vec f$ comme au point précédent, et soit $c$ une chambre, alors il existe un appartement contenant $c$ et un scp de $f$.
\end{enumerate}
\end{prop}
   
\demo
\begin{enumerate}
\item Soit $\delta$ une demi-droite caractéristique de $R$, $\Gamma$ une galerie tendue la contenant dans son adhérence.  D'après \ref{cor:inclusion d'une galerie ou d'une facette et d'une chambre dans un appart}, il existe $Z\in\A$ contenant $c$ et une sous-galerie de $\Gamma$, alors $Z$ contient $c$ et une scp de $R$.\\
\item Soit $\delta_1=x+\vec\delta_1$ et $\delta_2=y+\vec \delta_2$ avec $\vec\delta_1\subset\vec f$ et $\vec\delta_2\subset\vec g$. Soit $B\in\A$ et $\delta_1'$, $\delta_2'$ des sous-demi-droites de $\delta_1$ et $\delta_2$ tels que $\delta_1'\cup \delta_2'\subset B$. Alors $\cl(\delta_1')\cup\cl(\delta_2')\subset B$. Mais nous savons que $\cl(\delta_1')=\cl(c_1+\vec h_1)$ où $c$ est la facette affine contenant le sommet de $\delta_1'$ et $\vec h_1$ est la facette vectorielle contenant la direction de $\delta_1'$, c'est-à-dire $\vec\delta_1$. Alors $\vec f\subset \vec h_1$. De plus, $c_1$ contient le sommet $s_1$ de $\delta_1'$, qui est un point de $f$. Donc au final, $\cl(c_1+\vec h_1)$ contient un scp $s_1+\vec f$ de $f$. De même, $\cl(\delta_2')$ contient un scp de $g$.\\
\item Ce point est semblable aux précédents.\end{enumerate}\cqfd\\

   Faisons tout de suite le lien avec les coeurs de cônes affines:
 \begin{cor} Soient $f\in \F_A$ et $g\in \F_B$ deux cônes affines. Alors il existe un appartement $Z$ contenant un scp de $\delta(f)$ et un scp de $\delta(g)$.\end{cor}
 Ceci provient du fait que $\cl(\delta(f))$ et $\cl(\delta(g))$ sont des cheminées.\\ 
   
\rema Tout ceci généralise le fait bien connu que deux quartiers contiennent des sous-quartiers inclus dans un même appartement (\cite{ronan} proposition 9.5 page 121).\\

\subsection{Systèmes d'appartements}

Rappelons que selon les conventions adoptées ici, un immeuble est muni par définition d'un système d'appartements, de sorte que rigoureusement, l'immeuble étudié devrait être dénoté par $(\I,\A)$, $\I$ n'étant que l'ensemble sous-jacent. Si la définition des facettes est indépendante du système d'appartements choisis, celle des quartiers et autres cônes lui est au contraire très liée.\\
Jusqu'ici, nous avons toujours supposé que $\A$ était le système complet d'appartements, mais il peut être intéressant de noter que les résultats précédents sont encore vrais pour d'autres systèmes d'appartements (la construction d'un compactifié de $\I$ ne sera cependant possible qu'en munissant $\I$ de son système complet d'appartements). Nous allons étudier ici le cas d'un système d'appartements $\A'$ qui permet la définition d'un immeuble sphérique à l'infini $(\I,\A')^\infini$, c'est-à-dire que deux quartiers de $(\I,\A')$ contiennent toujours des sous-quartiers contenus dans un même appartement de $\A'$. Un tel système d'appartements pourra être qualifié de "cohérent" ou "double" car c'est le type de système d'appartement qui provient d'un système de Tits double selon la terminologie de \cite{bruhat-tits} 5.1.3.\\

Fixons avant tout un peu de formalisme ($\A$ est toujours le système complet d'appartements):
\begin{defin}
Un système d'appartements pour $(\I,\A)$ est un sous-ensemble $\A'$ de $\A$ tel que deux chambres de $\I$ sont toujours incluses dans un élément de $\A'$.
\end{defin}
La proposition suivante est alors immédiatement vérifiée:
\begin{prop}\begin{itemize}
\item Si $\A'$ est un système d'appartements pour $(\I,\A)$, alors $(\I,\A')$ est un immeuble.
\item Soit $A\in \A$ et $\mathcal B$ une partie bornée de $A$. Alors il existe $A'\in\A'$ contenant $\mathcal B$.
\end{itemize}
\end{prop}

\begin{defin}
Un système d'appartements $\A'$ pour $(\I,\A)$ est dit double, ou cohérent, si pour tous quartiers $Q_1$ et $Q_2$ de $(\I,\A')$, il existe des sous-quartiers $Q_1'$ et $Q_2'$ ainsi qu'un appartement $A\in\A'$ tels que $Q_1'\cup Q_2'\subset A$.
\end{defin}
D'après la proposition \ref{prop: facettes dans appart}, le système complet d'appartements $\A$ est double.\\
Pour tout système double d'appartements, on peut définir l'immeuble sphérique à l'infini $(\I,\A')^\infini$, voir \cite{weiss}.\\

On fixe à présent un système double d'appartements. La proposition suivante est la proposition 10.28 de \cite{weiss} (où un système d'appartements est par définition double, et où un système "plein" est un système où la conclusion de la proposition suivante est vraie), elle affirme que le deuxième point du corollaire \ref{cor:inclusion dans un appartement} est encore vrai dans $(\I,\A')$:
\begin{prop}
 Soit $\alpha$ un demi-appartement de $(\I,\A')$, et $c$ une chambre ayant une cloison dans $\partial \alpha$, alors il existe $A\in\A'$ tel que $\alpha\cup c\subset A$.
\end{prop}
Voici les grandes lignes de la démonstration:\\
Soit $m$ la cloison de $c$ incluse dans $\partial \alpha$, on peut supposer que $m$ contient un sommet spécial $o$. Alors l'immeuble à l'infini $(\I,\A')^\infini$ est isomorphe à l'ensemble des facettes de quartier de sommet $o$, et l'application $p:(\I,\A')^\infini\rightarrow o^*$ est un morphisme surjectif.\\
 Soit $Z$ contenant $\alpha$, soient $x$ et $y$ les chambres de $Z$ contenant $m$, avec $x\in\alpha$, $y\in -\alpha$. Soit $d$ la chambre opposée à $y$ dans $o^*\cap Z$.\\
  Soient $Q_x$, $Q_y$, $Q_d$ les quartiers inclus dans $Z$, de sommet $o$, contenant $x$, $y$ et $d$ respectivement, ces quartiers définissent des chambres à l'infini $x^\infini$, $y^\infini$ et $d^\infini$. L'intersection de $x^\infini$ et $y^\infini$ est la cloison $m^\infini$ de $(\I,\A')^\infini$ correspondant au cône inclus dans $\partial \alpha$, de sommet $o$ et contenant $m$.\\

Nous voulons à présent montrer qu'il existe une chambre $c^\infini\in (m^\infini)^*$ telle que $p(c^\infini)=c$. Ceci provient de \cite{weiss} 29.53, on reproduit ici l'argument.\\
Soit $Y\in \A'$ contenant $x$ et $c$, soit $x'$ l'opposée de $x$ dans $Y$, donc il existe une galerie minimale $\Gamma$ de $x$ à $x'$ passant par $c$, de longueur égale au diamètre $d$ de $o^*$ et de $(\I,\A')^\infini$. Soit $x'^\infini\in (\I,\A')^\infini$ telle que $p(x'^\infini)=x'$. Comme $p$ est un morphisme de complexes de chambres, on a $d(x,x')\leq d(x^\infini, x'^\infini)$, mais comme $d(x,x')=d$ est déjà maximale, on a en fait égalité. Donc $x^\infini$ et $x'^\infini$ sont opposées, et il existe une galerie $\Gamma^\infini$ reliant $x^\infini$ à $x'^\infini$ de même type que $\Gamma$. Cette galerie est alors envoyée par $p$ sur $\Gamma$, par conséquent on peut choisir pour $c^\infini$ sa deuxième chambre.\\

Ensuite, il existe $A$ tel que l'appartement à l'infini $A^\infini$ contient $c^\infini$ et $d^\infini$, donc $A$ contient des scp de $Q_c$ et $Q_d$. Comme ces quartiers sont opposés en $o$, $A$ contient en fait $Q_c\cup Q_d$ et donc en particulier $c$. D'autre part, $A^\infini$ contient toute galerie tendue de $d^\infini$ à $c^\infini$, et il en existe une qui passe par $x^\infini$. Donc $A$ contient un scp de $Q_x$, mais comme $A$ contient le sommet $o$ de $Q_x$, on obtient $Q_x\subset A$. Alors par clôture de $A\cap Z$, $\alpha\subset A$.\cqfd\\

\begin{cor} Le résultat de \ref{prop:inclusion d'une partie et d'une chambre dans un appart} est encore vrai pour un immeuble muni d'un système double d'appartements.
\end{cor}
Les corollaires de la proposition \ref{prop:inclusion d'une partie et d'une chambre dans un appart} sont eux aussi encore vrais pour un système double d'appartements, si on part d'une galerie tendue ou d'une facette de quartier déjà incluse dans un appartement de $\A'$:
\begin{cor}
\begin{itemize}
\item Soit $\Gamma$ une galerie tendue incluse dans un appartement de $\A'$ et $c$ une chambre, alors il existe $A\in\A'$ contenant $c$ et une sous-galerie de $\Gamma$.
\item Soit $f$ une facette de quartier dans un appartement de $\A'$ et $c$ une chambre, alors il existe $A\in\A'$ contenant $c$ et un scp de $f$.
\end{itemize}
\end{cor}

La proposition \ref{prop:cheminees dans appart} reste elle aussi vraie dans un système double d'appartements:
\begin{prop}
Soient $R_1$ et $R_2$ deux cheminées de $(\I,\A')$, alors il existe $A\in\A'$ contenant des sous-cheminées pleines de $R_1$ et $R_2$.
\end{prop}

\demo\\
Soient $A_1,A_2\in\A'$ des appartements contenant $R_1$ et $R_2$, respectivement.
Quitte à raccourcir $R_1$ et $R_2$, il existe un appartement $Z$ dans le système complet $\A$ d'appartements tel que $R_1\cup R_2\subset Z$. Par le lemme \ref{lemme: cheminee de base une chambre}, on peut supposer que $R_1$ et $R_2$ ont pour bases des chambres, disons $R_1=\cl(c_1+\vec f_1)$ et $R_2=\cl(c_2+\vec f_2)$. Soient $\delta_1=s_1+\R^+\vec v_1$ et $\delta_2=s_2+\R^+\vec v_2$ des droites caractéristiques pour $R_1$ et $R_2$; le fait de supposer que les bases de ces cheminées sont des chambres entraîne que $R_i$ contient un voisinage de $\delta_i$, pour $i\in\{1,2\}$. En particulier, $Z\cap A_i$ contient un ouvert, ce qui fournit un isomorphisme canonique entre $\vec Z$ et $\vec A_i$.\\
Notons $M^i_t=s_i+t\vec v_i$, pour $i\in\{1,2\}$ et $t\in\R^+$.
Quitte à déplacer $s_2$ à l'intérieur de $c_2$, on peut supposer que pour $t$ assez grand, $\fl{M^1_tM^2_t}=\fl{s_1s_2}+t(\vec v_2-\vec v_1)$ est dans une chambre $\vec D$ de $\vec Z$. En raccourcissant $\delta_1$ et $\delta_2$, on peut supposer que c'est le cas pour tout $t\in\R^+$.\\
 Soit $\vec D_i$ la chambre correspondant à $\vec D$ dans $\vec A_i$. Alors $s_1+\vec D_1$ et $s_2-\vec D_2$ sont deux quartiers de $(\I,\A')$, et il existe $A\in\A'$ contenant des sous-quartiers $s_1'+\vec D_1$ et $s_2'-\vec D_2$. Montrons que $A$ contient $\delta_1\cup\delta_2$.\\
Soit $t\in\R^+$. Il existe $N^1_t,N^2_t\in ]M^1_t,M^2_t[$ tels que $[M^1_t,N^1_t[\subset A_1\cap Z$ et $]N^2_t,M^2_t]\subset A_2\cap Z$. On peut donc prolonger de manière unique le segment $[M^1_t,M^2_t]$ en un ensemble $D_t$, dans l'appartement $A_1$ d'un côté et dans $A_2$ de l'autre, de sorte que $D_t\cap A_1$ et $D_t\cap A_2$ soient des demi-droites. Comme $\fl{M^1_tN^1_t}\in \vec D_1$ et $\fl{M^2_tN^2_t}\in -\vec D_2$, on voit que $D_t\cap (s_1'+\vec D_1)$ et $D_t\cap(s_2'-\vec D_2)$ sont deux demi-droites (non vides).\\

Montrons que $D_t$ est une droite. On note $D_i=D_t\cap A_i$ et $D_Z=D_t\cap Z$, donc $D_i\cap D_Z$ est un segment non trivial. Soit $i\in\{1,2\}$, on commence par montrer que $D_i\cup D_Z$ est une géodésique.\\
 Soit $c$ une chambre de $Z\cap A_i$ telle que $\bar c$ contient un segment non trivial de $D_i\cap D_Z$, et soit $\rho$ la rétraction $\rho_{c,Z}$. Alors $\rho(D_i\cup D_Z)$ est une demi-droite dans $Z$. 
 Soient $x,y\in D_i\cup D_Z$, notons $l_{D_t}(x,y)$ (resp. $l_{\rho(\D_t)}(\rho(x),\rho(y))$ ) la longueur de la partie de $D_t$ entre $x$ et $y$ (resp. de $\rho(D_t)$ entre $\rho(x)$ et $\rho(y)$). Alors, en utilisant successivement que $\rho(D_i\cup D_Z)$ est une demi-droite, que $\rho$ induit une isométrie sur $D_i$ et sur $D_Z$  et que $\rho$ diminue les distances: $d(\rho(x),\rho(y))=l_{\rho(\D_t)}(\rho(x),\rho(y))=l_{D_t}(x,y)\geq d(x,y)\geq d(\rho(x),\rho(y))$. Par conséquent, $d(x,y)=d(\rho(x),\rho(y)$ pour tous $x,y\in D_i\cup D_Z$. Comme $\rho(D_i\cup D_Z)$ en est une, ceci prouve que $D_i\cup D_Z$ est une géodésique, grâce à l'égalité $[x,y]=\ens{z}{d(x,z)+d(z,y)=d(x,y)}$.\\
 Ainsi, $D_1\cup D_Z$ et $D_2\cup D_Z$ sont des géodésiques. Soit $d$ une chambre de $Z$ telle que $\bar d$ contient un segment non trivial de $D_Z$, alors $\rho_{d,Z}(D_1\cup D_Z\cup D_2)$ est une droite de $Z$. Le même raisonnement que précédemment prouve maintenant que $D_t=D_1\cup D_Z\cup D_2$ est une géodésique, donc une droite.\\

Comme $A$ contient deux demi-droites opposées de $D_t$, il contient $D_t$, donc en particulier $[M^1_t,M^2_t]$. Donc $\delta_1\cup\delta_2\subset A$, d'où $R_1\cup R_2\subset A$.\cqfd\\

\begin{cor}\begin{itemize}
\item Soit $R$ une cheminée de $(\I,\A')$, et $c$ une chambre, alors il existe $A\in\A'$ contenant $c$ et une scp de $R$.
\item Soient $f$ et $g$ deux cônes de direction incluse dans une facette de Weyl, alors il existe $A\in\A'$ contenant un scp de $f$ et un scp de $f$.
\item Soient $\Gamma_1$ et $\Gamma_2$ deux galeries, chacune incluse dans un appartement de $\A'$. Alors il existe $A\in\A'$ contenant des sous-galeries de $\Gamma_1$ et $\Gamma_2$.
\end{itemize}
\end{cor}

\demo Les deux premiers points forment l'analogue de \ref{prop: facettes dans appart}, ils sont clairs. Le troisième point est l'analogue de \ref{th:deux galeries dans appart} à ceci près qu'il faut bien supposer $\Gamma_1$ et $\Gamma_2$ incluse à priori dans des appartements de $\A'$. Pour le prouver, il suffit de montrer le lemme suivant:
\begin{lemme}
Soit $\Gamma=c_0,c_1,...$ une galerie tendue infinie dans un appartement $A$. Alors il existe une cheminée $R$ contenant une sous-galerie de $\Gamma$ et telle que toute scp de $R$ contient une sous-galerie de $\Gamma$.
\end{lemme}
\pv\\
Pour toute racine vectorielle $\vec\alpha$, une et une seule des trois possibilités suivantes est vérifiée (les racines sont vues ici comme des demi-appartements):\begin{enumerate}
\item Il existe une racine $\alpha$ de direction $\vec\alpha$ contenant une sous-galerie de $\Gamma$, et toute racine de direction $-\vec\alpha$ ne contient qu'une nombre fini de chambres de $\Gamma$.
\item Il existe une racine $-\alpha$ de direction $-\vec\alpha$ contenant une sous-galerie de $\Gamma$, et toute racine de direction $\vec\alpha$ ne contient qu'une nombre fini de chambres de $\Gamma$.
\item Il existe une racine $\alpha$ de direction $\vec\alpha$ telle que $\alpha$ et $-\alpha-1$ contiennent des sous-galeries de $\Gamma$. Une telle racine $\alpha$ est alors unique.
\end{enumerate}
Soit $N$ l'ensemble des racines affines obtenues dans le cas $3.$, soit $i_0$ le premier indice tel que $c_{i_0}\subset \bigcap_{\alpha\in N}(\alpha\cap (-\alpha-1))$. Pour chaque racine vectorielle $\vec\beta$ dans le cas $1.$, soit $\beta$ la plus petite racine contenant $\{c_i\}_{i\geq i_0}$, soit $P$ l'ensemble de ces racines. Soit la facette vectorielle $\vec f=\{\vec \alpha=0\et \vec\beta>0,\: \forall \alpha\in N,\: \beta\in P\}$. Alors la cheminée $R:= \cl(c_{i_0}+\vec f)=\bigcap_{\alpha\in N\cup P} \alpha$ convient.\cqfd\\

   \subsection{Parallélisme}
   
   %On suppose toujours $\I$ affine.
 La notion de cheminée permet, grâce à la proposition \ref{prop:cheminees dans appart}, d'étendre la notion de parallélisme à des cônes qui ne sont pas forcément dans le même appartement:
 \begin{defin} Soient $f$ et $g$ deux cônes chacun inclus dans une facette de quartier, de sorte qu'il existe un appartement contenant des scp $f'$ et $g'$ de $f$ et $g$. On dit que $f$ et $g$ sont parallèles lorsque $f'$ est parallèle à $g'$ dans cet appartement.\end{defin}
 
 On vérifie facilement que la définition est indépendante des scp et de l'appartement choisi: si $f$ et $g$ sont parallèles, alors dès que deux scp $f'$ et $g'$ sont inclus dans un même appartement $A$, $f'$ est un translaté de $g'$.\\
 
 Nous savons qu'il existe toujours un appartement contenant des scp de $f$ et $g$ lorsque $f$ et $g$ sont (inclus dans) des facettes de quartiers. Lorsque $f\parallele g$, on peut encore améliorer ce résultat: il existe un appartement contenant un des deux cônes en entier. On peut même remplacer un des deux cônes par une cheminée.

 \begin{prop}\label{prop:parallele} Soient $f,g$ deux facettes de quartier parallèles, et $R=\cl(c+\vec f)$ où $c$ est une chambre $c$ contenant $s(f)$ dans son adhérence (donc $f\subset \bar R$). Alors il existe un appartement contenant $R$ et un scp de $g$.
 \end{prop}
 \demo\\
 Soit $Z$ un appartement contenant une scp $R'=\cl(c'+\vec f)$ de $R$ et un scp $g'$ de $g$. Soit $f'$ un scp de $f$ inclus dans $R'$. Comme $\vec g'=\vec f'$, et que $f'^*_Z$ est un cône d'intérieur non vide, $g'$ doit couper $f'^*_Z$. Soit $y\in g'\cap f'^*_Z$. Soit $\vec C\subset \vec Z$ une chambre vectorielle, contenant $\vec f'$ dans son adhérence et telle que $y\in f'+\vec C$. Alors $y+\vec g'=y+\vec f'\subset f'+\vec C$. Or, par la proposition \ref{prop:cone et chambre dans un appart}, il existe un appartement $X$ contenant $c\cup (f'+\vec C)$, donc en particulier, $y+\vec g'$ est un scp de $g$ inclus dans $X$. Si $A$ est un appartement contenant $f$, la clôture de $A\cap X$ implique que $R=\cl(c\cup (x'+\vec f'))$ est également incluse dans $X$.\cqfd\\

\rema Si $f$ et $g$ ne vérifient pas $\vec f=\vec g$ mais seulement $\vec f\subset \overline{\vec g}$, alors le résultat de la proposition est encore vrai.\\
 
 \begin{prop} Le parallélisme est une relation d'équivalence.\end{prop}
 \demo\\
 Il est évident que le parallélisme est une relation réflexive et symétrique, montrons la transitivité. Soit $f\parallele g$ et $g\parallele h$ trois cônes inclus dans des facettes de quartiers. On peut supposer $f$ et $g$ inclus dans un appartement $A$, et $h$ inclus dans un appartement $B$ contenant un scp $g'$ de $g$.  Montrons que $f\parallele h$. On remarque que puisque le parallélisme à l'intérieur d'un appartement fixé est une relation d'équivalence, s'il existe un appartement contenant un scp de chacun des trois cônes $f$, $g$ et $h$ alors le résultat est vrai.\\
  Soit $x$ le sommet de $f$ et $z$ celui de $h$. Soit $\Gamma=c_1,...,c_l$ une galerie minimale entre une chambre de $A$ contenant $x$ dans son adhérence et une chambre de $B$ contenant $z$ dans son adhérence. On raisonne par récurrence sur la longueur $l$ de $\Gamma$.\\

 Si $l=1$, alors $x$ et $z$ sont dans une même chambre fermée $c\subset A\cap B$. Alors $A\cap B\supset \{x\}\cup g'$. Mais la clôture de ceci dans $A$ contient $f$, d'où $f\subset A\cap B$. Alors $B$ contient $f$, $h$ et $g'$, et dans $B$, $f$ est translaté de $g'$ et $g'$ est translaté de $h$. Ceci implique que $f\parallele h$.\\
 
 Supposons $l>1$. Soit $x'\in c_1\wedge c_2$ et $f'$ le cône de $A$ parallèle à $f$ (et donc à $g$) de sommet $x'$. Soit $A_2$ contenant $c_2$ et un scp $g_2$ de $g$. Alors $A\cap A_2\supset g_2\cup\{x'\}$, d'où par convexité, $f'\subset A_2$. Alors la galerie $c_2,...,c_l$ est de longueur $l-1$, elle part d'une chambre de $A_2$ contenant le sommet de $f'$, et $f'\subset A_2$, nous pouvons utiliser l'hypothèse de récurrence pour obtenir que $f'\parallele h$. Enfin, si $R$ est la cheminée $\cl(c_1+\vec f)$ dans $A$, il existe $Z\supset R\cup h_2$, avec $h_2$ un scp de $h$ (proposition \ref{prop:parallele}). Par convexité de $A\cap Z$, $f'\subset Z$. Et par transitivité du parallélisme à l'intérieur de $Z$, $f\parallele h$.\cqfd\\

   \section{Construction de $\ib$}
   \label{section:ensemble}
   On procède de la même manière que pour la construction de $\overline{A_0}$: l'ensemble $\ib$ sera un ensemble de cônes quotienté par la relation d'équivalence signifiant que deux cônes ont une intersection de dimension maximale. Ensuite ces cônes, légèrement élargis, fourniront une base de voisinage de la topologie.

 \subsection{Préliminaires}
Voici regroupés quelques résultats techniques qui serviront plusieurs fois dans la suite.

\begin{lemme} Soient $A,Z,Z'\in\A$, soient $\phi:Z\isom A$ et $\phi':Z'\isom A$ deux isomorphismes qui coïncident sur une facette $b\subset Z\cap Z'$. Alors il existe $w\in W_A$ qui fixe $\phi(b)$ et tel que $\forall x\in Z\cap Z'$, $\phi(x)=w.\phi'(x)$.\end{lemme}

\pv\\
Soit $\xi:Z'\isom Z$ fixant $Z'\cap Z$. Alors pour tout $x\in Z\cap Z'$, on a $\phi(x)=\phi\circ\xi\circ\phi^{' -1}\circ\phi'(x)$. On pose $w=\phi\circ\xi\circ\phi^{' -1}$: $w$ fixe bien $\phi(b)=\phi'(b)$, donc $w$ convient.\cqfd\\

\begin{prop} \label{lemme:partie stable par une suite de retractions}

%--------------- L'hypothèse 0\in\barre{\vec g} est peut-être inutile----------------

 Soient $A\in\A$, $x\in A$, $\vec\gamma\subset \vec A$ un cône inclus dans une facette de Weyl, et $b$ une partie d'une facette de $A$. On suppose que $b$ est incluse dans $x+\vec \gamma$, contient une base affine de $\aff(x+\vec\gamma)$ et que $x\in\bar b$.\\
 On considère un point $t\in\I$ inclus dans un appartement $Z$ tel que $b\subset A\cap Z$. On suppose que $t$ est envoyé par un isomorphisme $\phi:Z\isom A$ fixant $b$ dans une partie de $A$ de la forme $x+\vec g$, où $\vec g\subset \vec A$ est convexe et incluse dans $\vec\gamma^*_A$. On suppose également que $0\in\barre{\vec g}$ et $\vec g$ est stable par le fixateur de $\vec\gamma$ dans $W_{\vec A}$, qu'on note $W_{\vec \gamma}$.\\
Si $b'$ est une autre partie incluse dans une facette, telle que $\vec\gamma\subset \fl{\aff(b')}$, et $b'\subset b-\vec\gamma^*$, si $Z'$ est un appartement contenant $b'$ et $t$ et si $\phi'$ est un isomorphisme $Z'\isom A$ fixant $b'$, alors $\phi'(t)\in x+\vec g$.\\

%Soit $\vec g\subset \vec A$ une partie convexe, incluse dans $\vec\gamma^*$ et stable par le fixateur de $\vec\gamma$ dans $W_{\vec A}$, qu'on note $W_{\vec \gamma}$. On suppose aussi que $b'\subset b-\vec\gamma^*$.\\
 %Soit enfin $Z\in \A$, $t\in Z$ et $\phi:Z\isom A$ tels que $b\subset A\cap Z$ et $\phi(t)\in x+\vec g$. Alors pour tout appartement $Z'$ contenant $\{t\}\cup b'$ et pour tout $\phi':Z'\isom A$ fixant $b'$, on a $\phi'(t)\in x+\vec g$.
\end{prop}

\rema Par continuité, on obtient directement que si $t\in \bar Z$ vérifie $\phi(t)\in\overline{x+\vec g}$, alors $\phi'(t)\in \overline{x+\vec g}$.\\

Dans le cas où $b$ et $b'$ sont des chambres (c'est le cas essentiel comme on va le voir dès le début de la preuve), cette proposition s'exprime beaucoup plus clairement:
\begin{cor} Soient $x$, $\vec\gamma$ et $\vec g$ comme dans la proposition. Soit $c$ une chambre de $A$ telle que $x\in \bar c$ et $c'$ une autre chambre telle que $c'\subset x-\vec\gamma^*$. Soit $t$ un point de $\I$ tel que $\rho_{c,A}(t)\in x+\vec g$ alors $\rho_{c',A}(t)\in x+\vec g$.\\
\end{cor}

\textit{Démonstration de la proposition}\\
Soit $\Gamma=c_1,...,c_k$ une galerie tendue de $b$ à $b'$, notons $\rho_i=\rho_{A,c_i}$ la rétraction sur $A$ centrée en $c_i$. Comme $b$ contient une base de $x+\vec\gamma$, $x+\vec g$ est stable par $\fix_{W_A}(b)$. On déduit alors du lemme précédent que $\rho_1(t)\in x+\vec g$, et que le résultat sera prouvé si nous prouvons que $\rho_k(t)\in x+\vec g$. Montrons par récurrence que pour tout $l\in\llbracket 1,k\rrbracket $, $\rho_l(t)\in x+\vec g$, le cas $l=1$ est déjà vu.\\
Soit $l\in\llbracket 2,k\rrbracket $, supposons $\rho_{l-1}(t)\in x+\vec g$. Soit $M_l$ le mur de $A$ contenant $c_l\wedge c_{l-1}$. Si $M_l$ sépare, pas forcément strictement, $c_l$ et $\rho_{l-1}(t)$, alors on est dans la disposition: $c_l\;|^{M_l}\;c_{l-1}, \, \rho_{l-1}(t)$. Alors la proposition \ref{prop:retractions adjacentes}  affirme que $\rho_l(t)=\rho_{l-1}(t)\in x+\vec g$.\\

Dans l'autre cas, on a $b',\, c_l,\,\rho_{l-1}(t)\;|^{M_l}\; c_{l-1},\, b$, et $\rho_{l-1}(t)\not\in M_l$. Notons $\sigma_l$ la réflexion selon $M_l$, les deux possibilités $\rho_l(t)=\rho_{l-1}(t)$ ou $\rho_l(t)=\sigma_l\circ\rho_{l-1}(t)$ sont alors possibles. Dans le premier cas, on a évidemment $\rho_l(t)\in x+\vec g$, étudions le second cas. On note $y=\rho_{l-1}(t)$, donc $\rho_l(t)=\sigma_l(y)$.\\
 Soit $\vec\alpha\in \vec A^*$ tel que $\vec M_l=\ker(\vec\alpha)$. Comme $M_l$ est un mur d'une galerie minimale entre $b$ et $b'$, il ne peut contenir $b'$, donc il existe $u\in b'$ tel que $\vec\alpha(\vec{xu})\not=0$, on peut supposer $\vec\alpha(\vec{xu})>0$. Par ailleurs, avec ce choix de signe pour $\vec\alpha$, on a $\vec\alpha(\vec{xy})>0$. Mais $\vec{xu}\in -\vec\gamma^*$ et $\vec{xy}\in\vec\gamma^*$. Ceci entraîne que $\vec M_l$ coupe $\vec\gamma^*$, et donc que $\vec\gamma\subset\vec M_l$. En particulier, $\vec\sigma_l\in W_{\vec\gamma}$ et donc $\vec\sigma_l(\vec g)=\vec g$.\\
Soit $z$ l'unique point de l'intersection $[x,y]\cap M_l$ (le cas où $\{x,y\}\subset M_l$ est trivial, on a alors $\rho_l(t)=\rho_{l-1}(t)$). Par convexité de $\vec g$ et comme $0\in\barre{\vec g}$, $z\in x+\barre{\vec g}$.  Nous savons que $\vec\sigma_l(\vec{xy})\in \vec g$, et donc que $x+\vec\sigma_l(\vec{xy})\in x+\vec g$. Ensuite $\sigma_l(y)=z+\vec\sigma_l(\vec{zy})=z+\lambda\vec\sigma_l(\vec{xy})$ pour le scalaire $\lambda\in[0,1]$ adéquat. C'est alors une simple application du théorème de Thalès que de vérifier que $\sigma_l(y)\in [y,x+\vec\sigma_l(\vec{xy})]$, d'où $\sigma_l(y)=\rho_l(t)\in x+\vec g$ par convexité.\cqfd\\

   \begin{lemme} \label{lemme:scp dans le coeur} Soit $Z$ un appartement, $f,g\in\F_Z$ avec $g$ un scp de $f$. Alors $g$ coupe le coeur $\delta(f)$, et donc il existe un scp $g'$ de $g$ dont le coeur est un scp de $\delta(f)$.\end{lemme}
 \pv\\
  Soit $x=s(f)$ et $y=s(g)$ les sommets des deux cônes. Soit $W_{\vec f}\subset Gl(\vec Z)$ le stabilisateur de $\vec f$ dans le groupe de Weyl de $\vec Z$. On définit une action de $W_{\vec f}$ sur $Z$ en imposant que $x$ est un point fixe. Alors $\delta(f)=f\cap \fix_Z(W_{\vec f})$. Soit $z=x+\sum_{w\in W_{\vec f}} w(\vec{xy})$. Le point $z$ est bien $W_{\vec f}-fixe$. Comme $\vec{xy}\in \vec f$, il est dans $f$, et l'écriture $z=x+\vec{xy}+\sum_{w\in W_{\vec f}-\{e\}} w(\vec {xy})$ prouve qu'il est aussi dans $g$. donc $z\in \delta(f)\cap g$, ce qui prouve le lemme.\cqfd\\

 \subsection{Cônes dans l'immeuble}
 
 \begin{defin} Soit $A$ un appartement, et $f\in \F_A$. Le cône de $\I$ correspondant au cône $f$ de $A$ est :
 \[\tilde f= \bigcup_{B\in \A, \delta(f)\subset B} \phi_B(f)\]
 où $\phi_B$ est un isomorphisme quelconque de $A$ sur $B$ fixant $\delta(f)$.\\
 On note $\F_\I$ l'ensemble de tous ces cônes pour tous les choix de $A$ et de $f$ possibles.
 \end{defin}
 
 Le choix de $\phi_B$ dans la définition n'importe pas puisque deux choix diffèrent d'un automorphisme de $A$ fixant $\delta(f)$ donc stabilisant $f$.\\
 Il est clair que si $f\in\F_A$ et $g\in\F_B$ sont des cônes tels que $\delta(f)=\delta(g)$, alors $\tilde f=\tilde g$: le coeur du cône affine suffit à déterminer le cône de l'immeuble. Dit autrement, si $B$ contient $\delta(f)$, alors $\tilde f=\tilde g$ où $g$ est le cône de $B$ de coeur $\delta(f)$. Nous verrons plus tard que réciproquement, un cône d'immeuble définit un unique coeur de cône affine.\\

   \begin{lemme}  Soit $F=\tilde f\in\F_\I$, avec $f\in\F_A$. Soit $Z$ un appartement contenant $\delta(f)$, et $\phi:A\rightarrow Z$ un isomorphisme fixant $\delta(f)$. Alors $Z\cap \tilde f=\phi(f)$.\end{lemme}
   \pv\\
   L'inclusion $\supset$ est vraie par définition de $\tilde f$. Montrons $\subset$.\\
   Soit $y\in Z\cap \tilde f$. Par définition, il existe un appartement $Y$ contenant $y$ et $\delta(f)$ et tout isomorphisme de $Y$ sur $A$ fixant $\delta(f)$ envoie $y$ dans $f$. Soit $\psi:Y\rightarrow Z$ un isomorphisme fixant $y$ et $\delta(f)$. Alors $\phi\inv\circ\psi:Y\rightarrow A$ est un isomorphisme fixant $\delta(f)$. Donc $\phi\inv\circ\psi(y)\in f$, donc $y=\psi(y)\in\phi(f)$.\cqfd\\

   \begin{lemme} \label{lemme:inclusion} Soient $f,g\in\F_A$ tels que $f$ est un scp de $g$. Alors $\tilde f\subset \tilde g$.\end{lemme}
   \pv\\

 Remarquons que pour montrer qu'un point $t\in\tilde f$ appartient à $\tilde g$, il suffit de trouver un appartement contenant $\{t\}\cup\delta(f)\cup \delta(g)$.\\

   Soit $x$ le sommet de $f$ et $y$ celui de $g$. Il n'y a qu'un nombre fini de murs parallèles à $\delta(\vec f)$ séparant strictement $x$ et $y$, soit $n$ ce nombre. On raisonne par récurrence sur $n$.\\
    Si $n=0$, soit $t\in \tilde f$, et $Y$ un appartement contenant $\{t\}\cup \delta(f)$. Comme $f\subset \delta(f)^*_A$, le lemme précédent indique que $\tilde f\cap Y\subset \delta(f)^*_Y$, en particulier, $t\in \delta(f)^*_Y$. Dès lors, on peut trouver une chambre $\vec C$ de $\vec Y$ telle que $\delta(f)\cup\{t\}\subset \overline{x+\vec C}$ et la proposition \ref{prop:cone et chambre dans un appart} prouve qu'il existe $Z\in\A$ contenant $\{t\}\cup\delta(f)\cup\{y\}$. Par convexité et fermeture de $Z\cap A$, $Z$ contient aussi $x+\delta(\vec f)=\delta(g)$ (car $\vec f=\vec g$ dans $\vec A$). On conclut que $t\in\tilde g$.\\

  Si $n>0$, soit $M$ un mur parallèle à $\delta(\vec f)$ séparant strictement $x$ et $y$. Soit $z=[x,y]\cap M$, alors $z+\vec f$ est scp de $g$, et $f$ est scp de $z+\vec f$. De plus, le nombre de murs parallèles à $\delta(\vec f)$ séparant strictement $y$ et $z$ est strictement inférieur à $n$, de même que le nombre de murs parallèles à $\delta(\vec f)$ séparant strictement $z$ et $x$. Donc par hypothèse de récurrence, $\widetilde{z+\vec f}\subset \tilde g$ et $\tilde f\subset \widetilde{z+\vec f}$, d'où le résultat.\cqfd\\
  
  \rema Il n'est pas vrai que $\tilde f\subset \tilde g$ si $f\subset g$ sans que $\vec f=\vec g$.\\

   \subsection{Coeur d'un cône d'immeuble}

   \begin{lemme}
   Soit $f\in\F_A$, et $\gamma$ un scp de $\delta(f)$. Soit $B$ un appartement contenant $\gamma$, soit $\xi:B\rightarrow A$ un isomorphisme fixant $\gamma$. Soit $x\in B\cap \tilde f$, alors $\xi(x)\in f$. Plus précisément, si $Z$ contient $x$ et $\delta(f)$, alors il existe un isomorphisme de $Z$ sur $A$, fixant $\delta(f)$ et envoyant $x$ sur $\xi(x)$.\end{lemme}
   
   \pv\\
   Soit $Z\supset\{x\}\cup\delta(f)$. Soit $\phi:Z\rightarrow B$ un isomorphisme fixant l'intersection $Z\cap B$. Alors $\xi\circ\phi$ est un isomorphisme de $Z$ sur $A$ qui vérifie $\xi\circ\phi(x)=\xi(x)$. Il suffit donc de prouver que $\xi\circ\phi$ fixe $\delta(f)$ pour prouver le lemme. Mais $\xi\circ\phi$ fixe $\gamma$, donc $\fl{\xi\circ\phi}$ fixe $\vec  \gamma=\delta(\vec f)\subset\vec Z\cap\vec A$. Donc $\xi\circ\phi$ fixe $Z\cap A\cap (\gamma+\vect(\delta(\vec f)))$, ce qui contient $\delta(f)$.\cqfd\\

\begin{prop} Si $f\in \F_A$ et $g\in \F_B$ sont deux cônes tels que $\tilde f=\tilde g$, alors $\delta(f)=\delta(g)\subset A\cap B$.\end{prop}

\demo\\
D'après la proposition \ref{prop:cheminees dans appart}, il existe un appartement $Z$ contenant des scp $\delta_1$ et $\delta_2$ de $\delta(f)$ et $\delta(g)$. Soient $f',g'\in \F_Z$ les cônes engendrés par $\delta_1$ et $\delta_2$. Nous allons commencer par montrer que $\vec g'=\vec f'$.\\

Supposons par l'absurde que ce n'est pas le cas, donc $\vec f'\cap\vec g'=\vide$. Il est alors impossible que $\vec f'\subset \overline{\vec g'}$ et $\vec g'\subset \overline{\vec f'}$, d'après les hypothèses mises sur $\F$. on peut donc supposer qu'il existe $\vec v\in \vec g'\setminus\overline{\vec f'}$.\\
Soit $y$ le sommet de $g$, $y'$ celui de $g'$. Soit $\phi:Z\isom B$ un isomorphisme fixant $\delta_2$, alors $\phi(g')$ est le cône de $B$ engendré par $\delta_2$, c'est donc un scp de $g$. Par le lemme \ref{lemme:inclusion}, $\widetilde {\phi(g')}\subset \tilde g$, mais $\widetilde {\phi(g')}=\widetilde{g'}$, d'où $\widetilde {g'}\subset \tilde g$.

% Il existe un appartement $Y$ contenant $y$ et $y'+(\vec \delta_2)^*_Z$, $Y$ contient donc $g'$, ainsi que $\delta(g)$ par convexité de $Y\cap B$. Soit $\phi:B\rightarrow Y$ un isomorphisme fixant $\delta(g)$, on a alors $g'\subset \phi(g)$. En utilisant le lemme précédent, on arrive à $\tilde g'\subset \tilde{\phi(g)}=\tilde g$.

 De la même manière, soit $\psi:Z\rightarrow A$ un isomorphisme fixant $\delta_1$. Alors $\psi(f')$ est un scp de $f$, et on obtient que $\tilde f'\subset \tilde f=\tilde g$.\\
Nous déduisons en particulier que $\forall \lambda\in \R^+$, $y'+\lambda\vec v\in \tilde g$.  Et comme $\psi(f')$ est un scp de $f$, $\vec \psi(\vec f')=\vec f$. Pour tout réel $\lambda$, $\psi(y'+\lambda\vec v)=\psi(y')+\lambda\vec\psi(\vec v)$. Comme $\vec v\not\in\overline{\vec f'}$, on a $\vec\psi(\vec v)\not\in \overline{\vec f}$. Donc pour lambda assez grand, $\psi(y'+\lambda \vec v)\not\in f$.\\
Comme $y'+\lambda\vec v\in\tilde g=\tilde f$, il existe un appartement $X$ contenant $\delta(f)\cup\{y'+\lambda\vec v\}$. Soit $\eta:X\rightarrow Z$ un isomorphisme fixant $X\cap Z$, en particulier $\eta$ fixe $\delta_1$ et $y'+\lambda\vec v$. Alors $\psi\circ\eta:X\rightarrow A$ est un isomorphisme qui fixe $\delta_1$, et donc $\aff_X(\delta_1)\cap A$, ce qui contient $\delta(f)$. Nous savons dans cette situation que $X\cap\tilde f=(\psi\circ\eta)\inv(f)$, d'où $y'+\lambda\vec v\in(\psi\circ\eta)\inv(f)$, mais $\psi\circ\eta(y'+\lambda\vec v)=\psi(y'+\lambda\vec v)\not\in f$, d'où la contradiction.\\

 Nous avons donc prouvé que $\vec f'=\vec g'$, nous allons maintenant montrer que les sommets $x$ et $y$ de $f$ et $g$ coïncident. En utilisant la proposition \ref{prop:parallele}, on peut supposer que $Z$ contient $\delta(f)$ en entier, donc $\tilde f=\tilde f'$, on peut même supposer que $Z=A$. On peut également supposer que $B$ contient un scp $\gamma_1$ de $\delta(f)$. En résumé, $A\cap B$ contient un scp $\gamma_1$ de $\delta(f)$ et un scp $\delta_2$ de $\delta(g)$, et ces deux scp sont parallèles. L'égalité $\delta(\vec f)=\delta(\vec g)$ a donc un sens (et est vraie), dans $\vec A\cap \vec B$.\\
  Soit $\xi:B\rightarrow A$ un isomorphisme fixant $A\cap B$. D'après le lemme précédent, $\xi(y)\in f$ et $\xi\inv(x)\in g$. La seconde égalité implique $x\in \xi(g)=\xi(y)+\vec\xi(\vec g)$. Or $\vec\xi(\vec g)$ est un cône vectoriel de coeur $\vec\xi(\delta(\vec g))=\delta(\vec f)$, puisque $\xi$ fixe $\delta_2$ et $\delta(\vec g)=\delta(\vec f)$. Donc $\vec\xi(\vec g)=\vec f$, et $x\in\xi(y)+\vec f$.\\
  Nous arrivons ainsi aux deux relations: $x\in\xi(y)+\vec f$ et $\xi(y)\in x+\vec f$. Par unicité du sommet de $f$, ceci implique $x=\xi(y)$. On utilise enfin la version précise du lemme précédent: soit $T$ un appartement contenant $\delta(f)$ et $y$, il existe un isomorphisme $\zeta:T\rightarrow A$ fixant $\delta(f)$ et tel que $\zeta(y)=\xi(y)$. Mais $\zeta\inv(x)=x$ puisque $x\in\overline{\delta(f)}$, d'où $x=y$.\\
 En particulier, $x=y\in A\cap B$, et par convexité de $A\cap B$, $\delta(f)\subset A\cap B$ et $\delta(g)\subset A\cap B$, et enfin on peut calculer, dans $A$ ou dans $B$: $\delta(f)=x+\delta(\vec f)=y+\delta(\vec g)=\delta(g)$.\cqfd\\

Cette proposition permet les définitions suivantes:
\begin{defins}\begin{itemize} 
\item Le coeur d'un cône $\tilde f$ de l'immeuble $\I$ est $\delta(f)$. On le note $\delta(\tilde f)$. 
\item On appelle sommet de $\tilde f$ le sommet de $\delta(f)$, on le note $s(\tilde f)$.
\item Deux cônes $F,G\in\F_\I$ sont dit parallèles lorsque $\delta(F)\parallele \delta(G)$, ceci définit une relation d'équivalence sur $\F_\I$.
\item Si $G,F\in\F_\I$ sont deux cônes parallèles avec $F\subset G$, on dit que $F$ est un sous-cône parallèle de $G$, et on abrège sous-cône parallèle en scp. Ceci définit une relation d'ordre sur $\F_\I$.
\end{itemize}
\end{defins}

\begin{lemme}\label{lemme:coeurs dans appart}
   Soient $F,G\in\F_\I$ tels que $G$ est scp de $F$. Alors il existe un appartement contenant $\delta(F)\cup\delta(G)$.
   \end{lemme}
  
  \pv\\
 Nous savons déjà  qu'il existe (proposition \ref{prop:parallele}) un appartement $A_1$ contenant $\delta(G)$ et un scp $\delta$ de $\delta(F)$. Il existe aussi $A_2$ contenant $\delta(F)\cup\{s(G)\}$, ceci car $s(G)\in F$. Alors $A_1\cap A_2\supset \{s(G)\}\cup \delta$, d'où par convexité, $A_1\cap A_2\supset \delta(G)$, et donc $A_2\supset \delta(F)\cup\delta(G)$.\cqfd\\

   \subsection{Équivalence de cônes, l'ensemble $\ib$}

   Nous définissons maintenant une relation d'équivalence sur l'ensemble des cônes de l'immeuble d'une manière tout à fait semblable à la définition de la relation d'équivalence sur l'ensemble des cônes affines de $A_0$.
   \begin{defin} Deux cônes $F,G\in\F_\I$ sont équivalents lorsqu'il existe $H\in\F_\I$ qui est un scp de $F$ et de $G$. On note alors $F\sim G$.
   \end{defin}
   
   Le lemme \ref{lemme:inclusion} montre que si $f,g\in\F_A$ alors $f\sim_A g\implique \tilde f\sim \tilde g$.\\

   \begin{lemme} Soient $F,G\in\F_\I$, alors $F\sim G$ si et seulement si $F\parallele G$ et $F\cap G\not=\vide$.
   \end{lemme}
   
   \pv\\
   Si $F\sim G$, il est clair que $F\cap G\not=\vide$ et que $F\parallele G$ (car $F\parallele H\parallele G$).\\
   Supposons maintenant que $F\cap G\not=\vide$ et que $F\parallele G$. Soit $x\in F\cap G$, soit $A$ et $B$ des appartements tels que $\delta(F)\cup\{x\}\subset A$ et $\delta(G)\cup\{x\}\subset B$. Nous pouvons alors construire deux cônes: $\delta_1:=x+\fl{\delta(F)}\subset A$ et $\delta_2:=x+\fl{\delta(G)}\subset B$. Ces deux cônes sont parallèles, donc il existe $Z$ contenant $\delta_1$ et un sous-cône parallèle $\delta_2'$ de $\delta_2$. Par convexité de $Z\cap B$, on voit qu'en fait $\delta_2\subset Z$. Donc $\delta_1$ et $\delta_2$ sont deux cônes parallèles, inclus dans $Z$ et de même sommet: $\delta_1=\delta_2$. Soit $H\in\F_\I$ tel que $\delta(H)=\delta_1$. Il est à peu près clair que $H\subset F\cap G$, mais je détaille quand même le raisonnement:\\
   Soit $f\in\F_A$ tel que $\delta(f)=\delta(F)$. Alors $x\in f$, puis $\delta_1\subset f$. Soit $h_1\in\F_A$ tel que $\delta(h_1)=\delta_1$, alors $H=\tilde{h_1}$ et $h_1\subset f$. Comme $F=\tilde f$, on conclut par le lemme \ref{lemme:inclusion} que $H\subset F$. On procède de même pour prouver que $H\subset G$.\cqfd\\

   \begin{prop} La relation "être équivalent" est une relation d'équivalence sur $\F_\I$.\end{prop}
   
   \demo\\
   Cette relation est symétrique et réflexive, montrons qu'elle est transitive. Soient $F_1,F_2,F_3\in\F_\I$ tels que $F_1\sim F_2$ et $F_2\sim F_3$. Soit $G_1$ un scp de $F_1$ et $F_2$, et $G_2$ un scp de $F_2$ et $F_3$. Il nous suffit de trouver un cône $H\in\F_\I$ qui soit un scp de $G_1$ et de $G_2$. Et comme nous savons déjà que $G_1\parallele G_2$, grâce au lemme précédent, il suffit de prouver que $G_1\cap G_2\not=\vide$.\\

   Quitte à remplacer $G_1$ par un scp, il existe un appartement $Z_1$ contenant $\delta(F_2)\cup\delta(G_1)$. Le lemme \ref{lemme:scp dans le coeur} appliqué à $F_2\cap Z_1$ et $G_1\cap Z_1$ dans l'appartement $Z_1$ prouve l'existence d'un scp $H_1$ de $G_1$ dont le coeur est un scp de $\delta(F_2)$. De même, il existe un scp $H_2$ de $G_2$ dont le coeur est aussi inclus dans $\delta(F_2)$.\\
   Ainsi $\delta(H_1)$ et $\delta(H_2)$ sont deux scp du même cône $\delta(F_2)$, leur intersection contient donc un troisième scp de $\delta(F_2)$ (il suffit de vérifier que $s(\delta(F_2))+ \fl{s(\delta(F_2)s(\delta(H_1))} +\fl{s(\delta(F_2)s(\delta(H_2))}\in \delta(H_1)\cap\delta(H_2)$ ). Si $H$ est le cône de $\I$ correspondant à ce coeur, alors $H$ est un scp de $H_1$ donc de $G_1$ donc de $F_1$, ainsi que de $H_2$ donc de $G_2$ donc de $F_3$. Ce qui prouve que $F_1\sim F_3$.\cqfd\\
   
   Nous pouvons enfin définir l'ensemble $\ib$:
   \begin{defin} On pose $\ib=\F_{\I} / \sim $. Si $F\in \F_\I$, on note $[F]$ la classe d'équivalence de $F$ pour $\sim$. \end{defin}

Notons enfin ce petit résultat:\\
   \begin{lemme} \label{lemme:scp de G dans A} Soit $A\in\A$, et $f\in \F_A$. Soit $G\in F_\I$ tel que $\tilde f\sim G$. Alors il existe $f'$ un scp de $f$ tel que $\tilde f'$ est un scp de $G$.\end{lemme}

\demo D'après la proposition \ref{prop:parallele}, il existe un appartement $Z$ contenant $\delta(f)$ et un scp de $\delta(G)$. Soit $g'$ le cône affine de $Z$ engendré par ce scp. D'après le lemme \ref{lemme:scp dans le coeur}, il existe $\delta'$ un scp de $\delta(f)$ qui est inclus dans $g'$. Alors le cône $f'$ engendré par $\delta'$ dans $A$ convient.\cqfd\\

   \subsection{Injections canoniques}
   
   \begin{defins}\begin{itemize}
   \item Soit \fonc{\iota_\I}{\I}{\ib}{x}{[\{x\}]}. C'est l'injection canonique de $\I$ dans $\ib$.
   \item  Pour tout appartement $A$, soit \fonc{\iota_A}{\bar A}{\ib} {\left[f\right]_A} {[\tilde f]}. C'est l'injection canonique de $\bar A$ dans $\overline{\I}$.
   \end{itemize}
   \end{defins}
   
   Quelques mots pour s'assurer que ces définitions sont licites:\begin{itemize}
   \item Pour $\iota_{\I}$, il faut vérifier que $\{x\}\in\F_\I$, $\forall x\in\I$. C'est le cas car si $x$ est dans l'appartement $A$, alors $\{x\}\in \F_A$, et $\tilde{\{x\}}=\{x\}$.
   \item Nous avons déjà vu que si $f\sim_A g$ alors $\tilde f\sim \tilde g$ ce qui montre que $\iota_A$ est bien définie.\end{itemize}
   
\begin{prop} Les injections canoniques sont injectives.\end{prop}
\demo\\
Il est clair que $\iota_\I$ est injective. Soit $A\in\A$, pour montrer que $\iota_A$ est injective, on considère deux cônes $f,g\in\F_A$ tels que $\iota_A([f]_A)=\iota_A([g]_A)$. C'est à dire que $\tilde f\sim \tilde g$. Donc il existe $H\in\F_\I$ qui est scp de $\tilde f$ et $\tilde g$. Par le lemme \ref{lemme:coeurs dans appart}, il existe $B\in\A$ contenant $\delta(f)\cup\delta(H)$. En utilisant le lemme \ref{lemme:scp dans le coeur} dans $B$, on obtient un nouveau cône $H'\in\F_\I$ qui est scp de $H$ et dont le coeur est inclus dans $\delta(f)$, donc dans $A$. Alors notant $h'=H'\cap A$, on voit que $h'$ est scp dans $A$ de $f$ et $g$, et donc $f\sim_A g$.\cqfd\\

   \section{Topologie sur $\overline {\mathcal I}$}
   \label{section:topologie}

   \subsection{Définition}

\begin{defin} Soit $F\in\F_\I$, soit $A\in\A(\I)$ contenant $\delta(F)$. On note $b(F)$ et on appelle base de $F$ la base de $\delta(F)$. De même, si $f\in\F_A$, on notera $b(f)$ la base de $\delta(f)$.\end{defin}
La définition de la base d'un cône dont la direction est dans une facette vectorielle est donnée en \ref{defin:base}.\\

 Maintenant, nous avons besoin d'un moyen de définir un voisinage de $0$ dans les espaces directeurs de tous les appartements à la fois. Nous avons déjà défini $\U$ comme étant l'ensemble des voisinages de $0$ dans $\fl{A_0}$, où $A_0$ est l'appartement de référence choisi au début, stables par le groupe de Weyl de $\fl{A_0}$, $W{\fl{A_0}}$. Alors si $U\in\U$, si $\phi:A_0\isom B$ est un isomorphisme quelconque, l'image $\vec\phi(U)\subset \vec B$ est indépendante de l'isomorphisme $\phi$ choisi. Nous pouvons ainsi identifier $U$ à un voisinage de $0$ dans l'espace directeur de n'importe quel appartement. Et si $\phi:A\isom B$ est un isomorphisme, nous pouvons écrire $\phi(x+U)=\phi(x)+U$, $\forall x\in A$.\\

\begin{defin}

 Soit $F\in\F_\I$, $A$ tel que $\delta(F)\subset A$, $f=A\cap F$, et $U\subset \vec A$, stable par $W(\vec A)$. On pose
 \[ \widehat{f+U}=\bigcup_{B\supset b(F)}\phi_B((f+U)\cap \delta(f)^*_A) \]
 \end{defin}

 L'union est prise sur tous les appartements $B$ qui contiennent $b(F)$, et $\phi_B:A\isom B$ est un isomorphisme qui fixe $b(F)$.\\
\rema Il aurait pu paraître plus naturel de prendre l'union des $\phi_B(f+U)$, mais le fait de prendre l'intersection avec $\delta(f)^*_A$ permettra d'utiliser la proposition \ref{prop:cone et chambre dans un appart}.\\

 Le sous groupe de $W_A$ qui fixe $b(f)$ fixe aussi $\delta(f)$, puisque $b(f)$ contient une base affine de $\aff_A\delta(f)$, et donc il stabilise $f$, et sa partie vectorielle est un sous-groupe de $W^v_A$ et donc stabilise $U$. Ceci permet de vérifier que le choix de $\phi_B$ n'importe pas.\\
Si $A'$ est un autre appartement contenant $\delta(F)$, et si $f'=A'\cap F$ alors $\widehat{f+U}=\widehat{f'+U}$. Ceci permet la définition:
 \[\widehat{F+U}:=\widehat{f+U}\].\\

\begin{lemme} \label{lemme:intersection chapeau appartement} On garde les notations précédente. Si $B$ est un appartement contenant $b(F)$, et si $\phi_B:A\isom B$ est un isomorphisme qui fixe $b(F)$, alors $B\cap \widehat{(f+U)}=\phi_B((f+U)\cap \delta(f)^*_A)$. En particulier, $A\cap \widehat{(f+U)}=(f+U)\cap \delta(f)^*_A$.\end{lemme}
\pv Si $x\in B\cap \widehat{(f+U)}$, alors il existe un appartement $C$ contenant $b(f)$ et un isomorphisme $\phi_C:A\isom C$ fixant $b(f)$ tel que $x\in \phi_C((f+U)\cap \delta(f)^*_A)$. On choisit alors $\psi:C\isom B$ fixant $b(f)\cup\{x\}$, comme $\psi\circ\phi_C:A\isom B$ fixe $b(f)$, on a $\psi\circ\phi_C((f+U)\cap \delta(f)^*_A)=\phi_B((f+U)\cap \delta(f)^*_A)$. D'où $x\in \phi_B((f+U)\cap \delta(f)^*_A)$.\cqfd\\

On peut maintenant définir les ensembles qui seront la base de voisinage de la topologie sur $\ib$: 
 \begin{defin}
 Soit $F\in\F_\I$, et $U\in\U$. On pose
 \[ \V_\I(F,U)=\{ x\in\overline{\I} \;| \; \text{un représentant de x est inclus dans }  \widehat{F+U}\} \]
 Lorsqu'il n'y a pas de risque de confondre $\V_\I (F,U)$ avec un voisinage $\V_A(f,U)$ d'un point dans un appartement compactifié, on pourra juste noter $\V(F,U)$.\\
 Pour $x\in\overline{\I}$, on note $\V_x=\ens{\V(F,U)}{x=[F],\; U\in\U}$.
 \end{defin}

 \begin{lemme} \label{lemme:chapeaux} Soient $F,G\in\F_\I$ et $U,V\in \U$ alors:
 
 \begin{itemize}
 \item $\tilde f\subset \widehat{f+U}$
 \item si $\widehat{F+U}\subset\widehat{G+V}$, alors $\V(F,U)\subset \V(G,V)$
 \item si $U\subset V$, alors $\widehat{F+U}\subset\widehat{F+V}$
 \item si $F$ est scp de $G$, alors $\widehat{F+U}\subset\widehat{G+U}$.
 \end{itemize}
 
 \end{lemme}
\pv\\
Les trois premiers points sont évidents.\\
 Pour le quatrième point, on choisit un appartement $A$ contenant $\delta(F)$ et  $\delta(G)$ (lemme \ref{lemme:coeurs dans appart}). On note $f=A\cap F$ le cône engendré par $\delta(F)$ dans $A$, et $g=A\cap G$ le cône engendré par $\delta(G)$, donc $f$ est un scp de $g$.  Il ne reste plus qu'à appliquer la même méthode que celle utilisée dans cette situation pour prouver que $\tilde f\subset \tilde g$ (en \ref{lemme:inclusion}).\\

On commence par supposer qu'aucun mur de $A$, parallèle à $\delta(G)$ ne sépare strictement $s(f)$ et $s(g)$. Soit $x\in \widehat{F+U}$, il existe $B\supset \{x\}\cup b(F)$ et $\phi:A\isom B$ fixant $b(F)$, tel que $x\in \phi((f+U)\cap f^*_A)$. En particulier, $x\in \phi(f)^*_B$. D'après la proposition \ref{prop:cone et chambre dans un appart}, il existe un appartement $Z$ contenant $b(G)\cup b(F)\cup\{x\}$. Soit $\psi:Z\isom A$ fixant $A\cap Z$, donc en particulier $b(F)\cup b(G)$. Comme $x\in \widehat{F+U}$, $\psi(x)\in (f+U)\cap f^*_A\subset (g+U)\cap g^*_A$, donc $x\in \widehat{G+U}$.\\

S'il existe un mur $M$ parallèle à $\delta(G)$ et séparant  strictement $s(f)$ et $s(g)$, soit $z=[s(f),s(g)]\cap M$. Par récurrence, on a $\widehat{F+U}\subset \widehat{\widetilde{z+\vec f}+U}$ et $\widehat{\widetilde{z+\vec f}+U}\subset \widehat{G+U}$.\cqfd\\

 \begin{prop} Il existe une unique topologie sur $\ib$ telle que pour tout $x\in\ib$, $\V_x$ soit une base de voisinages de $x$.
 \end{prop}
 
 \demo\\
 Il faut montrer pour tout $x\in\ib$, que  $\V_x\not=\vide$, que $\forall V\in\V_x$, $x\in V$, et que $\forall x\in\ib$, $V_1,V_2\in\V_x$, $\exists V_3\in\V_x$ tel que $V_3\subset V_1\cap V_2$. Les deux premières assertions sont claires.\\
 Soit donc $x\in\ib$ et $V_1,V_2\in\V_x$. Soient $F_1,F_2\in x$ et $U_1,U_2\in \U$ tels que $V_1=\V(F_1,U_1)$ et $V_2=\V(F_2,U_2)$.\\
 Comme $F_1\sim F_2$, il existe $H$ qui est scp de $F_1$ et de $F_2$. Soit $U=U_1\cap U_2$, on a bien $U\in\U$. Alors, $\V(H,U)\in\V_x$ et d'après le lemme, $\V(H,U)\subset \V(F_1,U_1)\cap \V(F_2,U_2)$.\cqfd\\

 Notons tout de suite ce résultat évident, mais important:\\
\begin{prop}\label{prop:aut de i}
Soit $g\in \aut(\I)$. Alors $g$ se prolonge en un homéomorphisme de $\ib$, en posant $g( [F])=[g(F)]$.\end{prop}

\subsection{Lien avec la topologie de $\bar A$}

 Nous allons maintenant faire le lien entre les voisinages du type $\V(F,U)$ dans l'immeuble et les $\V_A(f,U)$ dans un appartement. La seule petite difficulté provient du fait que dans la définition de $\chap{f+U}$ on fait intervenir $(f+U)\cap f^*$ et non juste $f+U$ comme dans les $\V_A(f,U)$.\\

\begin{defin} \label{defin:Vprime} Soit $A\in\A$, $U\in\U$ et $f\in\F_A$, on note \[\V'_A(f,U)=\ens{ x\in \bar A}{ \text{un représentant de $x$ est inclus dans } (f+U)\cap f^*}\]
\end{defin}

\begin{prop} Soit $A\in\A$. Alors l'ensemble des $\V'_A(f,U)$, pour $f\in\F_A$ et $U\in\U$ est une base de voisinage de la topologie sur $\bar A$\end{prop}
\demo Déjà, quels que soient $f\in\F_A$ et $U\in\U$, on a $\V'_A(f,U)\subset \V_A(f,U)$. Il reste à vérifier que $\V'_A(f,U)$ est un voisinage de $[f]_A$.\\
Soit $U'\in\U$ borné et inclus dans $U$. Comme $f^*$ est un cône d'intérieur non vide, il existe $x\in f$ tel que $x+U'\subset f^*$. Alors $x+\vec f+U'\subset (f+U)\cap f^*$ donc $\V_A(x+\vec f,U')\subset \V'_A(f,U)$, et $[x+\vec f]_A=[f]_A$ donc $\V_A(x+\vec f,U')$ est bien un voisinage de $[f]_A$.\cqfd\\

 Étant donnés un cône $F\in\F_\I$ et un appartement $A$ contenant $b(F)$, il existe un unique élément de $\F_A$ dont le coeur est image de $\delta(F)$ par un isomorphisme d'appartements fixant $b(F)$. Nous noterons ce cône $\hat F\cap A$ car il est en fait égal à $\chap{F+\{0\}}\cap A$. Si $A$ contient un cône $f$ tel que $F=\tilde f$, c'est-à-dire si $A$ contient $\delta(F)$, alors $\hat F\cap A$ n'est autre que $f$. La propriété principale de $\hat F\cap A$ est que pour tout $U\in\U$, $\chap{F+U}=\chap{(\hat F\cap A) +U}$.\\

Nous voici en mesure d'écrire un voisinage de base de $\I$ en fonction des voisinages de base des appartements:
\begin{prop}\label{prop:voisinages dans I et dans A}
 Soit $F\in\F_\I$, $U\in\U$, borné. Alors
\[ \V_\I(F,U)=  \bigcup_{A\supset b(F)} \V'_A(\hat F\cap A,U)  \]
L'union est prise sur tous les appartements $A$ de $\I$ qui contiennent $b(F)$.
\end{prop}

\demo\\
On commence par l'inclusion "$\supset$". Soit $A\in\A$ tel que $b(F)\subset A$, notons $f=\hat F\cap A$. Soit $t\in \V'_A(f,U)$, soit $g\in\F_A$ un représentant de $t$ inclus dans $(f+U)\cap f^*$. Comme $\chap{F+U}=\chap{f+U}$, tout ce que nous devons montrer est que $\tilde g\subset \chap {f+U}$.
\begin{lemme}
Soit $A\in\A$, $U\subset \vec A$ borné et stable par $W(\vec A)$, et $f,g\in \F_A$  tels que $g\subset (f+U)\cap f^*$. Alors il existe un scp $g'$ de $g$ tel que $\tilde g'\subset \chap{f+U}$.
\end{lemme}
\pv On va appliquer la proposition \ref{lemme:partie stable par une suite de retractions}, avec $\vec \gamma=\delta(\vec g)$, $x=s(g')$ et $b=b(g')$. Il nous faut choisir $b'$ de manière à ce que $b'\subset s(g')-\vec g^*$, $\fl{\aff(b')}\supset \delta(\vec g)$. De plus, pour que la conclusion de la proposition \ref{lemme:partie stable par une suite de retractions} nous soit utile, il faut que $b(f)\subset b'$.\\

 Comme $\vec g\subset \barre{\vec f}$, il existe une facette de Weyl $\vec h$ tel que $\delta(\vec f)\cup\delta(\vec g)\subset \barre{\vec h}$. Soit alors $b'$ la base du cône $s(f)+\vec h$, de sorte que $\fl{\aff(b')}\supset \delta(\vec g)$ et $b(f)\subset b'$. Soit $g'$ un scp de $g$ tel que $b'\subset s(g')-\vec g^*$.\\
 Soit $z\in \tilde g'$, soit $Z\in\A$ contenant $b'\cup\{z\}$, et $\phi:Z\isom A$ fixant $b'$. Alors la proposition \ref{lemme:partie stable par une suite de retractions} prouve que $\phi(z)\in g'$, en particulier, $\phi(z)\in (f+U)\cap f^*$. Comme $\phi$ fixe $b'$ qui contient $b(f)$, ceci prouve que $z\in \chap{f+U}$.\cqfd\\

On passe à l'autre inclusion. Soit $t\in\V_\I(F,U)$. Il existe un représentant $G\in\F_\I$ de $t$ inclus dans $\chap{F+U}$ et un appartement $A$ contenant $b(F)\cup \delta(G)$ (corollaire \ref{cor:inclusion d'une galerie ou d'une facette et d'une chambre dans un appart}). Soit $f=\hat F\cap A$, alors $A\cap \chap{F+U}=(f+U)\cap f^*$, (par le lemme \ref{lemme:intersection chapeau appartement}) d'où $G\cap A\subset (f+U)\cap f^*$. Et comme $t=[G\cap A]_A$, on a bien $t\in\V'_A(f,U)$.\cqfd\\

\subsection{$\ib$ est séparé}

\begin{prop} L'espace topologique $\ib$ est séparé.\end{prop}

\demo\\
Soient $x$ et $y$ deux points distincts de $\ib$. D'après la proposition \ref{prop: facettes dans appart}, il existe $A\in\A$, et $f,g\in\F_A$ tels que $x=[\tilde f]$ et $y=[\tilde g]$. D'après le lemme \ref{lemme:inclusion}, $f\not\sim_A g$. Comme $\bar A$ est séparé et que $[f]_A\not=[g]_A$, on peut, quitte à restreindre $f$ et $g$ à des scp, supposer qu'il existe $U\in \U$ tel que $f+U\cap g+U=\vide$. Quitte à remplacer $U$ par une boule incluse dans $U$, on peut supposer que $U$ est convexe.\\
On peut encore restreindre $g$ à un scp pour s'assurer que $b(f)\subset s(g)-\delta(\vec g)^*$. Alors nous pouvons utiliser la proposition \ref{lemme:partie stable par une suite de retractions} pour obtenir que $\widehat{f+U}\cap\widehat{g+U}=\vide$: soit $t\in\widehat{g+U}$, soit $Z\in\A$ contenant $b(f)$ et $t$, soit $\phi:Z\isom A$ fixant $b(F)$, alors la proposition \ref{lemme:partie stable par une suite de retractions} affirme que $\phi(t)\in g+U$, donc $t\not\in \widehat{f+U}$. Ceci implique directement que $\V(\tilde f,U)\cap\V(\tilde g,U)=\vide$: nous avons obtenu un voisinage de $x$ et un de $y$ qui sont disjoints.\cqfd\\

\subsection{$\ib$ est à base dénombrable d'ouverts}

 Lorsque $\I$ est localement fini, ou plus généralement lorsque chaque cloison n'est pas incluse dans plus qu'un nombre dénombrable de chambres, la topologie qu'on vient de définir est à base dénombrable d'ouverts. Voici en quelques mots pourquoi.\\
 Tout d'abord, pour calculer $\widehat{f+U}$, seuls $U$, $s(f)$ et $b(f)$ comptent. Or, pour un sommet fixé, il n'existe qu'un nombre fini de bases le contenant. L'hypothèse sur $\I$ implique qu'il n'y a qu'un nombre dénombrable de points de $\I$ à coordonnées rationnelles (en fixant une chambre de référence, et une base affine dans cette chambre). Donc en prenant tous les sommets à coordonnées rationnelles dans $\I$, toutes les bases correspondantes à ces sommets, et tous les voisinages de $0$ dans $\vec A_0$ de la forme $B(0,\frac{1}{n})$, on obtient une base de voisinages pour la topologie de $\ib$ qui est dénombrable.

   \subsection{Injections canoniques}

\subsubsection{L'injection $\iota_\I$}

\begin{prop} L'injection canonique $\iota_\I:\I\rightarrow \ib$ est telle que, pour tout $U\in\U$, $F\in\F_\I$ on a $\iota_\I\inv(\V_\I(F,U))=\widehat{F+U}$. De plus, elle est continue, ouverte et d'image dense.\\
 \end{prop}

\demo\\

Soient $U\in U$, $F\in\F_\I$ et $x\in \iota_\I\inv(\V(F,U))$. On a alors $\{x\}\subset \chap{F+U}$, donc $x\in \chap{F+U}$.\\
 L'autre direction est (encore plus) claire.\\

Il est maintenant clair que l'image de $\iota_\I$ est dense dans $\ib$.\\

Montrons que $\iota_\I$ est ouverte. Soit $O$ un ouvert, $x\in O$. Il existe un voisinage $B$ de $x$ dans $\I$, inclus dans $O$ qui est une boule dans $\I$, $B=B(x,\epsilon)$. La boule dans $\vec A_0$ $B_{\vec A_0}(0,\epsilon)$ est stable par le groupe de Weyl de $\vec A_0$, donc $B_{\vec A_0}(0,\epsilon)\in \U$. Il est de plus clair que $B=\widehat{\{x\}+B_{\vec A_0}(0,\epsilon)}$. (Ici, $b(\{x\})=\{b\}$ et $\{x\}^*_A=A$).
Nous allons vérifier que $\iota_\I(B)=\V_\I(\{x\}, B_{\vec A_0}(0,\epsilon))$: ce sera alors bien un voisinage de $\iota_\I(x)$ dans $\ib$.\\
 L'inclusion $\subset $ découle directement de l'égalité $B=\widehat{\{x\}+B_{\vec A_0}(0,\epsilon)}$. Pour l'autre inclusion, il suffit de remarquer que puisque $B$ est borné, si $y=[\tilde f]\in \iota_\I(B)$ alors $\vec f=\{0\}$.\\

Enfin, montrons que $\iota_\I$ est continue. Soit $O$ un ouvert de $\ib$, soit $x\in\iota_\I\inv(O)$. Il existe un voisinage de $\iota_\I(x)$ inclus dans $O$ de la forme $\V(\{x\},U)$. Alors $x\in\iota_\I\inv(\V(\{x\},U))=\chap{\{x\}+U}\subset\iota_\I\inv(O)$. Il reste à prouver que $\chap{\{x\}+U}$ contient un voisinage de $x$ dans $\I$. Or $U$ contient une boule $B_{A_0}(0,\epsilon)$, donc $\chap{\{x\}+U}\supset \chap{\{x\}+B_{A_0}(0,\epsilon)}=B_\I(x,\epsilon)$ qui est un voisinage de $x$ dans $\I$.\cqfd\\

Nous pouvons donc identifier $\I$ à l'ouvert dense $\iota_\I(\I)$. La relation $\iota_\I\inv(\V(F,U))=\chap{F+U}$ devient alors $\V(F,U)\cap\I=\chap{F+U}$\\

  \subsubsection{Les injections $\iota_A$}

\begin{prop} Soit $A\in\A$, alors l'injection canonique $\iota_A:\overline{A}\rightarrow \ib$ est telle que, pour tout $f\in\F_A$, $U\in \U$ borné, il existe $f'\in\F_A$ un scp de $f$ tel que $\V_A(f',U)\subset \iota_A\inv(\V_\I(\tilde f,U)) \subset\V_A(f,U)$. Cette injection est de plus continue et fermée.\\
En particulier, c'est un homéomorphisme de $\bar A$ sur son image qui est un fermé compact de $\ib$.\end{prop}
  
\demo\\
Soit $f\in \F_A$, $U\in\U$, borné. Soit $[g]_A\in \iota_A\inv(\V_\I(\tilde f,U))$, ce qui signifie que quitte à remplacer $g$ par un scp, on a $\tilde g\subset \widehat{\tilde f+U}$. Alors $g\subset \widehat{\tilde f+U}\cap A=(f+U)\cap \delta(f)^*\subset f+U$. Donc $[g]_A\in \V_A(f,U)$.\\
Pour l'autre inclusion, notons $x=s(f)$. On choisit $f'$ un scp de $f$ tel que $f'+U\subset \delta(f)^*$, ceci existe car $U$ est borné et $\delta(f)^*$ est un cône convexe et d'intérieur non vide. Alors $f'+U\subset (f+U)\cap f^*$. Vérifions que $f'$ convient.

 Soit $[g]_A\in \V_A(f',U)$, on peut supposer que $g\subset f'+U$ d'où $g\subset (f+U)\cap f^*$. On peut aussi imposer que $b(f)\subset b(g)-\delta(\vec g)^*$. Alors en utilisant la proposition \ref{lemme:partie stable par une suite de retractions} (prendre $b=b(g)$, $b'=b(f)$, $\vec\gamma=\delta(\vec g)$), on vérifie que $\tilde g\subset \widehat{f+U}$, d'où $\iota_A([g]_A)\subset \V(\tilde f,U)$.\\

La continuité de $\iota_A$ est alors immédiate. Le fait que $\iota_A$ est fermée provient alors de la compacité de $\bar A$ et de la séparation de $\ib$.\cqfd\\

Nous pouvons dsormais identifier un appartement compactifié $\bar A$ à un fermé compact de $\ib$. La proposition \ref{prop: facettes dans appart} nous permet de prouver une généralisation de l'axiome indiquant que deux chambres doivent être incluses dans un même appartement:\\
 Soit $c$ une chambre de $\I$ et $x=[F]\in\partial \I$, alors il existe $A\in A$ contenant $c$ et un scp de $\delta(F)$. Donc $c\cup\{x\}\subset \bar A$. De même, si $x=[F]$ et $y=[G]$ sont deux points du bord de $\I$, il existe $A\in\A$ contenant un scp de $\delta(F)$ et un scp de $\delta(G)$, et alors $\{x\}\cup\{y\}\subset \bar A$.
\begin{lemme}\label{lemme:deux points ou facettes}
Soient $x$ et $y$ des facettes de $\I$ ou des points de $\partial \I$, alors il existe un appartement $A$ tel que $x\cup y\subset \bar A$. (Ici, on identifie $x$ à $\{x\}$ dans le cas où $x$ est un point.)
\end{lemme}

%subsection{Première propriétés immobilières de $\ib$}

%Il s'agit ici d'identifier quelques propriétés typiques des immeubles qui sont encore vérifiées pas $\ib$.

\subsection{Rétractions}

 Nous avons vu que si $A$ est un appartement et $w\in W_A$, alors $w$ s'étend de manière unique en un homéomorphisme de $\bar A$ sur lui-même, puis que si $\phi:A\isom B$ est un isomorphisme, alors $\phi$ s'étend de manière unique en un homéomorphisme de $\bar A$ sur $\bar B$, qu'on note encore $\phi$. Les rétractions peuvent aussi s'étendre à $\ib$:

\begin{prop} Soit $A\in\A$, et $c$ une chambre de $A$, soit $\rho=\rho_{c,A}$ la rétraction sur $A$ de centre $c$. Alors $\rho$ s'étend de manière unique en une fonction continue $\bar\rho$ de $\ib$ sur $\bar A$.
\end{prop}

\demo\\
 Soit $x\in \partial \I$, nous savons qu'il existe $Z\in\A$ tel que $\{x\}\cup c\subset \bar Z$. Soit $\phi:Z\isom A$ fixant $c$, alors $\rho_{|Z}=\phi$, il est donc naturel de poser $\rho(x)=\phi(x)$. Pour que cette définition soit valide, il faut s'assurer que si $Z'$ est un autre appartement tel que $\{x\}\cup c\subset \bar Z'$ et si $\phi':Z'\isom A$ fixe $c$, alors $\phi(x)=\phi'(x)$.\\
 Soient $f\in\F_Z$ et $f'\in\F_{Z'}$ des représentants de $x$. Alors $\rho(f)=\phi(f)\in\F_A$ et $\rho(f')=\phi'(f')\in\F_A$, et $\phi(x)=[\rho(f)]_A$, $\phi'(x)=[\rho(f')]_A$. Il nous faut donc vérifier que $\rho(f)\sim_A \rho(f')$.\\
 Quitte à raccourcir $f$ et $f'$, il existe un appartement $Y$ contenant $\delta(f)\cup\delta(f')$.
 Quitte à raccourcir encore $f$, on peut supposer, dans $Z$: $c\subset s(f)-\delta(\vec f)^*_Z$. On peut raccourcir maintenant $f'$, pour obtenir que, dans $Y$, $\delta(f')\subset \delta(f)^*$.\\

 Supposons dans un premier temps qu'aucun mur de $Y$ ne sépare strictement $\delta(f)$ de $\delta(f')$. Soit $\beta=\cl(\delta(f)\cup\delta(f'))$. Alors, en notant $\psi_Y$ un isomorphisme de $Z$ sur $Y$ fixant $\delta(f)$, si $M$ est un mur de $Z$ séparant strictement $s(f)$ et $c$, $\psi_Y(M)$ ne peut séparer strictement deux points de $\beta$. La proposition \ref{prop:inclusion d'une partie de A et d'une galerie dans un appart} permet donc de trouver un appartement $X$ contenant $\delta(f)\cup\delta(f')\cup c$. Les $\F$-cônes engendrés par $\delta(f)$ et $\delta(f')$ dans $X$ sont équivalents. Soit $\chi:X\isom A$ fixant $c$, alors $\rho(\delta(f))=\chi(\delta(f)$ et $\rho(\delta(f'))=\chi(\delta(f'))$. Il s'ensuit que les cônes engendrés par $\rho(\delta(f))$ et  $\rho(\delta(f'))$ sont équivalents dans $A$. Mais ces cônes sont exactement $\rho(f)$ et $\rho(f')$.\\

 Pour traiter le cas général, nous allons construire dans $Y$ une suite $\delta(f)=\delta_0,...,\delta_u=\delta(f')$ de coeurs qui ne sont séparés par aucun mur, et qui engendrent des cônes $f_0,...,f_k$ équivalents.
 La projection $p_{\delta(f)}:Y\rightarrow Y_{\delta(f)}$ de $Y$ sur la façade $Y_{\delta(f)}$ envoie $\delta(\vec f)$ et $\delta(\vec f')$ dans deux facettes bien déterminées. Soit $\Gamma=g_0,...,g_k$ une galerie dans $Y_{\delta(\vec f)}$ entre ces deux facettes.
 On note $f_Y$ et $f'_Y$ les cônes engendrés par $\delta(f)$ et $\delta(f')$ dans $Y$. D'après \ref{subsubsection:cx cx vectoriel de la facade}, après avoir choisi un système de générateurs $S$ de $W(\vec Y)$ adapté, on a une partition $I=I_1\sqcup I_2$, où $I$ est le type de $\delta(\vec f_Y)$, $I_1$ est l'ensemble des réflexions de $S$ fixant $\vec f^\perp$ et $I_2$ le type des réflexions fixant $\vec f_Y$. Si $\pi$ est une cloison de $Y_{\delta(f)}$ de type inclus dans $I_2$, $M$ le mur de $Y$ correspondant, alors $\vec M\supset \vec f_Y$. Par conséquent, $M$ ne peut séparer $\delta(f)$ de $\delta(f')$, car $f_Y\sim_Y f'_Y$. Ainsi le type de la galerie $\Gamma$ est inclus dans $I_1$.\\
 On définit à présent la suite de coeurs $\delta_0,\delta_1,...,\delta_{2k}$ et la suite de cônes $f_0,...,f_{2k}$ par récurrence. En plus des conditions déjà citées, on requiert que $p_{\delta(f)}(\delta_i)\subset g_{\frac{i}{2}}$ lorsque $i$ est pair.\\
On pose $\delta_0=\delta(f)$, $f_0=f_Y$. Ensuite si $\delta_0,...,\delta_i$ et $f_0,...,f_i$ sont construits pour $i$ un indice pair, soit $\sigma_i$ la réflexion, dont la direction fixe $\delta(\vec f_Y)$ qui relève la réflexion de $Y_{\delta(f)}$ définie par la cloison $g_i\wedge g_{i+1}$. La cloison $g_i\wedge g_{i+1}$ est de type un élément de $I_1$, donc $\vec\sigma_i\in W_{I_1}$, $\vec\sigma_i$ stabilise $\vec f_Y$ et fixe $\vec f_Y^\perp$. Soit $s$ le sommet de $f_i$, par hypothèse de récurrence, $f_i\sim_Y f_Y$, donc $f_i=s+\vec f_Y$, et $\sigma_i(f_i)=\sigma_i(s)+\vec\sigma_i(\vec f_Y)=\sigma_i(s)+\vec f_Y$. Soit $M_i$ le mur fixé par $\sigma_i$, donc $\vec f_Y^\perp\subset \vec M_i$, et soit $h$ le projeté orthogonal de $s$ sur $M$. Alors $\sigma_i(s)=s+2\vec{sh}$, et $\vec{sh}\in \vec M^\perp\subset \vec f_Y^{\perp \perp} = \vect(\vec f_Y)$.\\
On pose alors $f_{i+1}:=h+\vec f_Y=f_i+\vec{sh}$ et $f_{i+2}=\sigma_i(f_i)=f_i+2\vec{sh}$. Ces deux cônes sont équivalents à $f_i$. On pose ensuite $\delta_{i+1}=\delta(f_{i+1})$ et $\delta_{i+2}=\delta(f_{i+2})$. Les deux coeurs $\delta_i$ et $\delta_{i+1}$ sont projetés dans la même chambre fermée $\bar g_{\frac{i}{2}}$ de $Y_{\delta(f)}$, ce qui signifie qu'aucun mur dont la direction contient $\delta(\vec f)$ ne sépare strictement $\delta_i$ de $\delta_{i+1}$. Mais comme ces deux coeurs sont parallèles de direction $\delta(\vec f)$, un mur dont la direction ne contient pas $\delta(\vec f)$ ne peut de toute façon pas les séparer, donc aucun mur ne sépare strictement $\delta_i$ de $\delta_{i+1}$. De même, $\delta_{i+1}$ et $\delta_{i+2}$ sont projetés dans la même chambre fermée $g_{\frac{i}{2}+1}$, donc aucun mur ne sépare strictement $\delta_{i+1}$ de $\delta_{i+2}$.\\
 Ainsi $\delta_{i+1}$ et $\delta_{i+2}$ vérifient les conditions requises. On construit ainsi les suites $\delta_0,...,\delta_{2k}$ et $f_0,...,f_{2k}$. Alors $\delta_{2k}$ est projeté dans la même chambre fermée $g_k$ que $\delta(f')$ et $f_{2k}\sim_Y f_Y\sim_Y f'_Y$, il ne reste donc plus qu'à poser $f_{2k+1}=f'_Y$ et $\delta_{2k+1}=\delta(f')$ pour obtenir une suite de cônes de $Y$ équivalents allant de $f_Y$ à $f'_Y$ telle que les coeurs de deux cônes consécutifs ne sont séparés strictement par aucun mur.\\
 Nous avons que dans ces conditions, deux coeurs consécutifs ont des scp dont les images par $\rho$ sont des coeurs de cônes de $\F_A$ équivalents. Donc $\delta(f)$ et $\delta(f')$ ont des scp $\gamma$ et $\gamma'$ envoyés par $\rho$ sur des coeurs de cônes de $\F_A$ équivalents. Mais $\rho(\gamma)=\phi(\gamma)$ et $\rho(\gamma')=\phi'(\gamma)$, on a donc prouvé que $\phi(f)=\rho(f)\sim_A \phi'(f')=\rho(f')$. Ceci prouve que $\rho$ est bien définie.\\

 Montrons qu'elle est continue. Comme $\I$ est ouvert dans $\ib$, et que $\rho$ est continue, la continuité de $\bar\rho$ est acquise en tout point de $\I$. Il s'agit donc de montrer, pour $y\in \partial A$, pour $x\in \bar\rho\inv(\{y\})$ et pour $\V_A(f,U)$ un voisinage de $y$ dans $\bar A$, que $\bar\rho\inv(\V_A(f,U))$ contient un voisinage de $x$.\\
Soit $Z\in \A$ tel que $\{x\}\cup c\in \bar Z$, soit $\phi:Z\isom A$ fixant $c$. Comme $f$ n'est pas réduit à un point ($y\in \partial A$), on peut le raccourcir pour s'assurer que $c\subset s(f)-\delta(\vec f)^*$. Un voisinage de $\{x\}$ dans $\ib$ est $\V(\phi\inv(f),U)$, vérifions qu'il est inclus dans $\rho\inv(\V_A(f,U))$. Soit $B\supset b(\phi\inv(f))\cup\{t\}$ et $\psi:B\isom Z$ fixant $b(\phi\inv(f))$ tels que $\psi(t)\in (\phi\inv(f)+U)\cap \phi\inv(f)^*$, il faut montrer que $\rho(t)\in f+U$.\\
Nous savons que $c\subset s(\phi\inv(f))-\phi\inv(f)^*_Z$, donc en utilisant la proposition \ref{prop:cone et chambre dans un appart}, on voit qu'il existe $B'\in \A$ contenant $\{t\}\cup b(\phi\inv(f))\cup c$. Soit $\psi':B'\isom Z$ fixant $B'\cap Z$, en particulier, $\psi'$ fixe $b(\phi\inv(f))$ donc $\psi'(t)\in \phi\inv(f)+U$. Ensuite $\phi\circ\psi'$ fixe $c$, donc $\rho(t)=\phi\circ\psi'(t)\in f+U$.\\

L'unicité de $\bar\rho$ découle de sa continuité et de la densité de $\I$ dans $\ib$.\cqfd\\

A présent, nous pouvons nous contenter de noter $\rho$ à la place de $\bar\rho$.\\
Les résultats obtenus pour les rétractions de $\I$ peuvent souvent s'étendre aux rétractions sur $\ib$. Voici ce que devient la proposition \ref{prop:retractions adjacentes} (l'énoncé est un peu différent de celui sur les rétraction dans $\I$ car il faut prendre en compte la possibilité que $\rho_2(x)=\sigma\circ\rho_2(x)$):
\begin{prop} \label{prop:retractions adjacentes dans ib} Soit $A\in\A$, $c_1$ et $c_2$ deux chambres adjacentes dans $A$, $M$ le mur les séparant, $\sigma$ la réflexion selon $M$. Soit $\rho_1=\rho_{A,c_1}$ et $\rho_2=\rho_{A,c_2}$. Soit $x\in \ib$, alors:
\begin{itemize}
\item Soit $\rho_1(x)=\rho_2(x)$, soit $\rho_2(x)=\sigma\circ\rho_1(x)$.
\item Si $\rho_1(x)\in\overline\D(M,c_1)$ alors $\rho_1(x)=\rho_2(x)$.
\item Si $\rho_1(x)\not\in\overline\D(M,c_1)$ et si $\rho_1(x)=\rho_2(x)$, alors il existe un appartement $B$ contenant $c_1$, $c_2$ et $x$. 
\end{itemize}
\end{prop}

\demo\\
Les deux premiers points découlent de la proposition \ref{prop:retractions adjacentes}, et de la continuité de $\rho_1$, $\rho_2$ et $\sigma$.\\
Pour le troisième point, on répète quasiment la preuve de la proposition \ref{prop:retractions adjacentes}. Supposons $\rho_1(x)\not\in\overline{\D(M,c_1)}$, soit $Z_1\in\A$ tel que $\{x\}\cup c_1\subset Z_1$. Supposons que $c_2\not\subset Z_1$ et montrons qu'alors $\rho_2(x)=\sigma\circ\rho_1(x)$. Soit $\phi_1:Z_1\isom A$ fixant $c_1$, alors $\rho_1(x)=\phi_1(x)$. Soit $d$ la chambre de $Z_1$, voisine de $c_1$ le long de $\phi_1\inv(M)$ (en fait, $d=\phi\inv(c_2)$), $d\not=c_2$ puisque $c_2\not\subset Z_1$. Il existe un appartement $Z_2$ contenant $\D_{Z_1}(\phi_1\inv(M),d)\cup c_2$. Cet appartement contient $x$ dans son adhérence car il contient $\phi_1\inv(\D(M,c_2))$, et il ne contient pas $c_1$. Soit $\phi_2:Z_2\isom A$ fixant $c_2$, on a $\rho_2(x)=\phi_2(x)$. On a aussi $\phi_2(d)=c_1=\sigma\circ\phi_1(d)$. Alors $\phi_2$ et $\sigma\circ\phi_1$ coïncident sur $\D(\phi_1\inv(M),d)$, puis sur $\overline{\D(\phi_1\inv(M),d)}$ donc en $x$.\cqfd\\

\subsection{Compacité}

Pour l'instant, la notation $\widehat{f+U}$ est définie pour $f$ un cône dans un appartement et $U$ une partie de $\vec A_0$. On définit à présent, si $c$ est une facette dans un appartement $A$ et si $\vec C$ est une partie de $\vec A$ stable par les parties vectorielles des éléments de $\fix_{W_A}(c)$, 
\[ \widehat{c+\vec C}=\bigcup_{B\supset c}\phi_B(c+\vec C) \]
L'union est prise sur tous les appartements contenant $c$, et $\phi_B$ est un isomorphisme de $A$ vers $B$ fixant $c$ dont le choix n'importe pas.\\

Avant de prouver la compacité de $\ib$, vérifions cette conséquence de la proposition \ref{prop:galerie sur droite}:
\begin{lemme}\label{lemme:pour compacite}
Soit $A$ un appartement, et $C$ une partie convexe et ouverte de $A$. Soit $E$ la réunion de toute les facettes de $A$ qui rencontrent $C$. Alors pour toute chambre $c$ de $E$, et tout point $y\in\bar C$, il existe une galerie tendue de $c$ à $y$ qui reste dans $E$.
\end{lemme}

\pv Soient $x\in C\cap c$, alors $[x,y[\subset C$. Comme $c\cap C$ est ouvert, on peut supposer que $[x,y]$ ne rencontre que des chambres et des cloisons (aucune facette de dimension inférieure) et que $\vec xy$ n'est inclus dans aucune direction de mur de $A$. Soit $\Gamma$ une galerie le long de $[x,y]$, commençant par $c$, donnée par la proposition \ref{prop:galerie sur droite}. L'adhérence de chaque chambre de $\Gamma$ coupe $[x,y]$. Soit $e$ une de ces chambres, $[x,y]$ rencontre soit $e$ soit une de ses cloisons. Mais si elle rencontre une de ses cloisons, disons $m$, on ne peut avoir $\vec xy\subset \vect(m)$ car $\vect(m)$ est la direction d'un mur. Alors $[x,y]$ rencontre aussi $e$ ( $[x,y]\cap m$ est forcément dans $]x,y[$). Comme $C$ est convexe, $[x,y[\subset C$ et $e\in E$.\cqfd\\

Il est temps de montrer que $\ib$ est compact
\begin{prout} Si l'immeuble $\I$ est localement fini, alors l'espace $\ib$ est compact.\end{prout}

\demo\\
Comme nous savons déjà que $\ib$ est séparé et à base dénombrable d'ouverts, il ne reste plus qu'à montrer qu'il est séquentiellement compact. Soit $(x_n)_{n\in\N}$ une suite dans $\ib$, trouvons lui une valeur d'adhérence.\\

\textit{Première étape}\\
\textit{Quitte à remplacer $(x_n)$ par une sous-suite, il existe $A\in \A$, $f\in\F_A$, $c$ une chambre de $A$, et $\vec C$ une chambre de $\vec A$ tels que $(\rho_{A,c}(x_n))_{n\in\N}$ converge vers $y_\infini:=[f]_A$, $\delta(f)\subset \overline{c+\vec C}$ et $\forall n\in\N$, $\rho_{A,c}(x_n)\in \overline{ s(f)+\vec C}$.}\\

Commençons par choisir $c_0$ dans $A_0$ au hasard, et notons $\rho_0=\rho_{A_0,c_0}$. Comme $\overline{A_0}$ est compact, on peut supposer que la suite $(\rho_0(x_n))_n$ converge vers une limite $y_0\in\overline{A_0}$. On peut aussi supposer que tous les $\rho(x_n)$ sont dans une même façade de $A_0$, c'est à dire qu'ils ont tous la même direction $\vec g$. Soit $f_0\in y_0$ un représentant de $y_0$, alors $\vec g\subset\overline{\vec f_0}$. On peut de plus choisir $f_0$ de sorte que $c_0\subset s(f_0)-\vec f_0^*$. Comme un voisinage de $y_0$ dans $A_0$ est inclus dans $\delta(f_0)^*$, on peut supposer que tous les $\rho_0(x_n)$ ont un représentant $g_n$ inclus dans $\delta(f_0)^*$. Mais $\delta(\vec f_0)^*$ est une réunion finie de chambres fermées, donc il existe $\vec C_0$ chambre de $\vec A_0$ telle que $s(f_0)+\overline{\vec C_0}$ contient une infinité de $g_n$, on peut supposer qu'il les contient tous. Bien sûr, cela implique que $\overline {s(f_0)+\vec C_0}$ contient tous les $\rho_0(x_n)$.\\
Pour conclure la première étape, il ne reste plus qu'à s'assurer que $\delta(f_0)\subset \overline{c_0+\vec C_0}$, ce qui, puisque $\delta(\vec f_0)\subset \overline{\vec C_0}$ revient à s'assurer que $\bar c_0$ coupe $\overline{s(f_0)-\vec C_0}$. Pour l'instant, $c_0$ est inclus seulement dans $s(f_0)-\delta(\vec f_0)^*$. Soit $d$ une chambre à distance minimale de $c_0$ dont l'adhérence coupe $\overline{s(f_0)-\vec C_0}$.\\
 Nous allons maintenant prouver par récurrence sur $k$ que s'il existe $A'$, $c'$, $\vec C'$, $f'$ et $d'$ dans la même configuration que $A_0$, $c_0$, $\vec C_0$, $f_0$ et $d$ le sont actuellement, avec $d(c',d')=k$, alors il existe $A$, $c$, $f$, et $\vec C$ comme requis pour conclure la première étape.\\

Si $k=0$ la conclusion est immédiate, $A'$, $c'$, $\vec C'$ et $f'$ conviennent.\\

Si $k>0$, soit $\Gamma=(c'=c_0,c_1,...,c_k=d')$ une galerie minimale entre $c'$ et $d'$, et soit $M=\vect(c_0\wedge c_1)$ le premier mur traversé, soit $\sigma_M$ la symétrie de $A'$ selon $M$, notons $\rho_0=\rho_{A',c'}$ et $\rho_1=\rho_{A',c_1}$.\\
Si $\rho_0(x_n)=\rho_1(x_n)$ pour une infinité d'indices $n$, alors on peut passer à une sous-suite pour supposer que $\rho_0(x_n)=\rho_1(x_n)$ pour tout $n$. Alors il suffit de remplacer $c'$ par $c_1$ pour pouvoir appliquer l'hypothèse de récurrence.\\

Par contre, dans l'autre cas, on peut supposer que $\rho_1(x_n)=\sigma\circ\rho_0(x_n)$ pour tout $n\in\N$. Alors la suite des $\rho_1(x_n)$ ne vérifie plus du tout les hypothèses demandées, il va falloir changer d'appartement.

 D'après la proposition \ref{prop:retractions adjacentes dans ib}, cette situation n'est possible que si pour tout $n$, $\rho_0(x_n)\not\in\overline{\D(M,c')}$. De plus, pour tout $n\in\N$, il existe alors une chambre $e_n$ voisine de $c'$ et de $c_1$, et un appartement dont l'adhérence contient $e_n$, $c'$ et $x_n$. Comme $\I$ est supposé localement fini, il n'y a qu'un nombre fini de chambres voisines à $c'$ et à $c_1$, on peut donc remplacer $(x_n)_n$ par une sous-suite pour qu'il existe une chambre fixe $e$, voisine de $c_1$ et $c'$ telle que pour tout $n\in\N$, il existe $B_n\in\A$ tel que $c'\cup e\cup\{x_n\}\in \bar B_n$.\\
Soit $A''$ un appartement contenant $c'$ et $e$, et soit $\phi:\bar A'\isom \bar A''$ fixant $A'\cap A''$. Nécessairement, $\phi(c_1)=e$. Alors $\forall n\in\N$, $\rho_{A'',e}(x_n)=\phi(\rho_0(x_n))$. Comme de plus $d(e,\phi(d))=k-1$, on peut appliquer l'hypothèse de récurrence à $A''$, $e=\phi(c_1)$, $\phi(d)$, $\phi(f)$, et $\vec\phi(\vec C)$. Ceci conclut la première étape.\\

\textit{Seconde étape:}\\
\textit{On note $\rho=\rho_{A,c}$ la rétraction obtenue dans la première étape. On montre à présent que soit la suite $(x_n)_n$ admet une valeur d'adhérence dans $\I$, soit il existe une galerie tendue $(d_n)_{n\in\N}$ dans $A$ issue de $d_0=c$ et une suite $(t_n)_n$ telle que pour tout $n$, pour tout $i\geq n$, $\rho(x_i)\in \overline{d_n+\vec C}$, $t_n\in d_n$ et $\lim_{n\rightarrow \infini}t_n=y_\infini$. On peut enfin imposer que $\forall n\in\N$, $d_n\cap (c+\vec C)\not= \vide $.}\\

On note $\rho=\rho_{A,c}$.\\
On commence par regarder le cas où $\vec f=\{0\}$. Dans ce cas, $y_\infini\in A$, et $K=\bigcup_{\bar c\supset y_\infini}\bar c$ est un voisinage de $y_\infini$ dans $A$, il contient donc tous les $\rho(x_n)$ à partir d'un certain rang. Mais $\rho\inv(K)$ est une union finie de chambres fermées, donc c'est un compact, et la suite $(x_n)_n$ a une valeur d'adhérence dans $\rho\inv(K)$ donc dans $\I$.\\

Supposons à présent que $\vec f\not=\{0\}$.
Par le lemme \ref{lemme:pour compacite}, il existe une galerie minimale $c=d_0,d_1,...,d_u$ de $c$ à $s(f)$ telle que $\forall i\in \llbracket 0,u\rrbracket $, $d_i$ coupe $c+\vec C$.\\
 On choisit ensuite une demi-droite $\vec l\subset \delta(\vec f)$, et on prolonge $d_0,...,d_u$ par une galerie tendue dont l'adhérence contient $l:=s(f)+\vec l$ (proposition \ref{prop:galerie sur droite}). Vérifions qu'on peut s'assurer que $d_n\cap (c+\vec C)\not=\vide$ pour $n>u$. Notons $\M_l$ l'ensemble des murs contenant $l$, alors la construction faite dans la proposition \ref{prop:galerie sur droite} est telle que chaque chambre de la galerie $(d_n)_{n\geq d_u}$ est incluse dans $\bigcap_{M\in\M_l}\D(M,d_u)$. Par ailleurs, $c+\vec C$ est un convexe, délimité par un nombre fini de ses hyperplans d'appui, et ces hyperplans ont pour direction un mur de $\vec A$. Notons $\mathcal H_l$ l'ensemble des hyperplans d'appui de $c+\vec C$ qui contiennent $l$.\\
 Il suffit de montrer que $d_i \cap (c+\vec C)\not=\vide$ pour les chambres $d_i$ contenant un segment de $l$ dans leur adhérence, car alors nous pourrons relier deux telles chambres par une galerie dont chaque chambre coupe $c+\vec C$, et la construction donnée dans la preuve de la proposition  \ref{prop:galerie sur droite} consiste exactement à définir d'abord les chambres contenant un segment de $l$ dans leur adhérence, puis à relier ces chambres par des galeries tendues arbitraires. Soit $d_i$ une de ces chambres, soit $]x-\vec\epsilon, x+\vec\epsilon[$ un segment de $l$ inclus dans $\bar d_i$. Les hyperplans d'appui de $d_i$ sont des murs en nombre fini, et ceux contenant $x$ doivent contenir $]x-\vec\epsilon, x+\vec\epsilon[$ et sont donc des éléments de $\M_l$. Par conséquent, il existe un voisinage $V$ de $x$ tel que  $V\cap \bigcap_{M\in \M_l}\D(M,d_u)\subset \bar d_i$. De la même manière, les hyperplans d'appui de $c+\vec C$ qui contiennent $x$ doivent contenir $l$ car $l\subset \overline{c+\vec C}$. Donc il existe un voisinage $U$ de $x$ tel que $U\cap \bigcap_{H\in \mathcal H_l}\D(H,c)\subset c+\vec C$. Il ne reste plus qu'à vérifier que $U\cap V\cap\bigcap_{H\in\mathcal H_l}\D(H,c)\bigcap_{M\in\M_l}\D(M,d_u)$ est d'intérieur non vide. Notons $S=\bigcap_{H\in\mathcal H_l}\D(H,c)\bigcap_{M\in\M_l}\D(M,d_u)$, nous savons que $S$ est d'intérieur non vide puisque $d_u$ contient un point disons $u$ de $c+\vec C\cap \bigcap_{M\in\M_l}\D(M,d_u)$. De plus, $l\subset \bar S$ et $S$ est convexe. Il ne reste plus qu'à relier $u$ à $x$ par une droite pour finalement obtenir un point intérieur à $U\cap V\cap S$.\\

Ceci permet d'obtenir que $d_n\cap (c+\vec C)\not=\vide$ pour tout $n$. Vérifions que la galerie $(d_n)_{n\in\N}$ est tendue. Si ce n'est pas le cas, c'est qu'il y a un mur $M$ qui sépare $d_u$ de $d_0$ et d'une chambre $d_n$ avec $n\geq u$. Mais $s(f)\in \bar{d_u}\cap (\overline{c+\vec C})$ et $\bar d_n\cap \overline{d_u+\vec l}\not=\vide$, d'où $l\subset M$. En particulier $s(f)\in M$, mais ceci est impossible puisque $d_0,...,d_u$ est une galerie minimale de $d_0$ à $s(f)$.\\
 Donc la galerie $(d_n)_{n\in\N}$ obtenue est bien tendue.\\

Pour tout $i\in\N$, soit $s_i\in l\cap\bar d_i$. La suite $(s_i)$ tend vers $y_\infini$. Donc pour tout $i$, les $\rho(x_n)$ sont inclus à partir d'un certain rang $N_i$ dans $\overline{s_i+\delta(\vec f)^*}$. Vérifions qu'en fait ils sont dans $\overline{s_i+\vec C}$. Soit $n\geq N_i$, un représentant $g_n$ de $\rho(x_n)$ est inclus dans $s_i+\vec F$ où $\vec F$ est une facette vectorielle dont l'adhérence contient $\delta(\vec f)$ et $\vec g$. Disons $g_n=s_i+\vec v+\vec g$ avec $\vec v\in \vec F$.

 On calcule alors: $g_n=s(f)+\fl{s(f)s_i} +\vec v+\vec g$. Or $\fl{s(f)s_i}\in \delta(\vec f)$ donc $\fl{s(f)s_i} +\vec v\in\vec F$ et $g_n\subset s(f)+\vec F$. Mais nous savons d'autre part que $\rho(x_n)\in \overline{s(f)+\vec C}$. Donc $\vec F\subset \overline{\vec C}$, et $\rho(x_n)\in\overline{s_i+\vec F}\subset \overline{s_i+\vec C}$.\\
 Ainsi, pour tout $i\in \N$, les $\rho(x_n)$ sont inclus dans $\overline{s_i+\vec C}$ à partir d'un certain rang. En passant à une sous-suite, on peut supposer que $\forall n\geq i$, $\rho(x_n)\in d_i+\vec C$.\\
Par ailleurs, en déplaçant un peu les $s_i$ on peut obtenir une suite $(t_i)$ convergeant vers $y_\infini$ et telle que $t_i\in d_i$. Ceci conclut la seconde étape.\\

\textit{Troisième étape}\\
\textit{On suppose maintenant que $(x_n)$ n'a pas de valeur d'adhérence dans $\I$. Montrons qu'il existe une sous-suite de $(x_n)$, et une galerie tendue $\Gamma=(c_i)_{i\in\N}$ telle que $c_0=d_0=c$, que pour tout $i$, $\rho(c_i)=d_i$ et que les $x_k$ pour $k\geq i$ sont dans l'adhérence de $\widehat{c_i+\vec C}$.}\\
Ici, on note $\widehat{c_i+\vec C}=\widehat{c_i+\vec\phi_i(\vec C)}$ pour un isomorphisme $\phi_i$ de $A$ sur un appartement $B_i$ contenant $c$ et $c_i$ qui fixe $c$. On remarque que puisque $c_i$ est une chambre, on a en fait $\widehat{c_i+\vec C}=\rho_{B_i,c_i}\inv (c_i+\vec\phi_i(\vec C))$, et l'adhérence de ceci vaut $\rho_{B_i,c_i}\inv (\overline{c_i+\vec\phi_i(\vec C)})$ (pour vérifier l'inclusion non évidente, utiliser le lemme \ref{lemme:deux points ou facettes} pour inclure un point et $c_i$ dans un même appartement).\\

On commence par montrer par récurrence que pour tout $i$, il existe $c_0,...,c_i$ et une sous-suite $(x^i_n)$ de $(x_n)$ satisfaisant aux conditions demandées.\\
 D'après la première étape, tous les $x_n$ sont dans $\rho\inv(\overline{c+\vec C})$. Donc $c_0=c$ convient bien.\\
Ensuite, supposons $c_0,...,c_i$ construits.
Comme $x_n\in \overline{\widehat{c_i+\vec C}}$ et $c_i$ coupe $\widehat{c_0+\vec C}$, en utilisant la proposition \ref{prop:cone et chambre dans un appart}, on trouve un appartement $B_n$ tel que $c_0\cup c_i\cup\{x_n\}\subset \bar B_n$. Soit $\phi_n:B_n\isom A$ fixant $c_0$, comme $\phi_n(x_n)=\rho(x_n)\in \overline{d_{i+1}+\vec C}$, pour $n\geq i+1$, on voit que $x_n\in \overline{\widehat{\phi_n\inv(d_{i+1})+\vec C}}$. Or $\phi_n\inv(d_{i+1})$ est une chambre adjacente à $c_i$, et il n'y a qu'un nombre fini de telles chambres. On peut donc passer à une sous suite $(x^{i+1}_n)$ de $(x^i_n)$ pour que tous les $x^{i+1}_n$ soient dans un même $\overline{\widehat{e+\vec C}}$ avec $e$ une chambre adjacente à $c_i$ dont l'image par $\rho$ est $d_{i+1}$. Il ne reste plus qu'à poser $c_{i+1}=e$.\\

Enfin, on prend $(y_n)$ la sous-suite de $(x_n)$ définie par $y_n=x^n_n$. Le couple $(\Gamma=(c_0,c_1,...), (y_n))$ satisfait aux conditions demandées par l'étape 3. Pour ne pas alourdir les notations, dans la suite on notera $x_n=y_n$.\\

\textit{quatrième étape}\\

 Soit $Z$ un appartement contenant $\Gamma$, soit $\phi:Z\isom A$ fixant $c_0$, soit $x_\infini=\phi\inv(y_\infini)$. Montrons que $x_n$ tend vers $x_\infini$.\\

Soit $\V(g,U)$ un voisinage de $x_\infini$. Comme $\V_A(\phi(g),U)\cap \phi(g)^*$ est un voisinage de $y_\infini$, il existe un entier $N$ tel que $d_N\cap (\phi(g)+U)\cap \phi(g)^*\not=\vide$ et tel que $\forall n\geq N$, $\rho(x_n)\subset \V_A(\phi(g),U)\cap \phi(g)^*$. Soit $n\geq N$, montrons que $x_n\in \V(g,U)$.\\
Comme $c_N$ coupe $\widehat{c_0+\vec C}$ et $x_n\in \overline{\widehat{c_N+\vec C}}$, la proposition \ref{prop:cone et chambre dans un appart} indique l'existence d'un appartement $B$ tel que $\{x_n\}\cup c_N\cup c_0\subset \bar B$. Dès lors, on a les égalités: $\rho_{Z,c_N}(x_n)=\rho_{Z,c_0}(x_n)= \phi\inv\circ\rho(x_n)$, d'où $\rho_{Z,c_N}(x_n)\in \V_Z(g,U)\cap g^*$. Mais comme $c_N$ coupe $b(g)+\vec g^*$, on peut appliquer la proposition \ref{lemme:partie stable par une suite de retractions}, qui prouve que si un appartement fermé $\bar X$ contient $b(g)$ et $x_n$, et si $\psi:\bar X\isom \bar Z$ fixe $b(g)$, alors $\psi(x_n)\in V_Z(g,U)\cap g^*$. Ceci entraîne bien que $x_n\in \V(g,U)$.\cqfd\\

\section{Unicité de la construction}
\label{section:unicite}
On suppose dans cette partie que $\I$, muni de son système complet d'appartements, est l'immeuble obtenu à partir d'une donnée radicielle valuée discrète comme dans \cite{bruhat-tits}. Rappelons brièvement quelques notations.\\
% $(U_{\alpha,\lambda})_{\alpha\in\phi,\lambda\in\Lambda_\alpha}$ d'un groupe $G$
On choisit un appartement $A_0$, ainsi qu'un sommet spécial $0\in A_0$, qui permet d'identifier $A_0$ et $\fl{A_0}$ ainsi que $\vec A_0^*$ et les formes affines sur $A_0$ s'annulant en $0$. Il existe un système de racine $\phi$, éventuellement non réduit, sur $\vec A_0$, et pour tout $\alpha\in\phi$, il existe $\Gamma_\alpha$ un sous-groupe discret de $(\R,+)$ tel que les murs de $A_0$ sont les $\M(\alpha,k):=\ens{x\in A_0}{\alpha(x)+k=0}$ pour $\alpha\in\phi$ et $k\in\Gamma_\alpha$.\\
 Un groupe $G$ agit sur $\I$, pour tout $\alpha\in\phi$ et $k\in\Gamma_\alpha$, il existe un sous-groupe $U_{\alpha,k}$ de $G$ qui fixe le demi-appartement $\D(\alpha,k):=\ens{x\in A_0}{\alpha(x)+k\geq 0}$ ainsi que toutes les chambres ayant une cloison dans son intérieur. Pour $\alpha\in\phi$ et $k\leq l\in \Gamma_\alpha$, on a $U_{\alpha,l}\subset U_{\alpha,k}$, donc la réunion $U_\alpha:=\bigcup_{k\in\Gamma_\alpha}U_{\alpha,k}$ est un groupe. Le groupe $U_{\alpha,k}$ agit transitivement sur l'ensemble des appartements contenant $\D(\alpha,k)$. Le groupe $G$ est engendré par les $U_{\alpha,k}$, pour $\alpha\in\phi$ et $k\in \Gamma_k$ ainsi que par un autre sous-groupe noté $H$ qui fixe $A_0$.\\
 Ces sous-groupes vérifient les axiomes d'une donnée radicielle valuée. On suppose de plus l'axiome sur les commutateurs des $U_{\alpha,k}$ vérifié au sens fort: pour tout $\alpha,\beta\in\phi$ avec $\beta\not\in \R^-.\alpha$, pour tout $k\in\Gamma_\alpha$, $l\in\Gamma_\beta$,
\[[U_{\alpha,k},U_{\beta,l}]\subset \eng{ U_{p\alpha+q\beta,pk+ql}}{ p,q\in\N^*\et p\alpha+q\beta\in\phi }\]

 Le groupe $G$ est muni de la topologie de la convergence simple sur $\I$, c'est-à-dire qu'une base de voisinages de $e$ est l'ensemble des fixateurs de parties bornées. Alors l'action de $G$ sur $\I$ est continue.\\

 Cette situation se retrouve en particulier lorsque $\I$ est l'immeuble de Bruhat-Tits d'un groupe réductif sur un corps local non archimédien.\\
 Nous allons montrer dans ce cas qu'une fois fixée une compactification de $A_0$ de la manière présentée dans la partie \ref{section:compactification de lappart}, il n'existe en fait qu'une seule manière de l'étendre en une compactification de $\I$ sur laquelle $G$ agit par homéomorphisme.\\

Le fixateur dans $A_0$ d'un produit d'éléments de divers sous-groupes radiciels $U_{\alpha,k}$ est dans certains cas facile à déterminer:
\begin{sprop}\label{prop: intersection des fixateurs}
Soient $\alpha_1,...,\alpha_k \in\phi$ des racines deux à deux non colinéaires, soit $(u_1,...,u_k)\in U_{\alpha_1}\times...\times U_{\alpha_k}$. On suppose que $\bigcap_i \fix_{A_0}(u_i)$ contient au moins une chambre. Alors $\fix_{A_0}(u_1...u_k)=\bigcap_{i=1}^k \fix_{A_0} (u_i)$.\\

\end{sprop}
\rema Pour des racines dans un sous-système réduit et positif de $\phi$, l'hypothèse "non colinéaires" devient juste "distinctes".\\

\demo\\
On note $g=u_1...u_k$. Pour chaque $i$, l'ensemble des points fixes de $u_i$ est un demi-appartement $\D_i$ délimité par un mur $M_i$ de direction $\ker(\alpha_i)$.\\
La première inclusion est claire.\\
 Soit $x\in A_0\setminus \bigcap_{i=1}^k \D_i$, montrons que $x$ n'est pas fixé par $g$.
On remarque pour commencer qu'il est impossible que $x$ soit fixé par $g$ et tous les $u_i$ sauf un.\\
Soit $c$ une chambre incluse dans $\bigcap_i \D_i$. Comme $c$ est ouvert dans $A_0$, on peut choisir $y\in c$ tel que $[x,y]$ ne rencontre pas d'intersection de murs.  Si $x$ était fixé par $g$, alors $[x,y]$ le serait aussi. Or le segment $[x,y]$ sort de $\bigcap_i \D_i$ en coupant un seul mur, disons $M_j$. Donc $[x,y]$ contient un point $z$ qui est dans $(\bigcap_{i\not= j} \D_i)\setminus \D_j$. Donc $z$ est fixe par tous les $u_i$ sauf $u_j$, ce qui est impossible.\cqfd\\

Rappelons ce résultat (voir \cite{bruhat-tits} 7.4.33), exprimant que lorsqu'un élément de $U_{\alpha,k}$ fixe une grande boule dans un appartement, alors il fixe une petite boule de même centre dans l'immeuble. On note $B_A(o,r)$ la boule dans $A$ de centre $o$ de rayon $r$, et $B_\I(o,r)$ la boule dans $\I$.\\

\begin{sprop} \label{prop:fixateur dans moufang} Il existe une constante $\lambda$ dépendant uniquement de $\phi$ telle que pour tout  $o\in A_0$, et $r\in\R^+$, si $\alpha\in\phi$ et $k\in\Lambda_\alpha$ sont tels que $B_{A_0}(o,\lambda r)\subset\D(\alpha,k)$ et si $u\in U_{\alpha,k}$, alors $u$ fixe $B_\I(o,r)$.\end{sprop}

\textit{preuve:}\\
Pour $\beta\in\phi$, notons $M_\beta$ le maximum, de $|\vec\beta |$ sur la sphère unité, et soit $\lambda=\max\ens{\frac{M_{\alpha}}{M_\beta}} {\alpha,\beta\in\phi}$. 
\\

Soient $\alpha\in \phi$ et $k\in \Gamma_\alpha$ tels que $B_{A_0}(o,\lambda r)\subset \D(\alpha,k)$.
% Alors pour $\vec v$ le vecteur qui maximise $\vec \alpha$ sur la sphère unité, on a $o-\lambda r\vec v\in \D(\alpha,k)$ d'où $(\alpha+k)(o)\geq \lambda r M_\alpha$. Alors pour tout $\beta\in\phi$ et $l\in \Gamma_\beta$ tels que $o\in\D(\beta,l)$, pour tout $p,q\in\N^*$ tels que $p\alpha+q\beta$ est une racine, on a $(p\alpha+q\beta+pk+ql)(o)\geq \lambda r M_\alpha$, et ceci entraine que la boule $B_{A_0}(o,\lambda r \frac{M_\alpha}{M_{p\alpha+q\beta}})$ est incluse dans $\D(p\alpha+q\beta,pk+ql)$. Comme $\lambda  \frac{M_\alpha}{M_{p\alpha+q\beta}}\geq 1$, on voit donc que pour tous tels $\beta, l,p$ et $q$, le demi-appartement $\D(p\alpha+q\beta,pk+ql)$ contient la boule $B_{A_0}(o,r)$, en particulier, pour tout $v\in U_{\beta,l}$, $v$ fixe $B_{A_0}(o,r)$.\\

 Soit $V$ le groupe engendré par les $U_{p\alpha+q_1\beta_1+...+q_n\beta_n, pk+q_1 l_1+...+q_n l_n}$ pour $n\in\N$, $\beta_1,...,\beta_n\in\phi$, $(l_1,...,l_n)\in \Gamma_{\beta_1}\times...\times\Gamma_{\beta_n}$ et $p,q_1,...,q_n\in \N^*$ tels que $p\alpha+q_1\beta_1+...+q_n\beta_n\in\phi$ et chaque $\D(\beta_i,l_i)$ contient $o$. Montrons que les éléments de $V$ fixent $B_{A_0}(o,r)$. Déjà, pour $\vec v$ le vecteur qui maximise $\vec \alpha$ sur la sphère unité, on a $o-\lambda r\vec v\in \D(\alpha,k)$ d'où $(\alpha+k)(o)\geq \lambda r M_\alpha$. Soit maintenant $\gamma=p\alpha+q_1\beta_1+...+q_n\beta_n$ et $j=pk+q_1 l_1+...+q_n l_n$ comme dans la définition de $V$. Alors $\gamma(o)+j\geq p(\alpha(o)+k)\geq \lambda r M_\alpha$. Ceci entraîne que la boule $B_{A_0}(o,\lambda  \frac{M_\alpha}{M_{\gamma}} r)$ est incluse dans $\D(\gamma,j)$. Comme $\lambda \frac{M_\alpha}{M_\gamma}\geq 1$, on obtient le résultat voulu.\\

A présent, soit $u\in U_{\alpha,k}$, et $d$ une chambre de $\I$ coupant $B_\I(o,r)$. Soit $\Gamma=c_0,...,c_n=d$ une galerie tendue avec $o\in\bar c_0\subset A_0$. Soit $\rho$ la rétraction sur $A_0$ de centre $c_0$. Il existe $\beta_1\in\phi$, $l_1\in \Gamma_{\beta_1}$ et $v_1\in U_{\beta_1,l_1}$ tels que $v_1.c_1=\rho(c_1)$ et $\D(\beta_1,l_1)$ contient $c_0$ donc $o$. Ensuite il existe $\beta_2\in\phi$, $l_2\in\Gamma_{\beta_2}$, $v_2\in U_{\beta_2,l_2}$ tels que $v_2v_1c_2=\rho(c_2)$, et $\D(\beta_2,l_2)$ contient $c_0\cup \rho(c_1)$ donc encore $o$. On obtient finalement $\beta_1,...,\beta_n$, $l_1,...,l_n$, $v_1,...,v_n$ tels que $\rho(d)=v_n...v_2v_1.d$. Notons $v=v_n...v_1$.\\
Comme $\rho(d)$ coupe $B_{A_0}(o,r)\subset B_{A_0}(o,\lambda r)$, on a $u.\rho(d)=\rho(d)=v.d=uv.d$. Le fait que chaque commutateur $(u,v_i)=uv_iu\inv v_i\inv$ soit dans $V$ permet de prouver par un calcul classique que $(u,v)\in V$, et donc que $(u,v)\inv$ fixe $v.d$. Alors $v.d=(u,v)vu.d$ d'où $(u,v)\inv v.d=vu.d$ puis $d=u.d$.\cqfd\\

%\begin{cor}
%Soit $c$ une chambre de $A_0$, soit $o\in c$, soit $r\in\R^{+*}$. Alors tout élément de $U(c)$ qui fixe $\B_{A_0}(o,\lambda r)$ fixe aussi $\B_\I(o,r)$, où $\lambda$ est la constante donnée par la proposition.
%\end{cor}
%\demo Soit $g\in U(c)$ qui fixe $\B_{A_0}(o,\lambda r)$. On décompose $g$ en $u_{\alpha_1}...u_{\alpha_k}$ comme dans la proposition \ref{prop: decomposition unipotent}. Alors par la proposition \ref{prop: intersection des fixateurs}, $\B_{A_0}(o,\lambda r)\subset \fix_{A_0}(g)=\bigcap_i \fix_{A_0}(u_i)$. Donc chaque $\fix_{A_0}(u_i)$ contient $\B_{A_0}(o,\lambda r)$, et donc par la proposition \ref{prop:fixateur dans moufang}, chaque $u_i$ fixe $\B_\I(o,r)$.\cqfd\\

\begin{scor}\label{cor:fixe boule}
Soient $A$ et $B$ deux appartements, on suppose qu'il existe $o\in A\cap B$ et $r>0$ tels que $A\cap B$ contient la boule $B_A(o,\lambda r)$ où $\lambda$ est la constante de la proposition. Alors il existe $g\in G$ tel que $g.A=B$, et $g$ fixe $(A\cap B)\cup B_\I(o,r)$.
\end{scor}
\demo\\
Par conjugaison, on peut supposer $A=A_0$. Le groupe engendré par les $U_{\alpha,k}$ tels que $\D(\alpha,k)\supset B_A(o,\lambda r)$ est transitif sur les appartements contenant $B_A(o,\lambda r)$ et il fixe $B_\I(o,r)$ d'après \ref{prop:fixateur dans moufang}.\cqfd\\

% Comme $r>0$, il existe une chambre $c$ qui rencontre $B_{A_0}(o,\lambda r)$. Comme $U(c)$ agit transitivement sur les appartements contenant $c$, il existe $g\in U(c)$ tel que $g.B=A_0$. Un tel $g$ fixe $c$ et donc $B\cap A_0$, donc en particulier $\B_{A_0}(o,\lambda r)$ et donc $\B_\I(o,r)$ par le corollaire précédent.\cqfd\\

Et voici la proposition annoncée:\\
\begin{sprop}
Soit $\hat \I$ un espace topologique compact contenant $\I$, tel que $\I$ est ouvert dense dans $\hat \I$.\\
On suppose que l'adhérence de l'appartement de référence $A_0$, qu'on note $\hat A_0$ est homéomorphe à $\bar A_0$ par un homéomorphisme fixant $A_0$. On note $\wedge:a\mapsto \hat a$ cet homéomorphisme.\\
 On suppose que l'action de $G$ sur $\I$ s'étend en une action continue sur $\hat\I$. Alors il existe un homéomorphisme $\chi:\ib \rightarrow \hat\I$qui est $G$-équivariant et qui fixe $\I$.\end{sprop}

\rema Par densité de $A_0$ dans $\bar A_0$, et par densité de $\I$ dans $\ib$, les applications $\wedge: a\mapsto \hat a$ et $\chi$ sont uniques si elles existent.\\

\demo\\
Si $A$ est un appartement de $\I$, on notera $\hat A$ son adhérence dans $\hat\I$. Les deux clôtures d'appartements $\bar A$ et $\hat A$ sont alors homéomorphes, par un unique homéomorphisme qui vaut $g_{|\hat A}\inv\circ\wedge\circ g_{|\bar A}$, où $g\in G$ est tel que $g.A=A_0$.\\
 Il est donc naturel de vouloir définir une fonction de la sorte:
\fonc{\chi}{\ib}{\hat \I} {g.a}{g.\hat a}, pour $g\in G$ et $a\in \bar A_0$.\\

 Montrons que $\chi$ est bien définie: soient $g,g'\in G$ et $a,a'\in \bar A_0$ tels que $ga=g'a'$ et montrons que $g\hat a=g'\hat a'$. En remplaçant $g$ par $(g')\inv g$, on se ramène au cas où $g'=1$. Soient $f$ un représentant de $a$ et $f'$ un représentant de $a'$. Quitte à les réduire, il existe $Z\in\A$ contenant $g(\delta(f))\cup\delta(f')$. Soient $d$ et $d'$ des demi-droites incluses dans $\delta(f)$ et $\delta(f')$ de directions incluses dans $\vec f$ et $\vec f'$, de sorte que $a$ est l'unique point dans $\bar\I$ de $\bar d\setminus d$ et $a'$ est l'unique point dans $\bar \I$ de $\bar d'\setminus d'$. Ainsi, $d'$ et $g.d$ ont le même point limite dans $\bar Z$, qui est $a'$. Comme $\bar Z$ est homéomorphe à $\hat Z$ par un homéomorphisme fixant $Z$, on en déduit que $d'$ et $g.d$ ont le même unique point limite dans $\hat Z$.\\
Le point limite de $d'$ dans $\hat A_0$ est $\hat a'$ et comme $\hat \I$ est séparé, c'est aussi l'unique point limite de $d'$ dans $\hat I$ donc en particulier dans $\hat Z$. De même, l'unique point limite de $d$ dans $\hat \I$ est $a$, donc par continuité de $g$, le point limite de $g.d$ dans $\hat \I$, donc dans $\hat Z$ est $g.\hat a$. Ceci prouve que $\hat a'=g.\hat a$.\\

Donc $\chi$ est bien définie. C'est une fonction $G$-équivariante qui fixe $\I$, elle induit l'unique homéo\-morphisme de $\bar A$ sur $\hat A$ pour chaque appartement $A$. Comme deux points quelconques de $\ib$ sont toujours inclus dans l'adhérence d'un même appartement (lemme \ref{lemme:deux points ou facettes}), ceci entraîne que $\chi$ est injective. Il ne reste plus qu'à montrer que $\chi$ est continue: elle sera alors fermée par compacité de $\ib$ et séparation de $\hat\I$, et son image sera un fermé de $\hat\I$ contenant $\I$ c'est-à-dire $\hat\I$ lui-même.\\

 Soit $x\in \ib$, montrons que $\chi$ est continue en $x$. Le cas où $x\in \I$ est évident, et par $G$-équivariance de $\chi$, on se ramène au cas où $x\in\bar A_0$. Au total, on peut supposer $x\in\partial A_0$.

 Soit $O$ un ouvert de $\hat \I$ contenant $\chi(x)$.
Notons $\hat\rho: G\times \hat\I\rightarrow \hat\I$ l'action de $G$ sur $\hat\I$. Par continuité de $\hat\rho$, et parce que $e.x=x\in \bar A_0$, il existe un voisinage ouvert $Q$ de $e$ dans $G$ et un voisinage ouvert $\hat V$ de $\hat x$ dans $\hat \I$ tel que $Q.\hat V=\rho(Q\times \hat V)\subset O$. Soit $V=\chi\inv(\hat V)$. Alors $x\in Q.V\subset\chi\inv(Q.\hat V)\subset \chi\inv(O)$. Il reste à prouver que $Q.V$ contient un voisinage de $x$ dans $\ib$.

 Soit $f\in\F_{A_0}$ un représentant de $x$. On peut supposer que $Q$ est le fixateur dans $G$ d'une partie bornée,  de la forme $B_\I(o,r)$, avec $r\in \R^{+*}$ et $o\in A_0$, on peut même imposer $o\in -\delta(f)$, où $-\delta(f)$ signifie $s(f)-\delta(\vec f)$.

% On choisit ensuite un scp $\delta_1$ de $\delta(f)$ tel que $B_A(o,\frac{r}{\lambda})\subset s(\delta_1)-\vec delta_1^*$, ou $\lambda$ est la constante donnée par la proposition \ref{prop:fixateur dans moufang}.\\

 Soit $\lambda$ la constante donne par la proposition \ref{prop:fixateur dans moufang}. Soient $\{\beta_i\}_{i\in I}$ l'ensemble des parties d'appartement égales à la clôture de $\delta(f)$ et d'une boule de centre $o$ de rayon $\lambda r$, c'est-à-dire:
\[ \{\beta_i\}_{i\in I}=\ens{\beta \text{ sous-complexe simplicial de } \I}{\exists B\in \A \tq \delta(f)\subset B,\; o\in B \et \beta=\cl(B_B(o,\lambda r)\cup\delta(f)) } \]
 Cet ensemble est fini car l'ensemble des sous-complexes simpliciaux de $B_\I(o,r)$ est fini.\\
On choisit alors $\{B_i \}_{i \in I}$ une famille d'appartements tels que $\forall i \in I$, $\beta_i\subset B_i$.\\
% On choisit, ceci simplifie la suite, $A_0$ parmis les $B_i$.

 Pour tout $i\in I$, $\chi$ induit un homéomorphisme entre $\bar B_i$ et $\hat B_i$. Comme $x\in \bar B_i$, $\chi(x)=\hat x\in \hat B_i$. De plus, $\hat V\cap \hat B_i$ est un voisinage de $\hat x$ dans $\hat B_i$, donc $\chi\inv( \hat V\cap \hat B_i)\subset V$ est un voisinage de $x$ dans $\bar B_i$. Il contient donc un ensemble de la forme $\V'_{B_i}(f_i,R_i)$, avec $R_i\in \U$ et $f_i\in \F_{B_i}$, $[f_i]=x$ ($\V'_{B_i}(f_i,R_i)$ est défini en \ref{defin:Vprime}). Comme $\delta(f)\subset B_i$, on peut supposer que $\delta(f_i)\subset \delta(f)\subset A$ (voir \ref{lemme:scp dans le coeur}). Soit $\gamma$ un scp de $\bigcap_{i\in I} \delta(f_i)$ tel que $B_A(o,\lambda r)\subset -\gamma^*_A$. Soit $R=\bigcap_{i\in I} R_i$. 
Alors $R\in\U$, et $\gamma$ est un coeur de cône affine. Soit $g\in \F_{A_0}$ le cône engendré par $\gamma$, c'est un scp de $f$ donc $[g]=x$ et $\V_{\I}(g,R)$ est un voisinage de $x$ dans $\I$. Nous savons que pour tout $i\in I$, $\V_{\I}(g,R)\cap B_i\subset V$. Montrons que $\V_\I(g,R)\subset Q.V$, cela conclura la preuve.\\

 Soit $t\in \V_\I(g,R)$, il suffit de montrer qu'il existe $h\in Q$ et $i\in I$ tel que $h.t\in \V'_{ B_i}(f_i,R_i)$ .\\
D'après la proposition \ref{prop:voisinages dans I et dans A}, il existe un appartement $Z$ contenant $b(g)$ et un représentant $\tau\in\F_Z$ de $t$, tel que $\tau\subset (\phi_Z(g)+R)\cap\phi_Z(g)^*$, où $\phi_Z:A\isom Z$ est un isomorphisme fixant $b(g)$. Quitte à réduire $\tau$, on peut supposer $\delta(\tau)$ inclus dans un des quartiers de la forme $\barre{b(g)+\vec C}$ avec $\vec C$ une chambre de $\phi_Z(\gamma)^*$. Comme $o$ a été choisi dans $-\delta(f)$, il est dans $-\gamma$, donc aucun mur de $A_0$ dont la direction contient $\vec\gamma$ ne sépare strictement $o$ et $b(g)$. La proposition \ref{prop:cone et chambre dans un appart} permet de conclure à l'existence d'un appartement $Y$ contenant $\{o\}\cup b(g)\cup \delta(\tau)$.\\
 Ensuite, la proposition \ref{prop:coeur et etoile dans un appart} prouve qu'il existe un appartement $X$ contenant $-\gamma^*_Y \cup \gamma$. En particulier, $X\cap A$ contient $o$ et $\gamma$, donc $o+\vec\gamma=o+\delta(\vec f)$, ce qui contient $\delta(f)$. De plus, nous avons supposé $B_A(o,\lambda r)\subset -\gamma^*_A$, ce qui entraîne que $B_Y(o,\lambda r)\subset -\gamma^*_Y\subset Y\cap X$.\\
Soit $i\in I$ tel que $B_i\cap X\supset B_X(o,\lambda r)\cup\delta(f)$.
Par le corollaire \ref{cor:fixe boule}, il existe $h\in G$ tel que $h.Y=B_i$, $h$ fixe $Y\cap B_i$ et $B_\I(o,r)$. En particulier, $h\in Q$.\\
 %Par forte transitivité de $G_1$ sur $\I$, il existe $h\in G_1$ fixant une chambre de $B_Y(o,\lambda r)$ et tel que $h.Y=B_i$. Alors $h$ fixe $Y\cap B_i$, donc en particulier $b(g)$ et $B_Y(o,\lambda r)$. Par la proposition \ref{prop:fixateur dans moufang}, $h$ fixe $B_\I(o,r)$ donc $h\in Q$.\\
 Enfin, si $\phi_Y:A_0\isom Y$ est un isomorphisme fixant $b(g)$, alors $h.\delta(\tau)\subset h.((\phi_Y(g)+R)\cap \phi_Y(g)^*)= (g_i+R)\cap g_i^*$, où $g_i$ est le cône engendré par $\gamma$ dans $B_i$. Comme $g_i\subset f_i$ et $R\subset R_i$, on obtient $h.\delta(\tau)\subset (f_i+R_i)\cap f_i^*$, d'où $h.t\in \V'_{B_i}(f_i,R_i)\subset V$.\cqfd\\

Ce résultat prouve que les compactifications de $\I$ définies par A. Werner dans \cite{werner} sont identiques à celle présentée ici pour les décompositions en cône $\F^J$ décrites en \ref{subsection:exemple}. En effet, nous avons déjà vu que la décomposition $\Sigma(\rho)$ définie dans \cite{werner} était égale à $\F^J$, et nous avons défini la compactification de $A_0$ à partir de $\F^J$ exactement de la même manière que dans \cite{werner}.\\
En particulier, les compactifications de \cite{landvogt} et \cite{guivarch-remy} coïncident avec celle présentée ici pour $\F=\F^\vide$ la décomposition en cônes de Weyl vectoriels.\\

%section:description de \ib
\section{Description de $\ib$}
\label{section:bord}

Dans cette partie, on vérifie principalement une propriété attendue de $\ib$, à savoir que son bord est réunion d'immeubles affines, dont les groupes de Coxeter sont des sous-groupes de Coxeter de $W$.\\

\subsection{$\ib$ est une réunion d'immeubles}

\begin{lemme} \label{lemme:isom de facade} Soient $A,B\in\A$, soient $x\in \bar A\cap \bar B$, $y\in \bar A$. Si $y$ est dans la même facette de $\bar A$ que $x$ alors $y\in \bar B\cap \bar A$ et $x$ et $y$ sont dans la même facette de $\bar B$.\\
 Notons $a$ cette facette, elle est donc incluse dans $\bar A\cap \bar B$. Il existe un isomorphisme $\phi:A\isom B$ dont le prolongement à $\bar A$ dans $\bar B$ envoie la façade de $A$ contenant $x$ sur la façade de $B$ contenant $x$ et fixe $a$.\\
Enfin, si $f\in\F_A$, $f'\in\F_B$ sont tels que $x=[f]_A=[f']_B$, alors $\vec \phi(\vec f)=\vec f'$.

\end{lemme}

\pv Si $x$ et $y$ dans une même facette de $\bar A$, cela signifie qu'il existe $f,g\in\F_A$ tels que $x=[\tilde f]$ et $y=[\tilde g]$, avec $\vec f=\vec g$ et tels qu'aucun mur dont la direction contient $\vec f$ ne sépare au sens large $f$ de $g$. C'est-à-dire qu'aucun mur dont la direction contient $\vec f$ ne sépare $f$ de $g$ ni ne contient l'un sans contenir l'autre.\\
 D'autre part, $x\in \bar B$ donc il existe $f'\in\F_B$ tel que $x=[\tilde f']$.\\

On choisit un appartement $Z$ contenant, quitte à les raccourcir, $\delta(f)$ et $\delta(f')$. Nous allons montrer que $y\in \bar A\cap \bar Z$, c'est-à-dire qu'il existe un scp de $g$ dont le coeur est inclus dans $A\cap Z$. Il suffit de montrer que $g\cap \cl(\delta(f))\not=\vide$. (Pour une compactification polygonale classique, où $\F=\F^\vide$, c'est automatique.)\\
 On décompose $\vec A$ en somme orthogonale de trois sous-espaces comme dans \ref{prop:cx cx de facade 2}: $\vec A=\vect(\delta(f))\oplus \vec E\oplus \vec f^\perp$, où $\vec E$ est le supplémentaire orthogonal de $\vect(\delta(\vec f))$ dans $\vect(\vec f)$.
 %Nous savons que $\vec f^\perp\sim A_{\vec f}$, que $\vec E$ est un autre complexe de Coxeter et que $\stab_{W(\vec A)}(\vec f)=\fix_{W(\vec A)}(\delta(\vec f))$ est le produit des groupes de coxeter de $\vec E$ et de $A_{\vec f}$.
 Notons $p_{\vec E}$ la projection orthogonale sur $\vec E$, alors $p_{\vec E}(\vec f)$ est un cône convexe qui engendre $\vec E$ et qui est stable par son groupe de Coxeter. Comme $\vec E$ est essentiel, ceci implique que $p_{\vec E}(\vec f)=\vec E$.\\
 Il existe donc un point $v\in g$ tel que $p_{\vec E}(v)=p_{\vec E}(s(f))$. En rajoutant à $v$ un élément assez grand de $\delta(\vec f)$, on peut aussi s'assurer que $v\in f^*_A$. Vérifions que $v\in\cl(\delta(f))$. Supposons qu'un mur $M=\alpha\inv(\{0\})$ sépare strictement $v$ et $\delta(f)$. Disons que $\alpha(\delta(f))> 0$ et $\alpha(v)<0$. Alors $\vec\alpha(\delta(\vec f))\geq 0$ et $\vec\alpha(\fl{s(f)v})<0$. Mais comme $\fl{s(f)v}\in\vec f^*$, ceci implique $\vec\alpha(\delta(\vec f))= 0$, c'est-à-dire $\delta(\vec f)\subset \vec M$.  Mais nous savons qu'un tel mur soit contient $\vec E$ soit contient $\vec f^\perp$. Le premier cas est impossible car il implique $\vec f\subset \vec M$ et nous avons supposé qu'aucun mur dont la direction contient $\vec f$ ne sépare $f$ et $g$. Quand au second cas, il implique que $\vec\alpha$ s'annule sur $\vect(\delta(\vec f))\oplus \vec f^\perp$. Mais comme $p_{\vec E}(v)=p_{\vec E}(s(f))$, il est alors impossible que $M$ sépare $v$ et $s(f)$.\\

 Ainsi, $v+\vec f$ est un scp de $g$ et son coeur est inclus dans $\cl(\delta(f)$ donc dans $A\cap Z$. Ceci prouve que $y\in \bar A\cap \bar Z$. En choisissant un isomorphisme $\phi_1: A\isom Z$ qui fixe $\delta(f)$, et donc $\delta(v+\vec f)$, on vérifie que $x$ et $y$ sont dans la même facette de $\bar Z$. Alors comme précédemment, on peut trouver un scp $g'$ de $\phi_1(v+\vec f)$ dont le coeur est dans $\cl(\delta(f'))$ donc dans $Z\cap B$. On obtient alors que $y\in \bar B$ et $x$ et $y$ sont dans la même facette de $\bar B$.\\
 Il ne reste plus qu'à choisir un isomorphisme $\phi_2:Z\isom B$ fixant $\delta(f')$, et donc $\delta(g')$, et à poser $\phi=\phi_2\circ\phi_1$ pour obtenir un isomorphisme de $A$ sur $B$ dont le prolongement à $\bar A$ fixe $x$. Il est immédiat que $\vec\phi(\vec f)=\vec f'=\vec g$, donc $\phi(A_f)=A_g$. La définition de $\phi$ est indépendante de $y$, et les raisonnements que nous venons d'effectuer prouvent que $\phi$ fixe toute la facette de $\bar A$ contenant $x$.\cqfd\\

\rema  Comme $\bar A$ et $\bar B$ sont fermés, $\bar A\cap \bar B$ contient en fait $\bar a$. Et comme $\phi$ est continu, il fixe $\bar a$.\\

 Lorsque $\bar A\cap \bar B=\vide$, mais que $A$ contient un cône $f$, $B$ un cône $f'$ avec $f\parallele f'$, on peut quand même, comme dans la preuve précédente, en passant par un appartement $Z$ contenant un scp de $\delta(f)$ et un scp de $\delta(f')$, prouver qu'il existe un isomorphisme de $A$ sur $B$ qui induit un isomorphisme entre $A_f$ et $B_{f'}$. Ceci prouve le
\begin{lemme}\label{lemme:isom de facade cas intersection vide}
Soient $A$ et $B$ deux appartements, $f\in\F_A$, $f'\in\F_B$ tels que $f\parallele f'$. Alors il existe un isomorphisme $\phi:A\isom B$ induisant un isomorphisme de $A_f$ sur $B_ {f'}$, tel que $\vec \phi(\vec f)=\vec f'$.\end{lemme}
\cqfd\\

 On définit l'ensemble des facettes de $\ib$ comme étant la réunion des ensembles de facettes de chaque appartement compactifié. Le lemme montre que deux facettes sont disjointes ou égales, et que si un appartement contient un point d'une facette, alors il contient toute la facette, et son adhérence.\\

Comme un immeuble, $\ib$ vérifie les assertions suivantes:

 \begin{prop}\label{prop:comme un immeuble}\begin{itemize}
\item $\ib=\bigcup_{A\in\A} \bar A$
\item Pour deux facettes $a$ et $b$, il existe $\bar A$ contenant les deux.
\item Si $\bar A\cap \bar B$ contient deux facettes $a$ et $b$, alors il existe un "isomorphisme d'appartements compactifiés" $\phi:\bar A\isom \bar B$ fixant $a$ et $b$.\end{itemize}
\end{prop}
Mais bien sûr, les $\bar A$ ne sont pas des complexes de Coxeter.\\

\demo\\
Le premier point vient directement de la définition de $\ib$. Le second est conséquence du lemme \ref{lemme:deux points ou facettes} et du fait que lorsqu'un appartement compactifié contient un point d'une facette, alors il contient toute cette facette.\\

 Prouvons le troisième point. En remplaçant $a$ et $b$ un point inclus dans $a$ et un point inclus dans $b$, on se ramène au cas où $a$ et $b$ sont deux points.\\
Soit $f_A\in\F_A$ un représentant de $a$ dans $A$, et $g_B\in\F_B$ un représentant de $b$ dans $B$. Soient $d_a$ une demi-droite incluse dans l'intérieur de $\delta(f_A)$ dans $\aff(\delta(f_A))$, et $d_b$ une demi-droite incluse dans l'intérieur de $\delta(g_B)$ dans $\aff(\delta(g_b))$. En choisissant des galeries, dans $A$ et $B$, le long de ces demi-droites (proposition \ref{prop:galerie sur droite}), puis en choisissant un appartement contenant des sous-galeries de ces galeries, on obtient un appartement $Z$ dont l'intersection avec $A$ contient un scp de $\delta(f_A)$ et au moins une chambre disons $c_A$. L'intersection de $Z$ avec $B$ contient un scp de $\delta(g_B)$ et une chambre $c_B$. Quitte à réduire $f_A$ et $g_B$, on suppose $Z\supset \delta(f_A)\cup\delta(g_B)$.\\
 A présent, soit $\phi_1:B\isom Z$ un isomorphisme fixant $B\cap Z$. Alors $\phi$ est la restriction à $B$ de la rétraction $\rho_{Z,c_B}$. Lorsqu'on étend $\phi_1$ et $\rho_{Z,c_B}$ par continuité à $\bar B$ et $\ib$, on obtient $\rho_{Z,c_B}(a)=a$ et $\rho_{Z,c_B}(b)=b$ car $a\cup b\subset \bar Z$, d'où $\phi_1(a)=a$ et $\phi_1(b)=b$.\\
 De la même manière il existe $\phi_2:\bar Z\isom \bar A$ qui fixe $a\cup b$, alors $\phi_2\circ\phi_1$ convient.\cqfd\\

\begin{defin}
Soit $F\in\F_\I$, on note $\I_F=\ens{[G]}{G\in\I,\; G\parallele F}$. On appellera $\I_F$ la façade de l'immeuble $\I$ de type $F$.
\end{defin}

  Bien sûr, $\I_F$ ne dépend que de la classe de parallélisme de $F$, et même que de la classe de parallélisme de $\delta(F)$. Si $A$ est un appartement contenant $\delta(F)$, si $f\in\F_A$ est le cône engendré par $\delta(F)$, on notera indifféremment $\I_F$, $\I_f$, ou même $\I_{\delta(F)}$.\\

\begin{prop} Pour tout $F\in\F_\I$, $\I_F$ est un immeuble affine. Ses appartements sont les $A_f$ pour tous les appartements $A$ contenant un cône $f$ tel que $\tilde f\parallele F$.\end{prop}

\demo\\
 Chaque $A_{f}$ est un complexe de Coxeter affine, et la réunion de ces complexes de Coxeter est bien $\I_F$.
Si $A$ et $B$ sont deux appartements contenant des cônes $f$ et $g$ tels que $\tilde f\parallele F\parallele \tilde g$, alors $\delta(f)\parallele \delta(g)$, et le lemme \ref{lemme:isom de facade cas intersection vide} prouve l'existence d'un isomorphisme entre $A_f$ et $B_g$. Ainsi, tous les (présumés) appartements de $\I_\F$ sont isomorphes.

 Soient $a$ et $b$ deux facettes de $\I_F$. Par la proposition \ref{prop:comme un immeuble}, il existe un appartement $A$ tel que $a\cup b\subset \bar A$. Alors $a$ et $b$ sont inclus dans des façades de $A$, disons $a\subset A_g$ et $b\subset A_h$, avec $g,h\in\F_A$. On peut choisir $g$ représentant un élément de $a$ et $h$ représentant un élément de $b$. Alors $\tilde g\parallele F\parallele \tilde h$, d'où $A_g=A_h=A_F$ est un appartement de $\I_F$ qui contient $a\cup b$.\\

Enfin, soient $A_F$ et $B_F$ deux (présumés) appartements de $\I_F$, supposons que $A_F\cap B_F$ contienne deux facettes $a$ et $b$, et cherchons un isomorphisme de $A_F$ sur $B_F$ fixant $a\cup b$. Par la proposition \ref{prop:comme un immeuble}, il existe $\phi:\bar A\isom \bar B$ un isomorphisme d'appartements compactifiés fixant $a\cup b$. Alors $\phi$ induit l'isomorphisme cherché entre $A_F$ et $B_F$ .\cqfd\\

 Si $A$ est un appartement contenant un cône $f$ tel que $\tilde f\parallele F$, alors $\vec f$ est unique (car le parallélisme est une relation d'équivalence: si $g$ est un autre cône avec $\tilde g\parallele F$, alors $f\parallele g$). Ceci autorise la
\begin{defin} Si $A$ est un appartement contenant un cône $f$ tel que $\tilde f\parallele F$, on note $A_F:=A_f$.\\
La projection de $A$ sur $A_f$ sera notée $p_{A,f}$ ou $p_{A,F}$, ou juste $p_A$, $p_f$, $p_F$ selon le niveau de précision requis.
\end{defin}
  Donc en résumé, il y a trois notations pour le même objet: $A_F$, $A_f$ et $A_{\vec f}$.

\rema L'immeuble $\I_F$ n'est pas forcément essentiel. Par exemple, si $\vec\delta(F)$ est contenu dans une chambre de Weyl vectorielle, alors $\I_F$ est un seul appartement, sans mur (voir la remarque de \ref{prop:cx cx de facade}).\\

Maintenant qu'on sait que les façades sont des immeubles, on peut montrer d'autres résultats dans la même veine, par exemple:
\begin{prop}
 Soit $F\in\F_\I$, soient $A,B\in \A$ deux appartements contenant des cônes parallèles à $\delta(F)$. Alors il existe un isomorphisme $\phi:\bar A\isom \bar B$ qui fixe $A_F\cap B_F$.
\end{prop}
\demo\\
Soit $x$ une facette de dimension maximale de $A_F\cap B_f$. Il existe $\phi:\bar A\isom \bar B$ fixant $x$ et envoyant $A_F$ sur $B_F$, par le lemme \ref{lemme:isom de facade}. Notons $\phi_F : A_F\isom B_F$ l'isomorphisme d'appartements de $\I_F$ induit par $\phi$. Comme $\phi_F$ fixe une facette de dimension maximale de $A_F\cap B_F$, il fixe nécessairement $A_F\cap B_F$ (proposition \ref{prop:isom fixant l intersection}). D'où le résultat.\cqfd\\

\subsection{Bord d'une façade}

 On vérifie ici que chaque façade $\I_F$ peut être compactifiée tout comme nous avons compactifié $\I$, et que le résultat est homéomorphe à une réunion que nous préciserons, de plusieurs façades de $\ib$.\\
 On fixe à présent un cône $F\in\F_\I$.\\

\subsubsection{Compactification de $\I_F$}
Choisissons $A$ un appartement contenant $\delta(F)$, soit $p=p_{A,F}$ la projection $A\rightarrow A_F$. Soit $f=F\cap A$, alors $A_F=A_f$ est un complexe de Coxeter affine et le groupe vectoriel associé est isomorphe via la projection $p$ à $\fix_{W(\vec A)}(\vec f)$.\\
On pose $\F^F=\ens{\vec p(\vec g)} {\vec g\in \F \text{ et } \vec f \subset \overline{\vec g} }$.\\

\begin{lemme}\label{lemme:cone stable par une face}
Soit $\vec g\in \F$ dont $\vec f$ est une face, c'est-à-dire $\vec f\subset \overline{\vec g}$. Alors $\vec g$ est stable par $\vec f$, c'est-à-dire $\vec g+\vec f\subset \vec g$.
\end{lemme}
\rema En fait, $\vec g+\vec f= \vec g$ car $\vec g$ est ouvert dans $\vect (\vec g)$ et $0\in\overline{\vec f}$.\\
\pv C'est immédiat si on fixe une description de $\vec g$ et de sa face $\vec f$ en termes d'équations et d'inéquations linéaires comme en \ref{subsection:consequences directes}, mais voici une preuve élémentaire.\\
 Soit $x\in\vec g$ et $v\in \vec f$. Comme $\vec g$ est ouvert dans $\vect(\vec g)$, il existe $r>0$ tel que la boule $B(x,r)$ dans $\vect(\vec g)$  soit incluse dans $\vec g$. Ensuite $v\in\overline{\vec g}$ donc il existe $\epsilon\in\vect(\vec g)$, de norme inférieure à $r$ tel que $v+\epsilon\in\vec g$. Alors $x+v=(x-\epsilon)+(v+\epsilon)\in\vec g$.\cqfd\\

\begin{prop}\label{prop:cones de la facade} L'ensemble $\F^F$ vérifie les conditions requises pour définir une compactification de $\I_F$. De plus, l'application $\vec g\mapsto \vec p(\vec g)$ est une bijection entre l'ensemble des cônes de $\F$ ayant $\vec f$ comme face et $\F^F$.\end{prop}

\demo\\
On commence par les points évidents: $\F^F$ est fini, il contient $\{0\}=\vec p(\vec f)$, et il est stable par le groupe de Weyl $W(\vec A_F)\simeq \fix_{W(\vec A)}(\vec f)$.\\

Montrons que $A_f=\bigcup_{\vec g\in \F^F}\vec g$. Comme le terme de droite est stable par homothéties, il suffit de montrer qu'il contient un voisinage de $0$. Par conséquent, il suffit de montrer que, dans $A$, $\bigcup\{\vec g\;|\; \vec f\subset\overline{\vec g} \}$ contient un voisinage d'un point de $\vec f$. Soit $x\in \vec f$. Si $\vec g\in \F$ est un cône ne contenant pas $\vec f$ dans son adhérence, alors $x\not\in \overline{\vec g}$ (car $\partial \vec g$ est une réunion disjointe de cônes dans $\F$). Il existe donc $U_{\vec g}$ un voisinage de $x$ dans $A$ qui ne coupe pas $\vec g$. Comme $\F$ est fini, le nombre des $U_{\vec g}$ ainsi définis pour $\vec g$ variant parmi les cônes dont $\vec f$ n'est pas une face est fini, et leur intersection est un voisinage de $x$ inclus dans $\bigcup_{\vec g\tq \vec f\subset\overline{\vec g} }\vec g$. Nous avons prouvé que $A_f=\bigcup_{\vec g\in \F^F}\vec g$.\\

Montrons que cette union est disjointe. Soit $\vec p(\vec g)$ et $\vec p(\vec h)$ deux éléments de $\F^F$, où $\vec g,\vec h\in\F$ sont des cônes ayant $\vec f$ comme face. Supposons que $\vec p(\vec g)\cap\vec p(\vec h)$ contient un point, écrivons ce point sous la forme $\vec p(x)$ avec $x\in\vec A$. Alors il existe $v,w\in\vect(\vec f)$ tels que $x+v\in\vec g$ et $x+w\in\vec h$. Mais il existe $v_1,v_2,w_1,w_2\in\vec f$ tels que $v=v_1-v_2$ et $w=w_1-w_2$. Alors le point $x+v_1+w_1=x+v+v_2+w_1=x+w+w_2+v_1$ est dans $\vec g\cap\vec h$, d'après le lemme. Ceci entraîne que $\vec g=\vec h$, donc $\vec p(\vec g)=\vec p(\vec h)$.\\
Ceci montre également que l'application $\vec g\mapsto \vec p(\vec g)$ est bijective.\\
Le petit raisonnement qu'on vient de faire fonctionne aussi dans le cas de cônes affines, il donne le résultat suivant, qu'on note pour utilisation ultérieure:
\begin{lemme}\label{lemme:intersection des releves}
Soit $h,g\in\F_A$ deux cônes affines tels que $\vec f$ borde $\vec g$ et $\vec h$. Si $p(h)\cap p(g)\not=\vide$, alors $h\cap g\not=\vide$.\\
\end{lemme}\cqfd\\

 Passons à la description des cônes à l'aide d'un système d'équations et d'inéquations. Soit $\vec g\in\F$ ayant $\vec f$ comme face. Soient $\{\alpha_i\}_{i\in I\sqcup J\sqcup K}\in \vec A^*$ telles que
$\begin{cases} 	\vec g=\{ \alpha_i>0,\; \alpha_k=0,\: i\in I\sqcup J,\; k\in K\} \\
				\vec f=\{\alpha_i>0,\; \alpha_k=0,\: i\in I,\; k\in J\sqcup K\}
\end{cases}$\\
 Les $\alpha_k$ pour $k\in J\sqcup K$ définissent des formes linéaires sur $\vec A_f$. Vérifions que $\vec p(\vec g)=\{\alpha_j>0,\; \alpha_k=0,\: j\in J,\; k\in K\}$. L'inclusion "$\subset$" est claire. Réciproquement, un point de l'ensemble de droite est de la forme $\vec p(x)$ avec $x\in\vec A$, $\alpha_j(x)>0$ pour $j\in J$ et $\alpha_k(x)=0$ pour $k\in K$, et nous voulons montrer qu'il existe $v\in\vect(\vec f)$ tel que $x+v\in\vec g$. On s'aperçoit qu'il suffit de choisir $v\in\vec f$ assez grand.\\
Nous avons donc obtenu une description de $\vec p(\vec g)$ en système d'équations et d'inéquations linéaires comme requis.\\

Passons aux deux conditions portant sur les faces d'un cône.\\
 Soit $\vec p(\vec g)\in\F^F$, vérifions que son bord est la réunion d'autres cônes de $\F^F$. Soit $\vec p(x)\in\partial \vec p(\vec g)$, soit $\vec p(\vec h)\in\F^F$ le cône contenant $\vec p(x)$. Nous allons montrer que $\vec h\subset\partial \vec g$. On choisit $x$ pour que $x\in\vec h$. Pour tout $n\in\N$, il existe un point de $\vec p(\vec g)$ à distance inférieure à $\frac{1}{n}$ de $\vec p(x)$. Cela signifie qu'il existe $\epsilon_n\in\vec A$, $v_n\in\vect(\vec f)$ tels que $x+v_n+\epsilon_n\in\vec g$ et $||\epsilon_n||<\frac{1}{n}$. On décompose $v_n=v_n^1-v_n^2$ avec $v_n^1,v_n^2\in\vec f$, on impose en outre que $||x+v_n^1||\geq 1$. Alors $x+v_n^1\in\vec h$ et $x+v_n^1+\epsilon_n\in\vec g$.\\
 La suite $(\frac{ x+v_n^1} {||x+v_n^1||} )_n$ a une valeur d'adhérence $y$, a priori dans $\overline{\vec h}$. Mais en fait, cette suite reste dans $x+\vec f$, et $\overline{x+\vec f}=x+\overline{\vec f}\subset \vec h$ par le lemme \ref{lemme:cone stable par une face} appliqué à $\vec f$ et à ses faces. Donc $y\in\vec h$. La suite $(\frac{ x+v_n^1+\epsilon_n }{ ||x+v_n^1|| })_n$ aussi admet $y$ comme valeur d'adhérence, car $||x+v_n^1||\geq 1$ et $\epsilon_n\rightarrow 0$. Donc $y\in\overline{\vec g}\cap \vec h$. Comme $\vec h\not=\vec g$, $\vec h$ est une face de $\vec g$: $\vec h\subset \partial \vec g$. D'où $\vec p(\vec h)\subset \vec p(\partial \vec g)\subset \overline{\vec p(\vec g)}$. Mais comme $\vec p(\vec h)\cap \vec p(\vec g)=\vide$, on arrive à $\vec p(\vec h)\subset \partial \vec p(\vec g)$.\\
Ainsi, $\partial \vec p(\vec g)$ est une réunion d'autres cônes de $\F^F$.\\
 Enfin, il reste à montrer que $\overline{\vec p(\vec h)}=\vect(\vec p( \vec h))\cap \overline{\vec p(\vec g)}$, et nous savons que $\overline{\vec h}=\vect(h)\cap\overline{\vec g}$. Par conséquent, il suffit de montrer que $\overline{\vec p(\vec h)}=\vec p(\overline {\vec h})$ et $\overline{\vec p(\vec g)}=\vec p(\overline {\vec g})$. Les deux égalités, pour $\vec h$ et $\vec g$ sont similaires, on ne traite que celle concernant $\vec g$. L'inclusion $\supset$ vient de la continuité de $\vec p$. Pour l'autre inclusion, il suffit de montrer que $\vec p(\overline {\vec g})$ est fermé, ce qui revient à montrer que $\overline{\vec g}+\vect(\vec f)$ est fermé. On reprend les $\{\alpha_i\}_{i\in I\sqcup J\sqcup K}$ comme au paragraphe précédent, on vérifie alors que $\overline{\vec g}+\vect(\vec f)=\{\alpha_j\geq,\; \alpha_k=0, \: j\in J,\; k\in K\}$, c'est bien fermé.\cqfd\\

 L'avant dernier paragraphe de la preuve prouvait en fait:
\begin{lemme} \label{lemme:projection des faces} Si $\vec g,\vec h\in\F$ sont bordés par $\vec f$ et si $\vec p(\vec h)$ est une face de $\vec p(\vec g)$, alors $\vec h$ est une face de $\vec g$.\end{lemme}

Nous pouvons donc compactifier $\I_F$ à l'aide de la décomposition $\F^F$ en cônes de l'appartement $A_F$. On note provisoirement $\widehat \I_F$ l'espace obtenu, si $B_F$ est un appartement de $\I_F$, on note $\hat B_F$ sa compactification.\\
Vérifions rapidement que $\chap{\I_F}$ ne dépend pas du choix de l'appartement $A$:

\begin{prop} Soit $B_F$ un appartement de $\I_F$, soit $g$ le cône de $B$ engendré par un coeur parallèle à $\delta(F)$, soit $p_B=p_{B,F}:B\rightarrow B_F$ la projection. Alors l'ensemble des cônes vectoriels de $B_F$ est $\ens{\vec p(\vec h)}{ \vec h\in \F\et \vec g\subset \barre{\vec h}}$.\end{prop}

\demo\\
Par définition, ceci est vrai pour $B=A$. Pour $B$ quelconque, l'ensemble des cônes vectoriels de $\vec B_F$ est $\ens{\vec \phi(\vec k)} {\vec k\in\F^F}$ où $\phi:A_F\isom B_F$ est un isomorphisme de complexes de Coxeter dont le choix n'importe pas.\\
 Comme $f\parallele g$, il existe $\xi:A\isom B$ tel que $\vec\xi(\vec f)=\vec g$ (lemme \ref{lemme:isom de facade cas intersection vide}). Alors $\xi$ induit un isomorphisme entre $A_f$ et $B_g$, on peut choisir $\phi$ comme étant égal à cet isomorphisme. On a alors $\phi\circ p_A=p_B\circ \xi$, et on vérifie que $\ens{\vec \phi(\vec k)} {\vec k\in\F^F}=\ens{\vec p_B(\vec h)} {\vec g\subset \barre{\vec h}}$.\cqfd\\

Pour finir, étudions le lien entre le coeur d'un cône $\vec g\in\F$ bordé par $\vec f$ et le coeur de $\vec p(\vec g)\in\F^F$:
\begin{prop}\label{prop:projection du coeur}
Soit $\vec g\in\F$ un cône bordé par $\vec f$, soit $B_F$ un appartement de $\I_F$. Alors \begin{itemize}
\item $\vec p_B(\delta(\vec g))\subset \delta(\vec p_B(\vec g))$.
\item $\vec p_B( \vec g^*)=\vec p_B(\vec g)^*$

\item Si $g\in\F_B$ est de direction $\vec g$, alors $p_B(b(g))\subset b(p_B(g))$.
\end{itemize}
\end{prop}

\demo\\
Notons $p=p_B$ et fixons $f\in\F_B$ tel que $\tilde f\parallele f$.\\
 Soit $x\in\delta(\vec g)$, montrons que $\vec p(x)\in\delta(\vec p(\vec g))$. Il faut donc montrer que $x$ est fixé par $\stab_{W(\vec B_F)}(\vec p(\vec g))$. Soit $w\in \stab_{W(\vec B_F)}(\vec p_B(\vec g))$, on identifie $w$ à un élément de $\fix_{W(\vec B)}(\vec f)$, ce qui assure déjà que $w(\vec g)$ est encore bordé par $\vec f$. Et comme $\vec p(\vec g)=w(\vec p(\vec g))=\vec p(w(\vec g))$, on obtient par la deuxième partie de la proposition \ref{prop:cones de la facade} que $w(\vec g)=\vec g$. Alors $w(x)=x$ par définition de $\delta(\vec g)$, et donc $w(\vec p(x))=\vec p(x)$.\\

Passons au second point. Si $\vec C$ est une facette de Weyl dans $\vec B$ dont l'adhérence contient $\delta(\vec g)$, alors $\vec p(\delta(\vec g))\subset \vec p(\barre {\vec C})\subset \barre{\vec p(\vec C)}$. Notons $\vec C_F$ la facette de Weyl de $\vec B_F$ qui contient $\vec p(\vec C)$, nous venons de prouver que son adhérence contient $\vec p(\delta(\vec g))$, elle contient donc toute la facette de Weyl contenant $\vec p(\delta(\vec g))$, et cette facette contient $\delta(\vec p(\vec g))$ par le premier point. Donc $\vec C_F\subset \vec p(\vec g)^*$ et $\vec p(\vec C)\subset \vec p(\vec g)^*$. Ce qui prouve $\vec p(\vec g^*)\subset \vec p(\vec g)^*$.\\
Montrons l'autre inclusion: soit $\vec C_F$ une facette de Weyl de $\vec B_F$ dont l'adhérence contient $\delta(\vec p(\vec g))$. Soient $\{\alpha_i\}_{i\in I\sqcup J}$ des racines de $\vec B_F$, identifiées à des racines de $\vec B$ s'annulant sur $f$, telles que $\vec C_F=\ens{x\in\vec B_F} {\alpha_i(x)=0\et \alpha_j(x)>0,\; \forall i\in I,\; j\in J}$. Alors $\barre{\vec p\inv(\vec C_F)}$ contient $\delta(\vec g)$ et  $\vec p\inv(\vec C_F)=\ens{x\in\vec B} {\alpha_i(x)=0\et \alpha_j(x)>0,\; \forall i\in I,\; j\in J}$.
 Il existe une facette $\vec C$ de $\vec B$ incluse dans $\vec p\inv(\vec C_F)$, de dimension maximale, et dont l'adhérence contient $\delta(\vec g)\cup\delta(\vec f)$ (proposition \ref{prop:coeur vect}, point 9). Alors $\barre{\vec C}$ contient un cône $\vec f_1$ inclus dans $\vec f$, tel que $\vect(\vec f_1)=\vect(\vec f)$.  Il existe des racines $\{\alpha_k\}_{k\in K}$ telles que $\vec C= \ens{x\in \vec B} {\alpha_i(x)=0\et \alpha_j(x)>0,\; \forall i\in I,\; j\in J\sqcup  K}$. On ne rajoute aucune condition du type $\alpha_k=0$ car $\vec C$ est choisi de dimension maximale. S'il existe $k\in K$ tel que $\alpha_k(\vec f)=0$, alors $\ker(\alpha_k)$ induit un mur de $\vec B_F$, qui ne peut couper $\vec C_F$. Donc $\alpha_k(\vec C_F)>0$, puis $\alpha_k(\vec p\inv(\vec C_F))>0$, ce qui signifie que la condition $\alpha_k>0$ est inutile pour définir $\vec C$. On peut donc supposer que pour tout $k\in K$, $\alpha_k(\vec f)\not=\{0\}$.
 Alors, pour tout $k\in K$, il existe $v_k\in\vec f_1$ tel que $\alpha_k(v_k)\not=0$. Comme $v_k\in\vec f_1\subset \barre{\vec C}$, on a automatiquement $\alpha_k(v_k)>0$ et $\alpha_l(v_k)\geq 0$ pour tout $l\in k$. En sommant tous ces $v_k$, on obtient $w\in\vec f$ tel que pour tout $k\in K$, $\alpha_k(w)>0$. Alors pour tout $x\in \vec p\inv(\vec C_F)$, il existe $\lambda\in\R^+$ tel que $x+\lambda w\in \vec C$, et ceci prouve que $\vec p(\vec C)=\vec C_F$. Comme $\delta(\vec g)\subset \barre{\vec C}$, on a $\vec C\subset \vec g^*$, donc $C_F\subset \vec p(\vec g^*)$.\\

 Enfin, pour le dernier point, il faut montrer que $p(b(g))$ est dans chaque demi-appartement de $B_F$ contenant un voisinage de $s(p(g))$ dans $p(g)$. Notons $\D_{B_F}(M,+)$ un tel demi-appartement, $M$ est un mur de $B_F$, identifié à un mur de $B$ dont la direction contient $\vec f$. Soit $O$ ouvert de $B_F$ tel que $O\cap g\subset \D_{B_F}(M,+)$, alors $p\inv(O)\cap g \subset p\inv(O)\cap p\inv(p(g))= p\inv(O\cap g)\subset p\inv(D_{B_F}(M,+))$. Ainsi $p\inv(D_{B_F}(M,+))$ est un demi-appartement de $B$ contenant le voisinage $p\inv(O)\cap g$ de $s(g)$ dans $g$. Ceci entraîne que $b(g)\subset p\inv(D_{B_F}(M,+))$, donc que $p(b(g))\subset D_{B_F}(M,+)$.\cqfd\\

 L'inclusion réciproque du premier point n'est en général pas vraie, mais c'est presque tout comme, car le fait que $\vec p(\delta(\vec  g))\subset \delta(\vec p(\vec g))$ implique que  $\vec p(\delta(\vec g))$ est inclus dans la même facette de Weyl que celle contenant $\delta(\vec p(\vec g))$. En terme de cônes affines, on obtient le
\begin{cor}\label{cor:projection du coeur} Soit $B_F$ un appartement de $\I_F$, soit $g\in\F_B$ un cône bordé par un cône parallèle à $f$. Alors $\cl( p_B(\delta(g))) =\cl(\delta(p_B(g)))$. \end{cor}

On ne peut pas dire directement la même chose pour les bases, car $b(g)$ est définie comme l'intersection de $g$ avec des demi-appartements fermés. Il se pourrait donc a priori que $p(b(g))$ soit inclus dans un bord de la facette de Weyl contenant $b(p(g))$. Mais ce n'est pas le cas. En effet, si $M$ est un mur contenant $p(b(g))$, alors le mur correspondant contient $b(g)$ et donc $g$, d'où finalement $b(p(g))\subset M$. On peut donc énoncer le
\begin{cor}\label{cor:projection de la base}
Soit $B_F$ un appartement de $\I_F$, soit $g\in \F_B$ bordé par un cône parallèle à $f$. Alors $\cl(p_B(b(g))) = \cl(b(p_B(g)))$.
\end{cor}

\subsubsection{$\chap{\I_F}$ comme réunion de façades de $\ib$}

 On montre ici que $\chap{\I_F}$ est homéomorphe à la réunion de $\I_F$ et d'autres façades de $\ib$.\\
On pose:
 \[\B(F)=\ens{G\in\F_\I}{ \exists A\in\A,\; F'\in\F_\I,\; G'\text{ un scp de }G \tq F'\parallele F, \; \delta(F')\cup \delta(G')\subset A\et A\cap F'\subset \barre{A\cap G'}}\]
On dira que $\B(F)$ est l'ensemble des cônes de $\I$ bordés par $F$.\\

\begin{prop} L'espace topologique $\chap{\I_F}$ est homéomorphe à la réunion des $\I_G$ pour $G\in\B(F)$.\\
L'homéomorphisme est donné par: 

\fonc{\chi}{ \chap{\I_F} }{ \bigcup_{G\in\B(F)} \I_G }{\left[ g \right]_{A_F} }{ \left[ g_a \right]_{A} } \\

 Un point $[g]_{A_F}$ d'un appartement compactifié $\bar A_F$ est envoyé sur le point de $\bar A$ égal à $[g_a]$, où $g_a$ est un relevé de $g$ dans $A$, c'est à dire $g_a\in \F_A$, $p_A(g_a)=g$ et $\vec f \subset \barre{\vec g_a}$.\\
De plus, $\chi$ fixe $\I_F$.\\
\end{prop}

\demo\\
\textit{$\chi$ est bien défini:}\\
 Pour commencer, soit $A_F$ un appartement de $\I_F$, soient $g$ et $h$ deux cônes équivalents dans $A_F$, soient $g_a, h_a\in\F_A$ des relevés, montrons que $g_a\sim_A h_a$. On a égalité des directions: $\vec p(\vec g_a)=\vec g=\vec h=\vec p(\vec h_a)$, donc $\vec g_a=\vec h_a$ par la deuxième assertion de la proposition \ref{prop:cones de la facade}. Et $g\cap h\not=\vide$ donc le lemme \ref{lemme:intersection des releves} s'applique: $g_a\cap h_a\not= \vide$. Donc $g_a\sim h_a$.\\

 Passons au cas général: soit $x\in\chap{\I_F}$, soient $A_F$ et $B_F$ deux appartements de $\I_F$, soit $g\in\F_{A_F}$ et $h\in\F_{B_F}$ tels que $x=[g]_{A_F}=[h]_{B_F}$. Soient $g_a\in \F_A$ et $h_b\in\F_B$ des relevés. Alors $\{\tilde g_a,\tilde h_b\}\subset \B(F)$

 Il existe $f_a\in\F_A$ et $f_b\in\F_B$ tels que $\delta(f_a)\parallele \delta(F)\parallele \delta(f_b)$, $f_a\subset \bar g_a$, $f_b\subset \bar h_b$, $s(f_a)=s(g_a)$, et $s(f_b)=s(h_b)$. D'après la proposition \ref{prop:coeur vect}, point 9, il existe une facette de quartier $k\subset A$ et une autre $l\subset B$ telles que $\delta(g_a)\cup \delta(f_a)\subset \bar k$ et $\delta(h_b)\cup\delta(f_b)\subset \bar l$. On choisit $k$ et $l$ minimales. Soit $Z\in \A$ contenant un scp de $k$ et un scp de $l$, montrons que $Z$ contient alors le coeur d'un scp de $g_a$ et le coeur d'un scp de $h_b$.\\
Il y a deux possibilités: soit $\delta(\vec f_a)\subset \barre{\delta(\vec g_a)}$ et alors $\vec k$ est la facette de Weyl contenant $\delta(\vec g_a)$ et donc $\delta(\vec g_a)\subset \vec g_a\cap \vec k$, soit $\delta(\vec f_a)$ est disjoint de $\barre{\delta(\vec g_a)}$, alors $\overline{\vec k}$ est la plus petite facette fermée contenant $\conv(\delta(\vec f_a)\cup \delta(\vec g_a))$, alors $]u,v[\subset \vec g_a\cap \vec k$ pour n'importe quel $u\in \delta(\vec g_a)$ et $v\in \delta(\vec f_a)$. Dans tous les cas, $\vec k\cap \vec g_a\not=\vide$. Donc tout scp de $k$ contient un point $t$ de $g_a$, alors $t+\delta(\vec g)\subset k$ par le lemme \ref{lemme:cone stable par une face}. Ceci prouve que $A\cap Z$ contient le coeur d'un scp de $g_a$. De même, $B\cap Z$ contient le coeur d'un scp de $h_b$. Notons $g_z$ et $h_z$ les cônes engendrés par ces coeurs, ainsi $\tilde g_z\sim_\I \tilde g_a$ et $\tilde h_z\sim_\I \tilde h_b$.\\

 Comme $Z$ contient également un cône parallèle à $f_a$, il contient un cône parallèle à $f$ et fournit donc un appartement $Z_F$ de $\I_F$. Nous voulons montrer que $\tilde g_a\sim\tilde h_b$, et pour cela il suffit de montrer que $g_z\sim_Z h_z$. En raison du cas particulier traité au début de la preuve, il suffit de prouver que $p_Z(g_z)\sim_{Z_F} p_Z(h_z)$. Ce sera fait si nous prouvons que $\tilde g=\tilda{p_A(g_a)}\sim_{\I_F} \tilda{p_Z(g_z)}$ et $\tilde h=\tilda{p_B(h_b)}\sim_{\I_F} \tilda{p_Z(h_z)}$. Bien sûr les deux égalités sont similaires, vérifions juste la première.\\
 Quitte à remplacer $g_a$ par le scp engendré par $\delta(g_z)$, et $g$ par la projection sur $A_F$ du cône obtenu, on peut supposer $k\subset A\cap Z$, $\delta(g_a)=\delta(g_z)$ et $p_A(g_a)=g$. L'intersection $\bar A\cap \bar Z$ est fermée, donc contient $\bar k$. Or $\bar k$ contient $p_A(\bar k\cap A)$ car $\delta(\vec f_a)\subset \barre{\vec k}$ donc $\bar k\cap A$ est stable par addition avec $\delta(\vec f_a)$. Donc $A_F\cap Z_F\supset \cl (p_A(\delta(g_a)))$. Et, d'après le corollaire \ref{cor:projection du coeur}, $A_F\cap Z_F\supset \delta(p_A(g_a))=\delta(g)$.\\
 On peut maintenant choisir $\phi:\bar A\isom \bar Z$ un "isomorphisme d'appartements compactifiés" qui fixe $A_F\cap Z_F$. Alors $\phi(g)=\phi\circ p_A(g_a)=p_Z\circ\phi(g_a)=p_Z(g_z)$. Ainsi, $p_Z(g_z)$ est l'image de $g$ par un isomorphisme d'appartements qui fixe $\delta(g)$, donc $p_Z(g_z)$ et $g$ ont le même coeur, et donc $\tilde g=\tilda{p_Z(g_z)}$.\\
Ceci achève de prouver que $\chi$ est bien définie.\\

Prenons note de ce résultat, que nous avons montré au passage car il resservira:
\begin{lemme}\label{lemme:tout dans un appart}
Soient $F',G,H\in\F_\I$ tels que $G\in \B(F)$ et $H\in\B(F')$. Alors il existe un appartement $Z$, il existe $g_z,f_z,f'_z,h_z\in\F_Z$ tels que $\tilde g_z$ et $\tilde h_z$ sont des scp de $G$ et $H$, $\tilde f_z$ et $\tilde f'_z$ sont parallèles à $F$ et $F'$, $f_z\subset \bar g_z$ et $f'_z\subset \bar h_z$.\\
 De plus, si $A$ est un appartement de $\I$ contenant $\delta(G)$ et un cône parallèle à $\delta(F)$, si $g_a=A\cap G$, alors $\tilda{p_A(g_a)}\sim_{\I_F} \tilda{p_Z(g_z)}$. Et similairement pour $H$ et $F'$.
\end{lemme}
\rema En appliquant ce lemme au cas où $F'=H$ et $G\sim H$, on vérifie que si $G\in\B(F)$ et $H\sim G$, alors $H\in\B(F)$.\\

\textit{$\chi$ est surjective:} Soit $x\in\bigcup_{G\in\B(F)} \I_G$. Il existe un appartement $A$ et un cône $g\in\F_A$ tel que $x=[g]_A$ et $F$ borde $\tilde g$. Alors $p_A(g)\in \F_{A_F}$ et $x=\chi([p_A(g)]_{A_F})$.\\

\textit{$\chi$ est injective:}
Soient $x,y\in \chap{\I_F}$, $x=[g]_{A_F}$, $y=[h]_{B_F}$, et supposons que $\chi(x)=\chi(y)$. Soient $g_a\in \F_A$, $h_b\in\F_B$ des relevés, alors $\tilde g_a\sim \tilde h_b$. Quitte à raccourcir $g_a$ et $h_b$, on peut supposer grâce au lemme \ref{lemme:tout dans un appart} qu'il existe $Z\in\A$, $g_z,h_z,f_z\in \F_Z$ tels que $\tilde g_z\sim\tilde g_a$, $\tilde h_z\sim \tilde h_b$ et $\tilde f_z\parallele F$. Alors $g_z\sim_Z h_z$, c'est-à-dire que $g_z\cap h_z$ contient un scp de $g_z$ et de $h_z$. Alors l'image par $p_Z$ de ce scp est un scp de $p_Z(g_z)$ et de $p_Z(h_Z)$, donc $p_Z(g_z)\sim_{Z_F} p_Z(h_Z)$. Mais le lemme affirme en outre que $\tilde g=\tilda{p_A(g_a)}\sim_{\I_F} \tilda{p_Z(g_z)}$ et similairement pour $h$. Ceci implique que $\tilde g\sim_{\I_F} \tilde h$ donc $x=y$.\\

\textit{$\chi$ est continue:}
 Soit $x=[g]_{A_F}$ un point d'un appartement compactifié $\bar A_F$ de $\I_F$. Fixons un relevé $g_a\in\F_A$ de $g$ et un ouvert $U\in\U$ de sorte que $V:= \V_\I(\tilde g_a,U) \cap (\cup_{G\in \B(F)}\I_G)$ est un voisinage de $\chi(x)$, et les voisinages obtenus de la sorte forment une base de voisinages de $\chi(x)$. Montrons que $\chi\inv( V)$ contient un voisinage de $x$, précisément, nous allons montrer que $\V_{\I_F}(g,\vec p(U))  \subset \chi\inv(V)$. Soit $f\in\F_A$ un cône tel que $\tilde f\parallele F$.\\
 Soit $t\in \V_{\I_F}(g,\vec p(U))$, soit $H\in\F_\I$ un représentant de $t$. Alors $H$ est bordé par $F$, donc d'après la proposition \ref{prop:parallele}, il existe $B\in \A$ contenant $b(g_a)+\delta(\vec f)$ et un scp de $\delta(H)$, on peut supposer $\delta(H)$ lui-même. L'intersection $\bar A\cap \bar B$ contient alors $p(b(g_a))$ et donc $b(g)$ par le corollaire \ref{cor:projection de la base}. Soit $\phi:\bar B\isom \bar A$ induisant un isomorphisme $B_F\isom A_F$ fixant $b(g)$. Alors $\phi(t)\in \V'_{A_F}(g,\vec p(U))$ (voir 
\ref{prop:voisinages dans I et dans A}). On a $\chi\circ\phi(t)=\phi\circ\chi(t)$, il reste donc à vérifier que $\chi(\phi(t))\in \V'_A(g_a,U)$.\\
 Soit $h\in A_F$ un représentant de $\phi(t)$ inclus dans $(g+\vec p_A(U))\cap g^*$. Comme $p_A(g_a+U)=g+\vec p(U)$, il existe $s\in g_a+U$ tel que $p_A(s)=s(g)$. Mais alors $s\in p\inv(g^*)$, donc d'après la proposition \ref{prop:projection du coeur}, il existe $\vec v\in\vec f$ tel que $s+\vec v\in g_a^*$. Alors $s+\vec v\in (g_a+U)\cap g_a^*$. On choisit le relevé $h_a$ de $h$ dont le sommet est $s+\vec v$, donc $h_a=s+\vec v+\vec h_a\subset (g_a+U)\cap g_a^*$, car $\vec h_a$ est dans le bord de $\vec g_a$ (lemme \ref{lemme:projection des faces}). Ceci prouve que $\chi(\phi(t))\in \V'_A(g_a,U)$, donc $\chi$ est continue.\\

\textit {$\chi$ est fermée} car $\chap{\I_F}$ est compact et $\ib$ est séparé.\\

\textit{$\chi$ fixe $\I_F$: } Soit $x\in\I_F$, soit $A_F$ un appartement contenant $x$, soit $f\in\F_A$ un représentant de $x$, on a $\tilde f\parallele f$. On a $x=[ \{x\}]_{A_F}$ et un relevé de $\{x\}$ dans $A$ est $f$. Donc $\chi(x)=[f]_A=x$.\cqfd\\

Ainsi, $\chi(\chap{\I_F})$ est compact, donc fermé dans $\ib$. De plus, $\I_F$ est dense dans $\chap{\I_F}$, donc $\chi(\I_F)=\I_F$ est dense dans $\chi(\chap{\I_F})$. D'où le 
\begin{cor} L'adhérence de $\I_F$ dans $\ib$ est $\chi(\chap{\I_F})=\bigcup_{G\in\B(F)}\I_G$.\end{cor}

On peut donc identifier, au moyen de $\chi$, la compactification $\chap{\I_F}$ de $\I_F$ à la clôture de $\I_F$ dans $\ib$.

\bibliographystyle{alpha}
\bibliography{compactification}
\vspace{.5cm}

\begin{flushright}
\begin{minipage}{7cm}
Cyril Charignon\\
Institut \'Elie Cartan\\
Unité mixte de recherche 7502\\
Nancy-Université, CNRS, INRIA\\
Boulevard des aiguillettes\\
BP 70239\\
54506 Vandoeuvre lès Nancy cedex (France)\\
cyril.charignon@iecn.u-nancy.fr
\end{minipage}\end{flushright}
   
\end{document}